\documentclass[reqno,11pt]{amsart}
\usepackage{amsmath,amsfonts,amssymb,amsxtra,latexsym,amscd,enumerate,amsthm,verbatim}

\usepackage[dvipsnames]{xcolor}
\usepackage{graphicx}

\usepackage{url}
\usepackage[margin=1.40in]{geometry}
\usepackage{booktabs}
\numberwithin{equation}{section}

\newcommand{\R}{\mathbb{R}}

\newcommand{\T}{\mathbb{T}}
\newcommand{\C}{\mathbb{C}}
\newcommand{\Z}{\mathbb{Z}}
\newcommand{\eps}{\epsilon}
\newcommand{\al}{\alpha}
\newcommand{\be}{\beta}

\newcommand{\wPsi}{\widetilde{\Psi}}
\newcommand{\wOmega}{\widetilde{\Omega}}
\newcommand{\Luv}{{L_{uv}}}

\newcommand{\wt}{\widetilde}
\numberwithin{equation}{section} %pour numeroter les equations par section

\newtheorem{theorem}{Theorem}[section]
\newtheorem{lemma}[theorem]{Lemma}
\newtheorem{proposition}[theorem]{Proposition}
\newtheorem{corollary}[theorem]{Corollary}
\newtheorem{definition}[theorem]{Definition}

\begin{document}

\title[Non-uniqueness of axi-symmetric swirl-free Navier--Stokes equations, II]{On the non-uniqueness of solutions of the axi-symmetric swirl-free Navier--Stokes equations, II}

\author{M. Cadiot}
\address{Ecole Polytechnique}
\email{matthieu.cadiot@polytechnique.edu}
\author{A. D. Ionescu}
\address{Princeton University}
\email{aionescu@math.princeton.edu}
\author{H. Jia}
\address{University of Minnesota--Twin Cities}
\email{jia@umn.edu}
\author{S. Palasek}
\address{Princeton University}
\email{spalasek@princeton.edu}

\begin{abstract}
We prove optimal non-uniqueness of Leray solutions of the unforced incompressible Navier--Stokes equations in 3D within the class of axially symmetric swirl-free flows. More precisely, we construct two solutions $u_1,u_2$ that have the same initial data $u_1(0)=u_2(0)$, are smooth for every positive time, satisfy the energy identity, and satisfy 
\[
u_1,u_2\in C_tH_x^\alpha\cap L_t^2H_x^{\alpha+1}\,\,\text{ for every }\,\,0\leq\alpha<1/2,\qquad u_1,u_2\in L_t^\infty L_x^{3,\infty}\cap L_t^{2,\infty}L_x^\infty.
\]
Thus non-uniqueness holds in optimal Sobolev spaces and critical Ladyzhenskaya--Prodi--Serrin weak spaces, even in a symmetry class where smooth solutions are globally regular.

The proof relies on a computer-assisted construction of a linearly unstable forward self-similar solution, starting from the numerical profiles identified in our preceding paper \cite{IoJiPa}. We combine explicit asymptotic expansions and global residual estimates with weak coercivity inequalities stabilized by finitely many functionals. Approximate solutions of auxiliary linear equations yield quantitative inverse bounds, allowing us to correct both the stationary profile and its unstable eigenpair by contraction arguments. The instability and localization theorem proved in \cite{IoJiPa} then produces the non-unique solutions. 

The approximate profiles, the functionals used in the coercivity estimates and the code for the interval-arithmetic verification are available in \cite{proof_code}.
\end{abstract}

\maketitle

\setcounter{tocdepth}{1}

\tableofcontents

\setcounter{tocdepth}{1}

\newcommand{\bydef}{\stackrel{\mbox{\tiny\textnormal{\raisebox{0ex}[0ex][0ex]{def}}}}{=}}

\setcounter{tocdepth}{1}

\section{Introduction}

We consider solutions $u$ of the incompressible Navier--Stokes equation in $\mathbb{R}^3$,
\begin{equation}\label{Nav1}
\partial_tu+u\cdot\nabla u+\nabla p-\Delta u=0,\qquad \nabla\cdot u=0.
\end{equation}
In 1934, Leray \cite{Leray} pioneered the analysis of \eqref{Nav1}, showing local well-posedness of strong solutions from regular initial data and global existence of weak solutions satisfying the energy inequality in the energy space $L_t^\infty L_x^2\cap L_t^2\dot H_x^1$. Two fundamental problems were raised: whether classical solutions can develop finite-time singularities, and whether Leray weak solutions are unique. By weak--strong uniqueness, a Leray solution with smooth initial data coincides with the regular solution for as long as the latter exists, so non-uniqueness from smooth data requires a finite-time singularity. For rough data, on the other hand, no classical solution exists back to the initial time, so instantaneous non-uniqueness of Leray solutions is not ruled out by the theory.

The purpose of this paper is to construct non-unique Leray--Hopf solutions of \eqref{Nav1} within the class of axially symmetric swirl-free flows. The solutions are smooth at every positive time, satisfy the energy identity, and have initial data just below the critical Sobolev regularity $H^{1/2}$. The mechanism is an instability of a forward self-similar solution emanating from a velocity field homogeneous of degree $-1$. We develop a quantitative argument around the profile discovered numerically in \cite{IoJiPa}, reducing its construction and the verification of its instability to explicit bounds on approximate solutions and finitely many auxiliary functions.

\subsection{The main result}

Our main result is the following.

\begin{theorem}\label{Main3D}
There are two different solutions 
\begin{equation}\label{pr9}
u_1, u_2\in C([0,1]:H^\alpha_x)\cap L^2_tH^{\alpha+1}_x\cap L^\infty_tL^{3,\infty}_x,\qquad \text{ for any }\alpha\in[0,1/2),
\end{equation}
of the incompressible Navier--Stokes equations \eqref{Nav1} in $\R^3$, with the same initial data $u_1(0)=u_2(0)=u_0\in H^\alpha\cap L^{3,\infty}$. These solutions are smooth in $(0,1]\times\R^3$, satisfy $t^{1/2}\|u_j(t)\|_{L^\infty}\lesssim 1$, $j\in\{1,2\}$, and satisfy the energy identity 
\begin{equation}\label{pr10}
\frac{1}{2}\|u_j(t_0)\|_{L^2}^2=\frac{1}{2}\|u_j(t_1)\|_{L^2}^2+\int_{t_0}^{t_1}\|\nabla u_j(s)\|_{L^2}^2\,ds
\end{equation}
for any $t_0\leq t_1\in[0,1]$ and $j\in\{1,2\}$.

In addition, the solutions $u_1, u_2$ are axially symmetric and swirl-free.
\end{theorem}

Theorem~\ref{Main3D} improves significantly on all the earlier three-dimensional non-uniqueness results, which are discussed in Subsection~\ref{ear} below, in two ways:

(1) Both solutions $u_1, u_2$ belong to the space of axially symmetric swirl-free vector fields, a space in which the equations are globally well-posed for sufficiently regular data (see \cite{La, UkYu, GaSv}), so the failure of uniqueness for rough data is already present in a setting where global regularity for smooth data is known. The possibility of carrying out the Jia--{\v S}ver{\'a}k program entirely within this class was discovered numerically in Part~I of this work \cite{IoJiPa}, and the reduction is also used in the subsequent work of Buckmaster--Coiculescu--Li~\cite{BuCoLi}.

(2) Both solutions $u_1,u_2$ belong to the space $C_tH^\alpha_x$ for $\alpha<1/2$, which is optimal in view of Fujita--Kato theory \cite{FujitaKato}, and are smooth in $(0,1]\times\R^3$. In particular, the energy identity \eqref{pr10} holds for any times $t_0\leq t_1\in[0,1]$. In terms of integrability, $u_1, u_2 \in L^\infty_tL^{3,\infty}_x\cap L^{2,\infty}_tL^\infty_x$, which is again optimal in view of the classical Ladyzhenskaya--Prodi--Serrin (LPS) regularity criterion. 

In our earlier paper \cite{IoJiPa}, we identified numerically an unstable axially symmetric swirl-free profile and proved an instability-to-non-uniqueness theorem giving the sharp Sobolev and LPS bounds in Theorem~\ref{Main3D}. By that theorem, a single profile with one positive real eigenvalue suffices, and here we construct the profile and an unstable eigenpair rigorously. We control the corrections using weak coercivity in each parity class and approximate solutions of finitely many auxiliary linear equations. The main difficulty is that the linearized operators are not coercive, which we overcome by the stabilized coercivity method outlined below. Since the equations are posed on the whole space and the energy estimates pair the equation with $\Delta^{-1}f$, the residuals must moreover be bounded in a negative Sobolev norm.

\subsection{Earlier and related results}\label{ear} The question of uniqueness of solutions of the Euler and the Navier--Stokes equations is a classical problem that has been studied extensively in recent years. The first non-uniqueness result was proved by Buckmaster--Vicol \cite{BuVi}, who used convex integration on $\T^3$ to construct non-unique finite-energy weak solutions with a small positive Sobolev exponent. Later on, Cheskidov--Luo \cite{ChLu1,ChLu2} obtained non-uniqueness in $L^p_tL^\infty_x$ for every $p<2$, and in two dimensions in $C_tL^p_x$ for $p<2$, thus approaching particular endpoints of the classical uniqueness criteria, while Coiculescu--Palasek \cite{CoPa} introduced a mechanism based on a discrete energy cascade to construct non-unique solutions from critical $\mathrm{BMO}^{-1}$ data, smooth for all positive times. These ideas were further developed by Cheskidov--Dai--Palasek \cite{ChDaPa}, who used an inverse-cascade construction in all dimensions $d\geq 2$. We emphasize that none of the solutions constructed in any of these works are in the Leray class.

To construct non-unique Leray solutions, Jia--{\v S}ver{\'a}k~\cite{JiSv1,JiSv2} proposed a different mechanism based on the instability of a forward self-similar
solution, and Guillod and {\v S}ver{\'a}k~\cite{GuSv} found numerical evidence
for this mechanism in the form of a symmetry-breaking bifurcation. Albritton--Bru{\'e}--Colombo \cite{AlBrCo} were the first to make this mechanism rigorous, constructing two distinct Leray solutions with zero initial velocity and the same singular external force. For the unforced equations, Hou--Wang--Yang~\cite{HoWaYa} were the first to propose a full computer-assisted proof of non-uniqueness of Leray solutions. Very recently, two analytic proofs of non-uniqueness were proposed by Hou--Qin--Song--Zhang~\cite{HoQiSoZh}, using a cascade of Gaussian swirls, and Buckmaster--Coiculescu--Li~\cite{BuCoLi}, using a discretely self-similar flow with an unstable Floquet mode in the axially symmetric swirl-free class.

We note also the important recent OpenAI result \cite{OpenAI} proving the existence of a finite-time blow-up for the three-dimensional Navier--Stokes equations with zero initial velocity and a smooth external force. This result concerns the \emph{forced equations}, while finite-time blow-up for the unforced equation \eqref{Nav1} remains open.

The non-uniqueness problem has also been investigated extensively for the incompressible Euler equations, and we mention here only a few of the main results obtained by convex integration. De Lellis--Sz\'ekelyhidi
\cite{DeLeSz09,DeLeSz10} introduced this approach and showed that certain
initial data admit infinitely many weak solutions satisfying the local
energy equality. Isett~\cite{Isett} proved the optimal version of Onsager's
conjecture on $\T^3$ by constructing, for each $0<\beta<1/3$, $C^\beta$
solutions that fail to conserve energy. Buckmaster--De Lellis--Sz\'ekelyhidi--Vicol
\cite{BuDeSzVi} constructed solutions of this regularity with any
prescribed energy profiles, while Giri--Radu~\cite{GiRa} used a Newton--Nash iteration to construct nonzero
$C^\beta$, $\beta<1/3$ solutions with compact support in time in 2D.

Our proof of Theorem~\ref{Main3D} belongs to a growing body of computer-assisted proofs in PDE. A common approach is to use a fixed-point argument and rigorous bounds on the linearization to validate a numerical approximation. See, for example, \cite{Na,Pl,BrMcPl} for applications to various elliptic problems and \cite{ZgMi,DaLeMi,CaLeNa} for spectral methods in dissipative equations and on unbounded domains. Applications to the Navier--Stokes equations include steady and time-periodic flows~\cite{HeNaXi,ArGaKo,BeBrLeVe} and the linear instability of plane Poiseuille flow~\cite{WaPlNa}.

In fluid mechanics, computer-assisted methods have been used to construct smooth inviscid SQG solutions~\cite{CaCoGo}, turning waves in the Muskat problem~\cite{GoGr}, and overhanging water waves~\cite{CaHa}. Examples of singularity formation include nearly self-similar blow-up for the Boussinesq and axi-symmetric Euler equations in domains with boundary~\cite{ChenHou2022,ChenHou2023} and imploding solutions for compressible fluids~\cite{BuCaGo}. While these results are based on classical numerical methods, another promising route is to discover the numerical object as a neural network; see for instance \cite{WaLaGoBu,PINN}. Beyond fluid mechanics, we refer to~\cite{DaFi,BuCaGoSh,CaIoPa} for computer-assisted constructions of singular solutions of nonlinear dispersive and parabolic equations. A mechanism close to the one used here appears in~\cite{DaSv}, where instability of a forward self-similar profile gives non-uniqueness after blow-up for the complex Ginzburg--Landau equation. For broader accounts of computer-assisted proofs in PDE, see~\cite{NakaoPlumWatanabe2019,Tucker2011,GoSe,Ru,vandenBergLessard2015}.

\subsection{Self-similar instability and the scalar equations}

To explain the Jia--{\v S}ver{\'a}k reduction, we introduce similarity variables
\[
u(t,x)=t^{-1/2}V(s,y),\qquad s=\log t,\quad y=x/\sqrt t.
\]
After applying the Leray projection $\Pi$, the equation becomes
\[
\partial_sV=\Delta V+\tfrac12y\cdot\nabla V+\tfrac12V
-\Pi(V\cdot\nabla V).
\]
A forward self-similar solution corresponds to a stationary point $V=U$. A positive eigenvalue $\widetilde\lambda$ of the linearized evolution at $U$ produces a perturbation with leading behavior $e^{\widetilde\lambda s}\widetilde U$. This perturbation tends to zero as $s\to-\infty$, that is, as $t\to0^+$. Thus an instability in similarity time is compatible with a vanishing perturbation at the initial time, and the nonlinear unstable evolution and the stationary profile can give different solutions with the same initial data. 

Following \cite{IoJiPa}, we work in the class of axially symmetric swirl-free solutions, in which the problem reduces to the scalar equations below. We write the velocity and vorticity in cylindrical coordinates as
\[
U=r\partial_z\Psi\,e_r-(r\partial_r\Psi+2\Psi)e_z,
\qquad \operatorname{curl}U=r\Omega\,e_\theta.
\]
The divergence constraint is then automatic.  For $k\in\Z_+$ let $C^k_{e,o}=C^k_{e,o}(\R^2)$ (respectively $C^k_{e,e}=C^k_{e,e}(\R^2)$) denote the spaces of $C^k$ functions on $\R^2$ which are even in $r$ and odd in $z$ (respectively even in $r$ and even in $z$). We define the operators
\begin{equation}\label{sf1.4}
L_0:=\partial_r^2+\partial_z^2+\frac{3}{r}\partial_r,\qquad L_1:=L_0+\frac{1}{2}\big(r\partial_r+z\partial_z+3\big),
\end{equation}
acting on $C^2$ functions that are even in $r$. Theorem~\ref{Main3D} is a consequence of the following proposition.

\begin{proposition}\label{PropCentral}
There is a real-valued solution $(\Psi,\Omega)\in (C^2_{e,o})^2$ of the nonlinear system
\begin{equation}\label{sf1}
L_0\Psi-\Omega=0,\qquad L_1\Omega+r\partial_r\Psi\partial_z\Omega-r\partial_z\Psi\partial_r\Omega+2\Psi\partial_z\Omega=0,
\end{equation}
which is linearly unstable, in the sense that 
\begin{equation}\label{sf1.2}
\begin{split}
&L_0\wPsi-\wOmega=0,\\
&L_1\wOmega+\big\{r\partial_r\Psi\partial_z\wOmega+r\partial_z\Omega\partial_r\wPsi-r\partial_z\Psi\partial_r\wOmega-r\partial_r\Omega\partial_z\wPsi+2\Psi\partial_z\wOmega+2\partial_z\Omega\wPsi\big\}-\widetilde{\lambda} \wOmega=0,
\end{split}
\end{equation}
for some non-trivial real eigenpair $(\widetilde{\lambda},\wPsi,\wOmega)\in [0.23,0.24]\times (C^2_{e,e})^2$. Moreover, the solutions $(\Psi,\Omega), (\wPsi, \wOmega)$ satisfy the asymptotic bounds
\begin{equation}\label{sf2}
|\wPsi(r,z)|+\Big|\Psi(r,z)-\frac{\sigma z}{r^2+z^2}\Big|\lesssim \langle r,z\rangle^{-2},\qquad|\wOmega(r,z)|+\Big|\Omega(r,z)+\frac{6\sigma z}{(r^2+z^2)^2}\Big|\lesssim  \langle r,z\rangle^{-4},
\end{equation}
for $(r,z)\in\R^2$ with $|r,z|\geq 1$, where $\sigma=130$.
\end{proposition}

The leading potential $\sigma z/(r^2+z^2)$ corresponds to the explicit homogeneous velocity
\[
U_0=\frac{\sigma}{(r^2+z^2)^2}
\big[r(r^2-z^2)e_r-2z^3e_z\big].
\]
We emphasize that the base scalar profiles are odd in $z$, whereas the unstable scalar profiles are even in $z$, so the instability breaks the reflection symmetry of the base velocity. We work with the single value $\sigma=130$ and look for the positive eigenvalue near $\widetilde\lambda_{app}=0.235059685902438$.

\subsection{Outline of the quantitative proof}

As explained above, the main task is to prove Proposition~\ref{PropCentral}. The proof has three parts: (1) constructing numerical approximations $(\Psi_{app},\Omega_{app})$ and $(\widetilde\lambda_{app},\wPsi_{app},\wOmega_{app})$ with small residuals in \eqref{sf1}--\eqref{sf1.2}, (2) upgrading these approximate solutions to exact solutions via a fixed-point argument, and (3) certifying, by interval arithmetic, the quantitative bounds used in the argument.

\subsubsection{Numerical construction of the approximate solutions}
Proposition~\ref{PropCentral2} provides global $C^4$ approximate solutions of the nonlinear and eigenvalue problems, with the fixed positive eigenvalue \eqref{lambda:app}. Starting from the numerical profile found in \cite{IoJiPa}, we refine the interior solutions by collocation in double-double arithmetic. In the coordinates $u=(1+12r^2)^{1/4}-1$, $v=(1+12z^2)^{1/4}-1$, the computational domain is $[0,100]^2$, corresponding to physical distances of approximately $3000$. The computational fields are represented by tensor-product splines of degree $13$, with their parities encoded as in \eqref{transported}.

The exterior fields are the explicit expansions $(\Psi_{\geq-9},\Omega_{\geq-11})$ and $(\wPsi_{\geq-7},\wOmega_{\geq-9})$ of Subsection~\ref{fixed_asymp}. We choose their free coefficients using moment identities for the interior and exterior fields. We glue the expansions to the interior fields between $84$ and $100$, using a polynomial cutoff function of degree nine in each transition variable. The residual calculation includes the derivatives of the cutoff and the infinite exterior. After eliminating the errors in the Poisson equations, the combined residuals satisfy
\begin{equation}\label{ilu0}
\|\mathcal R\|_{Z'}\leq6\times10^{-12},\qquad
\|\wt{\mathcal R}\|_{Z'}\leq1.5\times10^{-12}.
\end{equation}
Table~\ref{close:profiles} collects these bounds and the profile estimates used in the fixed-point argument.

\subsubsection{Weak coercivity and exact solutions}
We describe the correction argument for the nonlinear problem, on real axi-symmetric functions odd in $z$. Lifting the equations to $\R^5=\R^4_x\times\R_z$ identifies $L_0$ with $\Delta_{\R^5}$. Let $\Omega^\ast$ be the lifted approximate vorticity and set $\Psi^\ast=\Delta^{-1}_{\R^5}\Omega^\ast$. After simple reductions, the problem is equivalent to solving the scalar equation for the vorticity correction $f$
\begin{equation}\label{ilu1}
\mathcal O_{\R^5}f+\mathcal L_2f+\mathcal L_3(\Delta^{-1}_{\R^5}f)
+\mathcal N_5(\Delta^{-1}_{\R^5}f,f)=\mathcal R.
\end{equation}
Here $\mathcal O_{\R^5}=\Delta_{\R^5}+\tfrac12S+\tfrac32$, with $S=x\cdot\nabla_x+z\partial_z$, is the Ornstein--Uhlenbeck operator and the remaining operators are defined in \eqref{cl3}. The linearized operator is $\mathcal S_0=\mathcal O_{\R^5}+K$, where $K:=\mathcal L_2+\mathcal L_3\Delta^{-1}_{\R^5}$ is a compact (but large) perturbation. The source $\mathcal R$ in \eqref{cl3.3} includes the effect of the Poisson residual on the vorticity equation.

The spaces we use to solve \eqref{ilu1} are adapted to our linear operator. With all norms taken on $\R^5$, they are defined by
\[
\begin{split}
&\|f\|_Z^2=\|f\|_{L^2}^2+\||\nabla|^{-1}f\|_{L^2}^2,\qquad
\|F\|_{Z'}=\||\nabla|^{-1}\langle\nabla\rangle^{-1}F\|_{L^2},\\
&\|f\|_{Z_0}=\|f\|_Z+\delta\|Sf\|_{Z'},\qquad 0<\delta\ll1,
\end{split}
\]
where $|\nabla|=(-\Delta_{\R^5})^{1/2}$ and $\langle\nabla\rangle=(1-\Delta_{\R^5})^{1/2}$. The $Z$ norm controls both the vorticity correction and the gradient of its potential. The source norm is chosen so that $\Delta_{\R^5}:Z\to Z'$ is bounded and $|\langle F,\Delta^{-1}_{\R^5}f\rangle|\leq\|F\|_{Z'}\|f\|_Z$. Its additional smoothing at high frequencies also makes the decaying lower-order perturbation $K:Z\to Z'$ compact. By including the scaling term in $Z_0$, we ensure that $\mathcal S_0:Z_0^{odd}\to(Z')^{odd}$ is a bounded operator.

The first analytic ingredient is that $\mathcal O_{\R^5}:Z_0^{odd}\to(Z')^{odd}$ is an isomorphism, proved by Fourier inversion in Lemma~\ref{OrnUhl}. Its energy identity is
\[
\langle\mathcal O_{\R^5}f,\Delta^{-1}_{\R^5}f\rangle
=\|f\|_{L^2}^2+\tfrac14\||\nabla|^{-1}f\|_{L^2}^2.
\]
Compactness of $K$ then makes $\mathcal S_0$ a Fredholm operator of index zero, and thus injectivity or surjectivity suffices for invertibility.

Strong coercivity fails for the linearized operator $\mathcal{S}_0$. To prove invertibility we use what we call the {\it{stabilized coercivity method}}: we find a finite-rank correction that makes the linearized operator coercive, and then establish approximate solvability for $\mathcal{S}_0$ in the corresponding stabilizing directions. A small-error argument then yields surjectivity of $\mathcal{S}_0$.

More precisely, we prove first a weak coercive bound of the form
\begin{equation}\label{ilu2}
\langle\mathcal S_0f,\Delta^{-1}_{\R^5}f\rangle
+\sum_{j=1}^{n}\ell_j(f)^2\geq\gamma\|f\|_{Z^{odd}}^2,
\qquad \ell_j(f):=\langle q_j,\Delta^{-1}_{\R^5}f\rangle,
\end{equation}
with $q_j\in(Z')^{odd}$, $\gamma=1/10$, and $n=26$. The point is to choose the functions $q_j$ suitably to account for a finite-dimensional defect of positivity. We construct these functionals in Section~\ref{sec:coercivity}, by showing that we can confine the potentially negative part of the quadratic form to bounded cells. On each cell, a weighted Neumann spectral decomposition isolates finitely many modes, and spectral bounds control the complement. 

The second step is approximate solvability in the stabilizing directions. We construct response functions $\Phi_j\in Z_0^{odd}$ satisfying
\begin{equation}\label{ilu3}
\mathcal S_0\Phi_j=q_j+e_j,
\end{equation}
with small errors $e_j\in Z'$, where $q_j$ are the functions defined in \eqref{ilu2}. 

We first invert the finite-rank perturbation $\mathcal S_+=\mathcal S_0+Q$, where $Qf:=\sum_j\ell_j(f)q_j$. Since $\langle Qf,\Delta^{-1}_{\R^5}f\rangle=\sum_j\ell_j(f)^2$, the weak coercive inequality \eqref{ilu2} is a strong coercive inequality for $\mathcal S_+$. Pairing $\mathcal S_+f=F$ with $\Delta^{-1}_{\R^5}f$ gives $\gamma\|f\|_Z\leq\|F\|_{Z'}$, hence injectivity. Meanwhile, the Fredholm alternative gives surjectivity. We notice that this energy estimate controls only the $Z$ norm, since the anti-selfadjoint part of the scaling operator does not contribute to the quadratic form. The remaining term $\delta\|Sf\|_{Z'}$ in the $Z_0$ norm is then bounded using the equation. Finally, to invert $\mathcal{S}_0$ we define $Gf:=\sum_j\ell_j(f)\Phi_j$ and $Ef:=\sum_j\ell_j(f)e_j$ and notice that
\[
\mathcal S_0(I+G)=\mathcal S_++E,\qquad
\|E\|_{Z\to Z'}\leq\sum_j\|q_j\|_{Z'}\|e_j\|_{Z'}\leq5\times10^{-4}.
\]
Since $\|E\|_{Z\to Z'}$ is small compared to $\gamma$, this factorization shows that $\mathcal S_0$ is surjective, and hence an isomorphism by the Fredholm alternative. Furthermore, since $\mathcal{S}_0^{-1}=(I+G)(\mathcal S_++E)^{-1}$, we can use the quantitative bounds on our solutions in Tables~\ref{close:profiles} and~\ref{close:responses} to estimate
\[
\|\mathcal S_0^{-1}F\|_{Z^{odd}_0}\leq164\|F\|_{(Z')^{odd}},\qquad
\|\mathcal N_5(\Delta^{-1}_{\R^5}f,g)\|_{Z'}\leq0.08\|f\|_Z\|g\|_Z,
\]
where $\mathcal{N}_5$ is the bilinear form in \eqref{ilu1}. Using also \eqref{ilu0}, the main equation \eqref{ilu1} can thus be solved by a fixed-point argument in the $Z^{odd}_0$ ball of radius $10^{-9}$.

The even problem follows the same scheme, with an additional scalar for the eigenvalue correction and a normalization fixing the eigenfunction. Finally, the bootstrap regularity and decay argument in Subsection~\ref{finish} allows us to prove sufficient smoothness and pointwise decay, which completes the proof of Proposition~\ref{PropCentral} and Theorem~\ref{Main3D}.

\subsubsection{Certification of the numerical bounds}
The approximate profiles are determined by their spline knots and double-double coefficients, both components of which enter the estimates. We carry out all arithmetic on these data in interval arithmetic, and we bound the truncation errors and the contribution of the exterior region by analytic estimates. Tables~\ref{close:profiles}, \ref{close:responses} and~\ref{certi1}, together with Proposition~\ref{c:prop:CAP}, collect the resulting bounds.

A principal difficulty is cancellation, since the small residuals are sums of terms which are individually much larger. On each spline cell we therefore clear the rational denominators and form the complete residual numerator by polynomial arithmetic, and only then truncate it and restore the reciprocal factors using Taylor models with explicit remainders. Integrating these models gives bounds on the whole cell; see Sections~\ref{sec:CAP} and~\ref{sec:residual}.

In Subsection~\ref{sec:main-residual} we use Hardy inequalities to bound the norm $Z'$ by weighted $L^2$ norms of the residuals, and the Newton kernel to estimate the contribution of the Poisson residual in the gluing region, which is far from the origin. In the exterior region the approximate profiles are polynomials in $r^2/R^2$, $R^{-1}$ and $\log R$, so the integrals reduce to explicit angular and radial moments.

The coercivity proof also uses bounds on the infinite-dimensional complement of the retained modes. In Appendix~\ref{app:eigenfunctions} we enclose the weighted Neumann eigenvalues and their normalized eigenfunctions and identify their indices. The complementary spectral gaps and the finite matrix bounds in Proposition~\ref{c:prop:CAP} then give the coercive estimate on each cell. In Subsection~\ref{sec:responses} we propagate the eigenfunction enclosures into the exact sources $q_j$ and $\widetilde q_j$, and use them to bound the errors $e_j$ in \eqref{ilu3} and their analogues for the even problem in $Z'$.

\subsection{Organization of the paper}

In Section~2 we develop the functional framework, prove the inverse and nonlinear estimates, and reduce the proof to a collection of quantitative estimates. In Section~3 we specify the global approximate profiles and their asymptotic expansions. We develop the weak coercivity estimates in the odd and even parity classes in Section~\ref{sec:coercivity}, where we also collect the required numerical inputs in Proposition~\ref{c:prop:CAP}. Section~\ref{sec:CAP} develops the local polynomial enclosures and integration estimates. In Section~\ref{sec:residual} we apply them to the profiles, the response equations, and the finite coercivity inequalities. In Appendix~\ref{app:eigenfunctions} we explain how to compute the local spectral data, and in Appendix~\ref{app:verification-formulas} we collect the auxiliary formulas used in the verification.

\subsection{Acknowledgements and AI disclosure} M.~C.\ was supported by the ANR project CAPPS (ANR-23-CE40-0004-01) and by the FMJH (ANR-22-EXES-0013). A. D. I. was supported in part by NSF-FRG grant DMS-2245228 and by NSF-UEFISCDI grant DMS-2407694. S. P. was supported by NSF Grant No. DMS-2424441. H. J. was partially supported by NSF DMS-2453270, NSF DMS-2245021 and a grant from the Simons Foundation. The authors acknowledge that this work was performed using the Princeton Research Computing resources, a consortium including the Princeton Institute for Computational Science and Engineering (PICSciE) and Research Computing at Princeton University. The authors would also like to thank Javier G\'{o}mez-Serrano for useful discussions related to this work. 

This paper and the code for the computer-assisted proof were written by the authors and were not generated by AI. However, Anthropic's Claude Fable 5.1 and OpenAI's ChatGPT Astra 6 assisted with writing some of the codes, language editing, checking calculations, and the verification of the rigorous numerics. The authors reviewed all AI-assisted revisions and take full responsibility for the paper.

\section{Proof of Proposition~\ref{PropCentral} and Theorem~\ref{Main3D}} Our construction starts from explicitly specified approximate solutions with quantitative residual and profile bounds.

\begin{proposition}[Global approximate solutions]\label{PropCentral2} Fix the exact value
\begin{equation}\label{lambda:app}
\widetilde\lambda_{app}:=0.235059685902438.
\end{equation}
There are real-valued profiles $(\Psi_{app},\Omega_{app})\in(C^4_{e,o})^2$ satisfying
\begin{equation}\label{sf3}
L_0\Psi_{app}-\Omega_{app}=E_{\Psi},\qquad L_1\Omega_{app}+r\partial_r\Psi_{app}\partial_z\Omega_{app}-r\partial_z\Psi_{app}\partial_r\Omega_{app}+2\Psi_{app}\partial_z\Omega_{app}=E_{\Omega},
\end{equation}
and non-trivial real-valued profiles $(\wPsi_{app},\wOmega_{app})\in(C^4_{e,e})^2$ satisfying
\begin{equation}\label{sf3.2}
\begin{split}
&L_0\wPsi_{app}-\wOmega_{app}=E_{\wPsi},\qquad L_1\wOmega_{app}+\big\{r\partial_r\Psi_{app}\partial_z\wOmega_{app}+r\partial_z\Omega_{app}\partial_r\wPsi_{app}-r\partial_z\Psi_{app}\partial_r\wOmega_{app}\\
&\qquad\qquad\qquad-r\partial_r\Omega_{app}\partial_z\wPsi_{app}+2\Psi_{app}\partial_z\wOmega_{app}+2\partial_z\Omega_{app}\wPsi_{app}\big\}-\widetilde{\lambda}_{app} \wOmega_{app}=E_{\wOmega},
\end{split}
\end{equation}
The profiles and residuals satisfy all the quantitative bounds in Tables~\ref{close:profiles} and~\ref{certi1}, with the lifted functions, corrected potentials, and combined residuals defined in \eqref{cl1}--\eqref{cl3.3}, the norms defined in \eqref{robu1}, and the test fields defined in \eqref{close:test-fields}. In particular,
\begin{equation}\label{sf3:smallness}
\|\mathcal R\|_{Z'}\leq6\times10^{-12},\qquad
\|\wt{\mathcal R}\|_{Z'}\leq1.5\times10^{-12},\qquad
\|\wt\Omega^\ast\|_{L^2(\R^5)}\in[1,1.01].
\end{equation}
Moreover, they satisfy the exact identities
\begin{equation}\label{sf3.3}
\Psi_{app}=\Psi_{\geq -9}, \qquad \Omega_{app}=\Omega_{\geq -11}, \qquad \wPsi_{app}=\wPsi_{\geq -7},\qquad \wOmega_{app}=\wOmega_{\geq -9}
\end{equation}
in the asymptotic region $\max\{|r|,|z|\}\geq 3000$, where $\Psi_{\geq -9},\Omega_{\geq -11}, \wPsi_{\geq -7}, \wOmega_{\geq -9}$ are the asymptotic solutions constructed explicitly in Subsection~\ref{fixed_asymp}, with the normalization described in Subsection~\ref{sec:computational-domain}. 
\end{proposition}

Setting $\Psi=\Psi_{app}+\psi$, $\Omega=\Omega_{app}+\omega$, $\wPsi=\wPsi_{app}+\widetilde{\psi}$, $\wOmega=\wOmega_{app}+\widetilde{\omega}$, and $\widetilde\lambda=\widetilde\lambda_{app}+\delta\lambda$, the perturbations $(\psi,\omega)$ and $(\delta\lambda,\widetilde{\psi},\widetilde{\omega})$ satisfy the systems
\begin{equation}\label{sf4}
L_0(\psi)-\omega=-E_{\Psi},\qquad L_1(\omega)+L_2(\omega)+L_3(\psi)+\mathcal{N}(\psi,\omega)=-E_{\Omega},
\end{equation}
and
\begin{equation}\label{sf4.3}
\begin{split}
&L_0(\widetilde{\psi})-\widetilde{\omega}=-E_{\wPsi},\qquad L_1(\widetilde{\omega})+L_2(\widetilde{\omega})+L_3(\widetilde{\psi})
 -\widetilde{\lambda}_{app}\widetilde{\omega}
 +\widetilde{L}_2(\omega)+\widetilde{L}_3(\psi)\\
&\qquad\qquad
 +\mathcal{N}(\psi,\widetilde{\omega})+\mathcal{N}(\widetilde{\psi},\omega)
 -\delta\lambda\big(\wOmega_{app}+\widetilde{\omega}\big)=-E_{\wOmega},
\end{split}
\end{equation}
where $L_0,L_1$ are defined as in \eqref{sf1.4}, and
\begin{equation}\label{sf5}
\begin{split}
L_2(f)&:=r\partial_r\Psi_{app}\cdot\partial_zf-r\partial_z\Psi_{app}\cdot\partial_rf+2\Psi_{app}\cdot \partial_zf,\\
L_3(f)&:=r\partial_z\Omega_{app}\cdot\partial_rf-r\partial_r\Omega_{app}\cdot\partial_zf+2\partial_z\Omega_{app}\cdot f,\\
\widetilde{L}_2(f)&:=r\partial_r\wPsi_{app}\cdot\partial_zf-r\partial_z\wPsi_{app}\cdot\partial_rf
+2\wPsi_{app}\cdot\partial_zf,\\
\widetilde{L}_3(f)&:=r\partial_z\wOmega_{app}\cdot\partial_rf-r\partial_r\wOmega_{app}\cdot\partial_zf
+2\partial_z\wOmega_{app}\cdot f,\\
\mathcal{N}(f,g)&:=r\partial_rf\cdot\partial_zg-r\partial_zf\cdot\partial_rg+2f\cdot\partial_zg.
\end{split}
\end{equation}
Notice that $L_2(f)=\mathcal{N}(\Psi_{app},f)$, $L_3(f)=\mathcal{N}(f,\Omega_{app})$,
$\widetilde{L}_2(f)=\mathcal{N}(\wPsi_{app},f)$, $\widetilde{L}_3(f)=\mathcal{N}(f,\wOmega_{app})$.

\subsection{Lifting to axi-symmetric functions on $\R^5$}\label{L0L1} To construct solutions of the systems \eqref{sf4}--\eqref{sf4.3} we need to invert the linear operators $L_0$ and $L_1$. It is convenient to lift our functions to the Euclidean space $\R^5=\R^{4+1}$, in order to have access to a Fourier transform. Precisely, for any continuous function $h:\R^2\to\C$ that is even in $r$ we define its extension
\begin{equation}\label{sf8}
E(h):\R^{4+1}\to\C:\qquad E(h)(x_1,x_2,x_3,x_4,z):=h\Big(\sqrt{x_1^2+\ldots+x_4^2},z\Big).
\end{equation}
Conversely, any axi-symmetric function on $\R^5$ can be identified with a function on $\R^2$ that is even in $r$.

The point is that the operator $L_0$ is equivalent to the usual Laplacian, under this identification,
\begin{equation}\label{sf9}
E(L_0(h))=\Delta_{\R^5}(E(h)).
\end{equation}
Moreover, the operator $L_1$ is equivalent to a shift of the classical Ornstein--Uhlenbeck operator,
\begin{equation}\label{sf9.1}
E(L_1(h))=\mathcal{O}_{\R^5}(E(h)),\qquad\mathcal{O}_{\R^5}:=\Delta_{\R^5}+\frac{1}{2}(x_j\partial_j+z\partial_z)+\frac{3}{2}.
\end{equation}
Here $x=(x_1,\ldots,x_4)\in\R^4$ and $(x,z)\in\R^5$. Throughout, $x_j\partial_j:=\sum_{j=1}^{4}x_j\partial_{x_j}$ denotes the dilation
vector-field in the first four variables.  The full dilation vector-field is
$S:=x_\al\partial_\al=x_j\partial_j+z\partial_z$, with the convention $x_5:=z$.

The operators $\Delta_{\R^5}$ and $\mathcal{O}_{\R^5}$ can be inverted using the Fourier transform. In fact
\begin{equation}\label{sf10}
(\Delta^{-1}_{\R^5}g)(y)=\frac{-1}{8\pi^2}\int_{\R^5}g(Y)\frac{1}{|y-Y|^3}\,dY,
\end{equation}
for $g\in\mathcal{S}(\R^5)$ and $y\in\R^5$. The shifted Ornstein--Uhlenbeck operator can also be inverted on suitable spaces, see Lemma~\ref{OrnUhl} below.

Using this identification, the equations \eqref{sf4}--\eqref{sf4.3} can be transferred to $\R^5$. Let
\begin{equation}\label{cl1}
\begin{split}
&f:=E(\omega),\quad \Psi_5:=E(\Psi_{app}),\quad\Omega^\ast:=E(\Omega_{app}),\quad
\mathcal{R}_{\Psi}:=E(E_{\Psi}),\quad\mathcal{R}_{\Omega}:=E(E_{\Omega}),\\
&\wt f:=E(\wt\omega),\quad \wt\Psi_5:=E(\wPsi_{app}),\quad\wt\Omega^\ast:=E(\wOmega_{app}),\quad
\mathcal{R}_{\wPsi}:=E(E_{\wPsi}),\quad\mathcal{R}_{\wOmega}:=E(E_{\wOmega}).
\end{split}
\end{equation}
The first equations in \eqref{sf4} and \eqref{sf4.3} are equivalent to
\begin{equation}\label{cl1.5}
E(\psi)=\Delta^{-1}_{\R^5}\big(f-\mathcal{R}_{\Psi}\big),\qquad
E(\wt\psi)=\Delta^{-1}_{\R^5}\big(\wt f-\mathcal{R}_{\wPsi}\big),
\end{equation}
so the second equations become the two scalar equations
\begin{equation}\label{cl2}
\mathcal{O}_{\R^5}(f)+\mathcal{L}_2(f)+\mathcal{L}_3\big(\Delta^{-1}_{\R^5}f\big)
+\mathcal{N}_5\big(\Delta^{-1}_{\R^5}f,f\big)=\mathcal{R},
\end{equation}
\begin{equation}\label{cl2.5}
\begin{split}
&\mathcal{O}_{\R^5}(\wt f)+\mathcal{L}_2(\wt f)+\mathcal{L}_3\big(\Delta^{-1}_{\R^5}\wt f\big)
-\wt\lambda_{app}\wt f
+\wt{\mathcal{L}}_2(f)+\wt{\mathcal{L}}_3\big(\Delta^{-1}_{\R^5}f\big)\\
&\qquad\qquad
+\mathcal{N}_5\big(\Delta^{-1}_{\R^5}f,\wt f\big)
+\mathcal{N}_5\big(\Delta^{-1}_{\R^5}\wt f,f\big)
-\delta\lambda\big(\wt\Omega^\ast+\wt f\big)=\wt{\mathcal{R}},
\end{split}
\end{equation}
where
\begin{equation}\label{cl3}
\begin{split}
&\mathcal{L}_2:=(x_j\partial_j\Psi^\ast+2\Psi^\ast)\partial_z-\partial_z\Psi^\ast\cdot x_j\partial_j,\qquad \mathcal{L}_3:=\partial_z\Omega^\ast\cdot x_j\partial_j-(x_j\partial_j\Omega^\ast)\partial_z +2\partial_z\Omega^\ast I,\\
&\wt{\mathcal{L}}_2:=(x_j\partial_j\wt\Psi^\ast+2\wt\Psi^\ast)\partial_z-\partial_z\wt\Psi^\ast\cdot x_j\partial_j,\qquad \wt{\mathcal{L}}_3:=\partial_z\wt\Omega^\ast\cdot x_j\partial_j-(x_j\partial_j\wt\Omega^\ast)\partial_z +2\partial_z\wt\Omega^\ast I,\\
&\mathcal{N}_5(g,h):=x_j\partial_jg\cdot \partial_zh-\partial_zg\cdot x_j\partial_jh+2g\cdot\partial_zh,
\end{split}
\end{equation}
\begin{equation}\label{cl3.3}
\begin{split}
&\Psi^\ast:=\Psi_5-\Delta^{-1}_{\R^5}\mathcal{R}_{\Psi}=\Delta^{-1}_{\R^5}\Omega^\ast,\qquad
\wt\Psi^\ast:=\wt\Psi_5-\Delta^{-1}_{\R^5}\mathcal{R}_{\wPsi}=\Delta^{-1}_{\R^5}\wt\Omega^\ast,\\
&\mathcal{R}:=-\mathcal{R}_{\Omega}+\mathcal{L}_3\big(\Delta^{-1}_{\R^5}\mathcal{R}_{\Psi}\big),\qquad
\wt{\mathcal{R}}:=-\mathcal{R}_{\wOmega}+\mathcal{L}_3\big(\Delta^{-1}_{\R^5}\mathcal{R}_{\wPsi}\big)
+\wt{\mathcal{L}}_3\big(\Delta^{-1}_{\R^5}\mathcal{R}_{\Psi}\big).
\end{split}
\end{equation}
Notice that $\mathcal{L}_2(h)=\mathcal{N}_5(\Psi^\ast,h)$, $\mathcal{L}_3(h)=\mathcal{N}_5(h,\Omega^\ast)$,
$\wt{\mathcal{L}}_2(h)=\mathcal{N}_5(\wt\Psi^\ast,h)$, $\wt{\mathcal{L}}_3(h)=\mathcal{N}_5(h,\wt\Omega^\ast)$, in agreement with \eqref{sf5}. The two
identities $\Psi^\ast=\Delta^{-1}_{\R^5}\Omega^\ast$ and
$\wt\Psi^\ast=\Delta^{-1}_{\R^5}\wt\Omega^\ast$ hold because
$\Delta_{\R^5}\Psi_5=\Omega^\ast+\mathcal{R}_\Psi$ and
$\Delta_{\R^5}\wt\Psi_5=\wt\Omega^\ast+\mathcal{R}_{\wPsi}$, by the first equations of
\eqref{sf3} and \eqref{sf3.2}. The shift from $\Psi_5$ to $\Psi^\ast$ in $\mathcal{L}_2$, and from
$\wt\Psi_5$ to $\wt\Psi^\ast$ in $\wt{\mathcal{L}}_2$, allows us to absorb the terms
$\mathcal{N}_5(\Delta^{-1}_{\R^5}\mathcal{R}_{\Psi},\cdot)$ and
$\mathcal{N}_5(\Delta^{-1}_{\R^5}\mathcal{R}_{\wPsi},\cdot)$ produced by \eqref{cl1.5}.

We solve \eqref{cl2}--\eqref{cl2.5} for the axi-symmetric corrections $f,\wt f$ and the scalar $\delta\lambda$ using the $Z'$ bounds on $\mathcal R$ and $\wt{\mathcal R}$. The corrections determine the solutions by
\begin{equation}\label{cl4}
\begin{split}
&\omega=E^{-1}f,\qquad \psi=E^{-1}\big[\Delta^{-1}_{\R^5}(f-\mathcal{R}_{\Psi})\big],\\
&\wt\omega=E^{-1}\wt f,\qquad \wt\psi=E^{-1}\big[\Delta^{-1}_{\R^5}(\wt f-\mathcal{R}_{\wPsi})\big],
\qquad\wt\lambda=\wt\lambda_{app}+\delta\lambda.
\end{split}
\end{equation}

We start by inverting the operators $\mathcal{S}_\lambda:=\mathcal{O}_{\R^5}+\mathcal{L}_2+\mathcal{L}_3\Delta^{-1}_{\R^5}-\lambda I$, $\lambda\geq 0$. For this we define the Banach spaces $Z'$ and $Z_0\subseteq Z$ of axi-symmetric functions on $\R^5$ by the norms
\begin{equation}\label{robu1}
\begin{split}
&\|F\|_{Z'}:=\||\nabla|^{-1}\langle\nabla\rangle^{-1}F\|_{L^2(\R^5)},\\
&\|f\|_{Z}:=\|(\Delta_{\R^5}-1)f\|_{Z'}=\||\nabla|^{-1}\langle\nabla\rangle f\|_{L^2(\R^5)},\quad \|f\|_{Z_0}:=\|(\Delta_{\R^5}-1)f\|_{Z'}+\delta\|Sf\|_{Z'},
\end{split}
\end{equation}
where $S:=x_\al\partial_\al$ is the scaling vector-field and $\delta\in(0,1]$ is a small parameter. Here $|\nabla|:=(-\Delta_{\R^5})^{1/2}$ and $\langle\nabla\rangle:=(1-\Delta_{\R^5})^{1/2}$. Unless another domain is indicated, Lebesgue and Sobolev norms and duality pairings in Sections~2 and~3 are taken on $\R^5$, with measure $dx\,dz$.

\begin{table}[!htbp]
\centering
\setlength{\arraycolsep}{12pt}
\setlength{\belowcaptionskip}{0pt}
\caption{Quantitative bounds for Proposition~\ref{PropCentral2}.}\label{close:profiles}
\[
\begin{array}{|c|c|}
\hline
\rule{0pt}{2.6ex}\text{nonlinear solution}&\text{eigenpair}\\[4pt] \hline
\|\mathcal R\|_{Z'}\leq6\times10^{-12}
&\|\wt{\mathcal R}\|_{Z'}\leq1.5\times10^{-12}\\[8pt]
\|\mathcal R_\Psi\|_Z\leq1.251\times10^{-9}
&\|\nabla\wPsi^\ast\|_{L^2}\in[0.88,0.89]\\[8pt]
\|r\Omega^\ast\|_{L^\infty}\leq13
&\big\|\sqrt{r^2+z^2}\,\wOmega^\ast\big\|_{L^2}\leq22.55\\[8pt]
\big\|(-r\partial_z\Psi_5,\,r\partial_r\Psi_5+2\Psi_5)\big\|_{L^\infty}\leq22.45
&\big(\|a_\ast\|_{L^2}^2+\|b_\ast\|_{L^2}^2\big)^{1/2}\leq14.4\\[8pt]
\hline
\end{array}
\]
\end{table}

The test fields $a_\ast,b_\ast$ used to estimate the eigenvalue correction are defined by
\begin{equation}\label{close:test-fields}
\begin{split}
a_\ast=&\wOmega^\ast
+r\partial_z\Psi^\ast\Big(\partial_r\wPsi^\ast+\frac2r\wPsi^\ast\Big)
-(r\partial_r\Psi^\ast+2\Psi^\ast)\partial_z\wPsi^\ast,\\
b_\ast=&\Big(\frac r2\wOmega^\ast-r\Omega^\ast\partial_z\wPsi^\ast,\,
\frac z2\wOmega^\ast+r\Omega^\ast
\Big(\partial_r\wPsi^\ast+\frac2r\wPsi^\ast\Big)\Big).
\end{split}
\end{equation}
All norms use the full measure on $\R^5$, including the exterior region;
for an even integrand this is $4\pi^2r^3\,dr\,dz$ on $r,z>0$.
Vector norms are Euclidean norms of the displayed components.

\begin{lemma}\label{OrnUhl}
(i) For any $\lambda \geq 0$ the operator $\mathcal{O}_{\R^5}-\lambda I:Z_0\to Z'$ is a bounded linear isomorphism.

(ii) Let $K:=\mathcal{L}_2+\mathcal{L}_3\Delta^{-1}_{\R^5}$. Then $K:Z\to Z'$ is a bounded operator and 
\begin{equation}\label{robu1.5}
(\mathcal{O}_{\R^5}-\lambda I)^{-1}K:Z_0\to Z_0\qquad\text{is a bounded compact operator.}
\end{equation}
\end{lemma}

\begin{proof} (i) The Ornstein--Uhlenbeck operator is well understood, see for example \cite[Appendix A]{GaWa}. We provide all the details here mostly for the sake of completeness.

Since $|\xi|^2\leq\langle\xi\rangle^2$ we have
$\|\Delta_{\R^5}f\|_{Z'}=\||\nabla|\langle\nabla\rangle^{-1}f\|_{L^2}\leq\|f\|_{Z}$ and
$\|f\|_{Z'}\leq\|f\|_{Z}$, while $\|Sf\|_{Z'}\leq\delta^{-1}\|f\|_{Z_0}$. Therefore $\mathcal{O}_{\R^5}-\lambda I:Z_0\to Z'$ is a bounded linear map,
\begin{equation}\label{sf12}
\|(\mathcal{O}_{\R^5}-\lambda I)f\|_{Z'}\lesssim_{\delta,\lambda}\|f\|_{Z_0}.
\end{equation}

Now let $f\in Z_0$ and $g=(\mathcal{O}_{\R^5}-\lambda I)f$. We prove that 
\begin{equation}\label{sf10.7}
\|f\|_{Z_0}\lesssim_{\lambda}\|g\|_{Z'}.
\end{equation}
In the Fourier space, we have
\begin{equation}\label{sf10.5}
-|\xi|^2\widehat{f}(\xi)-\frac{1}{2}\,\xi\cdot\nabla_\xi\widehat{f}(\xi)-(1+\lambda)\widehat{f}(\xi)
=\widehat{g}(\xi).
\end{equation}
In polar coordinates, letting $\xi=s\omega$, $s=|\xi|\in(0,\infty)$, $\omega\in\mathbb{S}^4$, this becomes
\begin{equation}\label{sf10.6}
s\partial_s\widehat{f}(s\omega)+2(s^2+1+\lambda)\widehat{f}(s\omega)=-2\widehat{g}(s\omega).
\end{equation}
Multiplying \eqref{sf10.6} by
$s^2\overline{\widehat{f}(s\omega)}$ and taking real parts, we obtain, for almost every $s,\omega$,
\begin{equation*}
\frac{1}{2}\partial_s\big(s^3|\widehat{f}|^2\big)
+s^2\Big(2s^2+\frac{1}{2}+2\lambda\Big)|\widehat{f}|^2
=-2s^2\,\mathrm{Re}\big(\widehat{g}\,\overline{\widehat{f}}\big).
\end{equation*}
Fix $\omega$ such that $\int_0^\infty s^2(1+s^2)|\widehat{f}(s\omega)|^2\,ds<\infty$, which holds for
almost every $\omega$. We integrate the identity above over $s\in(\epsilon_n,R_n)$ along sequences
$\epsilon_n\to0$, $R_n\to\infty$ for which the boundary terms vanish, and then let $n\to\infty$ and
integrate in $\omega$, to get
\begin{equation}\label{sf10.8}
\int_{\mathbb{S}^4}\!\!\int_0^\infty\! s^2\Big(2s^2+\frac{1}{2}+2\lambda\Big)|\widehat{f}(s\omega)|^2\,ds\,d\omega
=-2\int_{\mathbb{S}^4}\!\!\int_0^\infty\! s^2\,\mathrm{Re}\big(\widehat{g}\,\overline{\widehat{f}}\big)(s\omega)\,ds\,d\omega.
\end{equation}
Since $2s^2+\tfrac{1}{2}+2\lambda\gtrsim_\lambda (1+s^2)$,
the left-hand side of \eqref{sf10.8} is $\gtrsim_\lambda\|f\|_{Z}^2$, while the
right-hand side is $\lesssim \|f\|_{Z}\|g\|_{Z'}$ by the
Cauchy--Schwarz inequality. Therefore $\|f\|_Z\lesssim_\lambda\|g\|_{Z'}$. Moreover, using \eqref{sf10.6} we have
$\widehat{Sf}=(2s^2-3+2\lambda)\widehat{f}+2\widehat{g}$. Therefore
$\|Sf\|_{Z'}\lesssim_\lambda \|f\|_{Z}+\|g\|_{Z'}$, and the desired bounds \eqref{sf10.7} follow.

In particular $\mathcal{O}_{\R^5}-\lambda I$ is injective. To prove that it is also surjective we multiply \eqref{sf10.6} by the integrating factor
$s^{1+2\lambda}e^{s^2}$ to obtain
\begin{equation*}
\partial_s\big\{s^{2+2\lambda}e^{s^2}\widehat{f}(s\omega)\big\}=-2s^{1+2\lambda}e^{s^2}\widehat{g}(s\omega).
\end{equation*}
We integrate and make the substitution $u=\alpha s$, thus
\begin{equation}\label{sf11}
\widehat{f}(s\omega)=A(\omega)s^{-2-2\lambda}e^{-s^2}
-2\int_0^1\alpha^{1+2\lambda}e^{-(1-\alpha^2)s^2}\widehat{g}(\alpha s\omega)\,d\alpha,
\end{equation}
for some function $A$ on $\mathbb{S}^4$. The first term on the right-hand side never belongs to $Z$
unless $A\equiv0$, since
$\int_{|\xi|\leq1}|\xi|^{-2}\langle\xi\rangle^2\,|\xi|^{-4-4\lambda}\,d\xi=\infty$.

We may therefore define
\begin{equation}\label{sf13}
\widehat{T_{\R^5}g}(\xi):=-2\int_0^1\alpha^{1+2\lambda}e^{-(1-\alpha^2)|\xi|^2}\widehat{g}(\alpha\xi)\,d\alpha,
\end{equation}
and let $\mathcal{A}$ denote the space of axi-symmetric functions $g$ with
$\widehat{g}\in C^\infty_c(\R^5\setminus\{0\})$, which is a dense subspace of $Z'$. If $g\in\mathcal{A}$, with $\widehat{g}$ supported in $\{a\leq|\xi|\leq b\}$ then
$\widehat{T_{\R^5}g}$ is smooth and vanishes in the ball $\{|\xi|<a\}$. Moreover, if $|\xi|\geq b$
then $\widehat{g}(\alpha\xi)\neq0$ forces $\alpha\leq b/|\xi|$, so $1-\alpha^2\geq1-b^2|\xi|^{-2}$
and $|\widehat{T_{\R^5}g}(\xi)|\leq2e^{b^2}\|\widehat{g}\|_{L^\infty}e^{-|\xi|^2}$. In particular 
$T_{\R^5}g\in Z_0$. Finally, \eqref{sf13} is equivalent to
\begin{equation*}
s^{2+2\lambda}e^{s^2}\widehat{T_{\R^5}g}(s\omega)=-2\int_0^su^{1+2\lambda}e^{u^2}\widehat{g}(u\omega)\,du,
\end{equation*}
so $\widehat{T_{\R^5}g}$ satisfies \eqref{sf10.6}, i.e. $(\mathcal{O}_{\R^5}-\lambda I)T_{\R^5}g=g$.

We may therefore apply the bound \eqref{sf10.7} with $f=T_{\R^5}g$, which gives
$\|T_{\R^5}g\|_{Z_0}\lesssim_\lambda\|g\|_{Z'}$ for any $g\in\mathcal{A}$. Since $\mathcal{A}$ is
dense in $Z'$, the operator $T_{\R^5}$ extends uniquely to a bounded operator $T_{\R^5}:Z'\to Z_0$,
and the identity $(\mathcal{O}_{\R^5}-\lambda I)T_{\R^5}=I$ extends to $Z'$ by continuity, using
\eqref{sf12}. In particular $\mathcal{O}_{\R^5}-\lambda I$ is surjective, so
$\mathcal{O}_{\R^5}-\lambda I:Z_0\to Z'$ is a bounded linear isomorphism, with bounded inverse
$T_{\R^5}$.

(ii) Let $G:=\Delta^{-1}_{\R^5}f$ and set
\begin{equation}\label{ou1}
\begin{split}
&U:=\big(-x_1\partial_z\Psi^\ast,\ldots,-x_4\partial_z\Psi^\ast,\ x_j\partial_j\Psi^\ast+2\Psi^\ast\big),\qquad
a:=2\partial_z\Psi^\ast,\\
&W:=\big(x_1\partial_z\Omega^\ast,\ldots,x_4\partial_z\Omega^\ast,\ -x_j\partial_j\Omega^\ast\big),\qquad
b:=-2\partial_z\Omega^\ast.
\end{split}
\end{equation}
A direct computation gives $\partial_\al U_\al=-a$ and $\partial_\al W_\al=-2b$, so \eqref{cl3} can be
rewritten in divergence form,
\begin{equation}\label{ou2}
\mathcal{L}_2(f)=U\cdot\nabla f=\partial_\al\big(U_\al f\big)+af,\qquad
\mathcal{L}_3(G)=W\cdot\nabla G+2\partial_z\Omega^\ast G=\partial_\al\big(W_\al G\big)+bG,
\end{equation}
and therefore
\begin{equation}\label{ou3}
Kf=\partial_\al\big(U_\al f+W_\al G\big)+af+bG.
\end{equation}
By \eqref{sf3.3} and the corresponding bounds on the first derivatives,
\begin{equation}\label{ou4}
|U|\lesssim\langle(x,z)\rangle^{-1},\qquad|W|\lesssim\langle(x,z)\rangle^{-3},\qquad
|a|\lesssim\langle(x,z)\rangle^{-2},\qquad|b|\lesssim\langle(x,z)\rangle^{-4},
\end{equation}
so that $U\in L^\infty$, $W\in L^\infty\cap L^5$, $a\in L^\infty\cap L^5$, $b\in L^\infty\cap
L^{5/2}$, and all four tend to $0$ at infinity.

We first show that $K:Z\to Z'$ is bounded. The definition \eqref{robu1} gives at once
\begin{equation}\label{ou5}
\|\partial_\al H_\al\|_{Z'}\leq\|\langle\nabla\rangle^{-1}H\|_{L^2}\leq\|H\|_{L^2},\qquad
\|F\|_{Z'}\leq\||\nabla|^{-1}F\|_{L^2}\lesssim\|F\|_{L^{10/7}},
\end{equation}
the last bound by the Hardy--Littlewood--Sobolev inequality in $\R^5$, and
\begin{equation}\label{ou6}
\|f\|_{L^2}\leq\|f\|_Z,\qquad\|\nabla G\|_{L^2}=\||\nabla|^{-1}f\|_{L^2}\leq\|f\|_Z,\qquad
\|G\|_{L^{10/3}}\lesssim\|f\|_Z,
\end{equation}
the last bound by the Sobolev embedding $\dot H^1(\R^5)\hookrightarrow L^{10/3}(\R^5)$. Combining
\eqref{ou3}--\eqref{ou6} with H\"older's inequality, we have
\begin{equation*}
\begin{split}
\|Kf\|_{Z'}&\lesssim\|Uf\|_{L^2}+\|WG\|_{L^2}+\|af\|_{L^{10/7}}+\|bG\|_{L^{10/7}}\\
&\lesssim\|U\|_{L^\infty}\|f\|_{L^2}+\|W\|_{L^5}\|G\|_{L^{10/3}}
+\|a\|_{L^5}\|f\|_{L^2}+\|b\|_{L^{5/2}}\|G\|_{L^{10/3}}\lesssim\|f\|_Z.
\end{split}
\end{equation*}

We now show that $K:Z\to Z'$ is compact. Let $f_k\rightharpoonup0$ weakly in $Z$ and
$G_k:=\Delta^{-1}_{\R^5}f_k$. By \eqref{ou6} the sequences $(f_k)$ and $(\nabla G_k)$ are bounded in
$L^2$ and converge weakly to $0$; in particular $G_k\to0$ strongly in $L^2(B_n)$ for every $n\geq1$,
by the Rellich--Kondrachov theorem, where $B_n:=\{(x,z)\in\R^5: |(x,z)|<n\}$.

The terms containing $G_k$ are estimated using \eqref{ou4}: for any $n\geq1$,
\begin{equation*}
\begin{split}
\|WG_k\|_{L^2}&\leq\|W\|_{L^\infty}\|G_k\|_{L^2(B_n)}+\|W\|_{L^5(\R^5\setminus B_n)}\|G_k\|_{L^{10/3}},\\
\|bG_k\|_{L^{10/7}}&\leq\|b\|_{L^5(B_n)}\|G_k\|_{L^2(B_n)}+\|b\|_{L^{5/2}(\R^5\setminus B_n)}\|G_k\|_{L^{10/3}}.
\end{split}
\end{equation*}
The second terms are arbitrarily small for $n$ large, uniformly in $k$, and the first terms vanish as
$k\to\infty$ for any fixed $n$.

For the terms containing $f_k$ we use the smoothing already present in \eqref{ou5}, together with the
standard fact that $\theta(\nabla)c(x,z)$ is a compact operator on $L^2(\R^5)$ whenever
$c,\theta\in L^\infty(\R^5)$ both tend to $0$ at infinity. Taking $\theta=\langle\xi\rangle^{-1}$ and $c=U_\al$ gives $\|\langle\nabla\rangle^{-1}(U_\al f_k)\|_{L^2}\to0$. For the term $af_k$ we split $a=a\chi_n+a(1-\chi_n)$, where $\chi_n$ is a cutoff to $B_n$. The second piece contributes
$\lesssim\|a\|_{L^5(\R^5\setminus B_n)}\|f_k\|_{L^2}$, which is small for $n$ large uniformly in $k$, while for the first piece we write $|\nabla|^{-1}\langle\nabla\rangle^{-1}=\theta_1(\nabla)+\theta_2(\nabla)$
with $\theta_1(\xi):=|\xi|^{-1}\langle\xi\rangle^{-1}\varphi(\xi)$ and $\varphi$ a cutoff to
$\{|\xi|\leq1\}$. The symbol $\theta_2$ is bounded and tends to $0$ at infinity, so
$\theta_2(\nabla)(a\chi_n)$ is compact; and $\theta_1\in L^2(\R^5)$, $a\chi_n\in L^2(\R^5)$, so
$\theta_1(\nabla)(a\chi_n\,\cdot)$ is a Hilbert--Schmidt operator on $L^2$. In both cases the contribution vanishes as $k\to\infty$.

Therefore $\|Kf_k\|_{Z'}\to0$, i.e. $K:Z\to Z'$ is compact. Since $Z_0\hookrightarrow Z$ is a bounded embedding and $(\mathcal{O}_{\R^5}-\lambda I)^{-1}:Z'\to Z_0$ is bounded by part (i), the operator
$(\mathcal{O}_{\R^5}-\lambda I)^{-1}K:Z_0\to Z_0$ is compact, which is \eqref{robu1.5}.
\end{proof}

In our applications, Lemma~\ref{OrnUhl} is used together with the following general lemma:

\begin{lemma}\label{lem:compactness}
Let $X$ and $Y$ be Banach spaces and suppose $S:X\to Y$ can be written as
\[
  S\;=\;L+K,
\]
where $L:X\to Y$ is a bounded linear isomorphism and $K:X\to Y$ is a linear relatively compact perturbation, i.e. the operator $L^{-1}K\colon X\to X$ is compact. Assume that $S$ is injective or $S$ is surjective.  Then $S$ is bijective from $X$ onto $Y$ and its inverse is bounded.
\end{lemma}

\begin{proof}
Notice that $S=L(I+C)$, where $C:=L^{-1}K:X\to X$.  The operator $C$ is compact and $I+C:X\to X$ is injective (respectively surjective). By the Fredholm alternative, $I+C$ is bijective and $(I+C)^{-1}$ is a bounded linear operator. Therefore $S$ is bijective with bounded inverse, as desired.
\end{proof}

\subsection{Weak coercivity and invertibility of the linear operators}  We solve \eqref{cl2}--\eqref{cl2.5} by inverting the linearized operators. Throughout this subsection we work with real-valued functions. Let $Z_0^{odd}$, $Z^{odd}$, $(Z')^{odd}$ (respectively $Z_0^{even}$, $Z^{even}$, $(Z')^{even}$) denote the subspaces of functions in $Z_0$, $Z$, and $Z'$ which are odd in $z$ (respectively even in $z$), with the inherited norms.

We start with a general lemma:

\begin{lemma}\label{strongCoercive}
Assume $\lambda\geq 0$, $\#\in\{odd, even\}$, and $T:Z^\#\to(Z')^\#$ is a bounded linear operator with the property that $(\mathcal{O}_{\R^5}-\lambda I)^{-1}T$ is a compact operator on $Z_0^\#$. Let $U:=(\mathcal{O}_{\R^5}-\lambda I)+T:Z_0^\#\to(Z')^\#$ and assume that $U$ satisfies a strong coercive bound of the form
\begin{equation}\label{robu2}
\left\langle Uf,\Delta^{-1}_{\R^5}f\right\rangle\geq\gamma\|f\|_{Z^\#}^2,\qquad\gamma\in(0,1/2],
\end{equation}
for any $f\in Z_0^\#$. If the constant $\delta$ in \eqref{robu1} is sufficiently small then $U$ is an isomorphism and 
\begin{equation}\label{robu6}
\|U^{-1}g\|_{Z_0^\#}\leq (0.99\gamma)^{-1} \|g\|_{(Z')^\#}\qquad\text{ for any }g\in (Z')^{\#}.
\end{equation}
\end{lemma}

\begin{proof}
Note that the right-hand side of \eqref{robu2} involves only $\|f\|_{Z^\#}$, not $\|f\|_{Z_0^\#}$, since the anti-selfadjoint part of the operator with respect to the energy pairing does not contribute to the coercive bound. We therefore recover the scaling derivative separately from the equation.

By \eqref{robu1}, for $g\in(Z')^\#$ and $h\in Z^\#$,
\begin{equation}\label{robu2.5}
\big|\langle g,\Delta^{-1}_{\R^5}h\rangle\big|
\leq \|g\|_{(Z')^\#}\|h\|_{Z^\#}.
\end{equation}
Combining this with \eqref{robu2} gives $\gamma\|f\|_{Z^\#}\leq\|Uf\|_{(Z')^\#}$ for every $f\in Z_0^\#$, so $U$ is injective. The Ornstein--Uhlenbeck inverse preserves parity, hence Lemma~\ref{OrnUhl} (i) applies on $Z_0^\#$ and $(Z')^\#$. The assumed compactness and Lemma~\ref{lem:compactness} therefore imply that $U$ is bijective.

To bound the full $Z_0^\#$ norm, we use the identity $Sf=2Uf-2\Delta_{\R^5}f-(3-2\lambda)f-2Tf$, thus
\[
\|Sf\|_{(Z')^\#}\leq2\|Uf\|_{(Z')^\#}+C\|f\|_{Z^\#},
\qquad C:=2+|3-2\lambda|+2\|T\|_{Z^\#\to(Z')^\#}.
\]
Consequently,
\[
\|f\|_{Z_0^\#}
=\|f\|_{Z^\#}+\delta\|Sf\|_{(Z')^\#}
\leq\frac{1+\delta(C+2\gamma)}{\gamma}\|Uf\|_{(Z')^\#}.
\]
The desired bounds \eqref{robu6} follow if $\delta\leq[99(C+1)]^{-1}$.
\end{proof}

In our applications it is not possible to prove strong coercive estimates like \eqref{robu2}, so we need to find suitable substitutes. Let
\begin{equation}\label{robu5}
\mathcal{S}_\lambda:=\mathcal{O}_{\R^5}-\lambda I+\mathcal{L}_2+\mathcal{L}_3\Delta^{-1}_{\R^5}:Z_0\to Z',
\end{equation}
where $\lambda\geq0$ and $\mathcal{L}_2,\mathcal{L}_3$ are defined in \eqref{cl3}. Since $\Psi^\ast$ and $\Omega^\ast$ are odd in $z$, both $K$ and $\mathcal{S}_\lambda$ preserve parity. For the nonlinear problem we use the following:

\begin{lemma}\label{weakCoercive1}
Assume that the operator $\mathcal{S}_0$ satisfies a weak coercive inequality of the form
\begin{equation}\label{robu7}
\left\langle \mathcal S_0f,\Delta^{-1}_{\R^5}f\right\rangle+\sum_{j=1}^n(\ell_j(f))^2\geq\gamma\|f\|_{Z^{odd}}^2,\qquad\gamma\in(0,1/2],
\end{equation}
for any \(f\in Z_0^{odd}\), where
\begin{equation}\label{robu8}
\ell_j(f)=\langle q_j,\Delta^{-1}_{\R^5}f\rangle,\qquad q_j\in (Z')^{odd}.
\end{equation}
Assume also that there are functions $\Phi_j\in Z_0^{odd}$ such that 
\begin{equation}\label{robu11}
\mathcal{S}_0\Phi_j=q_j+e_j,
\end{equation}
for $j\in\{1,\ldots,n\}$. Assume that
\begin{equation}\label{robu11.1}
\delta\leq \frac{0.01}{6+2\|K\|_{Z^{odd}\to(Z')^{odd}}+2\sum_{j=1}^n\|q_j\|_{(Z')^{odd}}^2},\qquad \sum_{j=1}^n\|q_j\|_{(Z')^{odd}}\|e_j\|_{(Z')^{odd}}\leq 0.01\gamma.
\end{equation}
Then $\mathcal{S}_0$ is a linear isomorphism between $Z_0^{odd}$ and $(Z')^{odd}$, and
\begin{equation}\label{robu11.2}
\|\mathcal{S}_0^{-1}g\|_{Z_0^{odd}}
\leq\frac{1+\sum_{j=1}^n\|\Phi_j\|_{Z_0^{odd}}\|q_j\|_{(Z')^{odd}}}{0.98\gamma}
\|g\|_{(Z')^{odd}}.
\end{equation}
\end{lemma}

\begin{proof}
For any $f\in Z_0^{odd}$,
\begin{equation}\label{robu9}
\begin{split}
&\left\langle \mathcal{S}_+f,\Delta^{-1}_{\R^5}f\right\rangle\geq\gamma\|f\|_{Z^{odd}}^2,\qquad \mathcal{S}_+:=\mathcal{S}_0+Q,\qquad Qf:=\sum_{j=1}^n\ell_j(f)q_j,
\end{split}
\end{equation}
since $\langle Qf,\Delta^{-1}_{\R^5}f\rangle=\sum_{j=1}^n(\ell_j(f))^2$.
The operator $Q:Z^{odd}\to(Z')^{odd}$ has finite rank and
$\|Q\|_{Z^{odd}\to(Z')^{odd}}\leq\sum_{j=1}^n\|q_j\|_{(Z')^{odd}}^2$.
Thus Lemma~\ref{strongCoercive}, Lemma~\ref{OrnUhl}, and the first condition in
\eqref{robu11.1} imply that $\mathcal{S}_+:Z_0^{odd}\to(Z')^{odd}$ is an isomorphism and
\begin{equation}\label{robu10}
\|\mathcal{S}_+^{-1}g\|_{Z_0^{odd}}\leq (0.99\gamma)^{-1} \|g\|_{(Z')^{odd}}\qquad\text{ for any }g\in (Z')^{odd}.
\end{equation}
Therefore
\begin{equation}\label{robu12}
\mathcal{S}_0(I+G)=\mathcal{S}_++E,\qquad G(f):=\sum_{j=1}^n\ell_j(f)\Phi_j,\qquad E(f):=\sum_{j=1}^n\ell_j(f)e_j,
\end{equation}
where $G:Z_0^{odd}\to Z_0^{odd}$ and $E:Z_0^{odd}\to(Z')^{odd}$. By \eqref{robu2.5} and the second condition in \eqref{robu11.1},
\[
\|Ef\|_{(Z')^{odd}}\leq
\sum_{j=1}^n|\ell_j(f)|\,\|e_j\|_{(Z')^{odd}}
\leq\sum_{j=1}^n\|q_j\|_{(Z')^{odd}}\|e_j\|_{(Z')^{odd}}\|f\|_{Z^{odd}}
\leq0.01\gamma\|f\|_{Z_0^{odd}}.
\]
Therefore $\mathcal{S}_++E:Z_0^{odd}\to(Z')^{odd}$ is a linear isomorphism, with
$\|(\mathcal{S}_++E)^{-1}g\|_{Z_0^{odd}}\leq(0.98\gamma)^{-1}\|g\|_{(Z')^{odd}}$
for any $g\in(Z')^{odd}$. Using \eqref{robu12}, it follows that the operator $\mathcal{S}_0:Z_0^{odd}\to (Z')^{odd}$ is surjective. In view of Lemma~\ref{lem:compactness} this operator is a linear isomorphism and for any $F\in (Z')^{odd}$ there is a unique $f\in Z_0^{odd}$ such that
\begin{equation}\label{robu14}
\mathcal{S}_0(f)=F,\,\,\,f=g+G(g),\,\,\,\|g\|_{Z_0^{odd}}\leq (0.98\gamma)^{-1} \|F\|_{(Z')^{odd}}.
\end{equation}
The desired bounds \eqref{robu11.2} follow.
\end{proof}

We consider now the corresponding linear operator for the eigenvalue problem, which is
\begin{equation}\label{robu20}
T(\mu,f):=\mathcal{S}_{\widetilde{\lambda}_{app}}f-\mu\wOmega^\ast.
\end{equation}
We define a linear functional $m_0:(Z')^{even}\to\R$ and the associated Banach spaces $Z_{0,\ast}^{even}$, $Z^{even}_\ast$, and $(Z')^{even}_\ast$ by
\begin{equation}\label{robu21}
\begin{split}
&m_0(F):=\frac{\langle F,\Delta^{-1}_{\R^5}\wOmega^\ast\rangle}{\langle \wOmega^\ast,\Delta^{-1}_{\R^5}\wOmega^\ast\rangle},\qquad (Z')^{even}_\ast:=\{F\in (Z')^{even}:\,m_0(F)=0\},\\
&Z_{0,\ast}^{even}:=Z_0^{even}\cap (Z')^{even}_\ast,\qquad Z_{\ast}^{even}:=Z^{even}\cap (Z')^{even}_\ast.
\end{split}
\end{equation}
We define the associated projection $P_0:(Z')^{even}\to(Z')^{even}_\ast$ by
\begin{equation}\label{robu22}
P_0(F):=F-m_0(F)\wOmega^\ast.
\end{equation}
By a slight abuse of notation we also let $P_0$ denote the restricted projections $P_0: Z_0^{even}\to Z_{0,\ast}^{even}$ and $P_0: Z^{even}\to  Z^{even}_\ast$. 

For the operator \eqref{robu20} we use the following analogue of Lemma~\ref{weakCoercive1}.

\begin{lemma}\label{weakCoercive2}
Assume that $\wOmega^\ast\in Z_0^{even}$ and that
\begin{equation}\label{robu31}
\left\langle\mathcal{S}_{\widetilde{\lambda}_{app}}f,\Delta^{-1}_{\R^5}f\right\rangle
+\sum_{j=1}^{\widetilde{n}}(\widetilde{\ell}_j(f))^2\geq\gamma\|f\|_{Z^{even}}^2,\qquad \gamma\in(0,1/2],
\end{equation}
for every $f\in Z_{0,\ast}^{even}$, where
\begin{equation}\label{robu32}
\widetilde{\ell}_j(f):=\langle \widetilde{q}_j,\Delta^{-1}_{\R^5}f\rangle,\qquad \widetilde{q}_j\in(Z')^{even}.
\end{equation}
Assume also that there are pairs $(\widetilde{\mu}_j,\widetilde{\Phi}_j)\in\R\times Z_{0,\ast}^{even}$ such that
\begin{equation}\label{robu34}
T(\widetilde{\mu}_j,\widetilde{\Phi}_j)=\mathcal{S}_{\widetilde{\lambda}_{app}}\widetilde{\Phi}_j-\widetilde{\mu}_j\wOmega^\ast
=\widetilde{q}_j+\widetilde{e}_j,
\end{equation}
for $j\in\{1,\ldots,\widetilde{n}\}$. Assume that
\begin{equation}\label{robu37}
\begin{split}
&\delta\leq
\frac{0.01}{6+2\|K\|_{Z^{even}\to(Z')^{even}}
+2\sum_{j=1}^{\widetilde{n}}\|\widetilde{q}_j\|_{(Z')^{even}}^2}
\left(1+\frac{\|\wOmega^\ast\|_{Z^{even}}+\|S\wOmega^\ast\|_{(Z')^{even}}}
{\||\nabla|^{-1}\wOmega^\ast\|_{L^2}}\right)^{-1},\\
&\sum_{j=1}^{\widetilde{n}}\|\widetilde{q}_j\|_{(Z')^{even}}\|\widetilde{e}_j\|_{(Z')^{even}}\leq0.01\gamma.
\end{split}
\end{equation}
Then $T:\R\times Z_{0,\ast}^{even}\to(Z')^{even}$ is a linear isomorphism. Moreover, if $T(\mu,f)=F$ then
\begin{equation}\label{robu38}
\|f\|_{Z_0^{even}}\leq
\frac{1+\sum_{j=1}^{\widetilde{n}}\|\widetilde{\Phi}_j\|_{Z_0^{even}}\|\widetilde{q}_j\|_{(Z')^{even}}}{0.98\gamma}\|F\|_{(Z')^{even}},
\end{equation}
and
\begin{equation}\label{robu39}
\begin{split}
|\mu|&\leq |m_0(F)|
+\frac{\|\wOmega^\ast\|_{L^2}+\frac12\|S\wOmega^\ast\|_{(Z')^{even}}
+\|K\|_{Z^{even}\to(Z')^{even}}\|\wOmega^\ast\|_{Z^{even}}}
{\||\nabla|^{-1}\wOmega^\ast\|_{L^2}^2}\|f\|_{Z^{even}},\\
|m_0(F)|&\leq
\frac{\|\wOmega^\ast\|_{Z^{even}}}{\||\nabla|^{-1}\wOmega^\ast\|_{L^2}^2}\|F\|_{(Z')^{even}}.
\end{split}
\end{equation}
\end{lemma}

\begin{proof}
{\textbf{Step 1.}} We prove that $T_0:\R\times Z_{0,\ast}^{even}\to(Z')^{even}$ is an isomorphism, where
\[
T_0(\mu,f):=(\mathcal{O}_{\R^5}-\widetilde{\lambda}_{app}I)f-\mu\wOmega^\ast.
\]

Given $F\in(Z')^{even}$, Lemma~\ref{OrnUhl} provides
\[
g=(\mathcal{O}_{\R^5}-\widetilde{\lambda}_{app}I)^{-1}F,
\qquad
h=(\mathcal{O}_{\R^5}-\widetilde{\lambda}_{app}I)^{-1}\wOmega^\ast,
\]
both in $Z_0^{even}$. The equation $T_0(\mu,f)=F$ is equivalent to $f=g+\mu h$.
Thus the normalization $m_0(f)=0$ requires
$m_0(g)+\mu m_0(h)=0$, and it remains to show that $m_0(h)\neq0$.
By symmetry of $\Delta^{-1}_{\R^5}$ and the Ornstein--Uhlenbeck energy identity,
\[
\begin{split}
m_0(h)\langle\wOmega^\ast,\Delta^{-1}_{\R^5}\wOmega^\ast\rangle
&=\langle h,\Delta^{-1}_{\R^5}\wOmega^\ast\rangle
=\langle(\mathcal{O}_{\R^5}-\widetilde{\lambda}_{app}I)h,\Delta^{-1}_{\R^5}h\rangle=\|h\|_{L^2}^2+\left(\frac14+\widetilde{\lambda}_{app}\right)
\||\nabla|^{-1}h\|_{L^2}^2.
\end{split}
\]
Hence $m_0(h)\neq0$, and the unique solution of the equation $T_0(\mu,f)=F$ is
\[
\mu=-\frac{m_0(g)}{m_0(h)},\qquad
f=g-\frac{m_0(g)}{m_0(h)}h.
\]
Thus $T_0$ is a bijective linear bounded map, so it is an isomorphism.

Moreover, $T(\mu,f)=T_0(\mu,f)+Kf$, and
$K:Z^{even}\to(Z')^{even}$ is compact by the proof of Lemma~\ref{OrnUhl}.
Consequently $T$ and all its finite-rank perturbations are compact perturbations of
$T_0$. By Lemma~\ref{lem:compactness}, each is an isomorphism as soon as it is
injective or surjective.

{\textbf{Step 2.}} We prove now that the stabilized operator
$T_+:\R\times Z_{0,\ast}^{even}\to(Z')^{even}$,  
\begin{equation}\label{robu26}
T_+(\mu,f):=T(\mu,f)+\widetilde{Q}f
=\mathcal{S}_{\widetilde{\lambda}_{app}}f+\widetilde{Q}f-\mu\wOmega^\ast,\qquad \widetilde{Q}f:=\sum_{j=1}^{\widetilde{n}}\widetilde{\ell}_j(f)\widetilde{q}_j,
\end{equation}
is also an isomorphism, and
\begin{equation}\label{robu30}
\|f\|_{Z_0^{even}}\leq(0.99\gamma)^{-1}\|F\|_{(Z')^{even}}\qquad\text{ if }T_+(\mu,f)=F.
\end{equation}
Indeed,  if $T_+(\mu,f)=F$ then we pair with $\Delta^{-1}_{\R^5}f$, and use \eqref{robu31}, \eqref{robu2.5}, and the assumption $m_0(f)=0$ to see that
\begin{equation}\label{robu26.1}
\gamma\|f\|_{Z^{even}}\leq\|F\|_{(Z')^{even}}.
\end{equation}
If $F=0$, then $f=0$ and subsequently $\mu=0$. Thus $T_+$ is injective.
Since $\widetilde{Q}$ has finite rank, Step 1 and Lemma~\ref{lem:compactness} imply that $T_+$ is an isomorphism.

To prove \eqref{robu30}, we still need to estimate the scaling derivative $Sf$.
We first bound the scalar $\mu$ in terms of $\|F\|_{(Z')^{even}}$. Indeed, integration by parts in $\R^5$ gives
\[
\langle Sf,\Delta^{-1}_{\R^5}\wOmega^\ast\rangle
=-\langle S\wOmega^\ast,\Delta^{-1}_{\R^5}f\rangle
-7\langle f,\Delta^{-1}_{\R^5}\wOmega^\ast\rangle.
\]
Consequently, since $m_0(f)=0$,
\begin{equation}\label{robu27}
m_0(\mathcal{S}_{\widetilde{\lambda}_{app}}f)
=\frac{\langle f,\wOmega^\ast\rangle
-\langle S\wOmega^\ast,\Delta^{-1}_{\R^5}f\rangle/2
+\langle Kf,\Delta^{-1}_{\R^5}\wOmega^\ast\rangle}
{\langle\wOmega^\ast,\Delta^{-1}_{\R^5}\wOmega^\ast\rangle}.
\end{equation}
Applying $m_0$ to $T_+(\mu,f)=F$, we obtain
\begin{equation}\label{robu27.5}
\begin{split}
|\mu|&\leq |m_0(\mathcal{S}_{\widetilde{\lambda}_{app}}f)|
+|m_0(\widetilde{Q}f)|+|m_0(F)|\\
&\leq\frac{\|f\|_{Z^{even}}}{\||\nabla|^{-1}\wOmega^\ast\|_{L^2}^2}
\Big[\|\wOmega^\ast\|_{L^2}+\frac12\|S\wOmega^\ast\|_{(Z')^{even}}+M\|\wOmega^\ast\|_{Z^{even}}\Big]+\frac{\|\wOmega^\ast\|_{Z^{even}}\|F\|_{(Z')^{even}}}
{\||\nabla|^{-1}\wOmega^\ast\|_{L^2}^2},
\end{split}
\end{equation}
where $M:=\|K\|_{Z^{even}\to(Z')^{even}}+\sum_{j=1}^{\widetilde{n}}\|\widetilde{q}_j\|_{(Z')^{even}}^2$. We use
\[
\begin{split}
&\|\Delta_{\R^5} f\|_{(Z')^{even}},\ \|f\|_{(Z')^{even}}\leq\|f\|_{Z^{even}},\\
&\|\wOmega^\ast\|_{(Z')^{even}}\leq\||\nabla|^{-1}\wOmega^\ast\|_{L^2},\qquad
\|\wOmega^\ast\|_{L^2}\leq\|\wOmega^\ast\|_{Z^{even}}.
\end{split}
\]
Since $|3-2\widetilde{\lambda}_{app}|\leq3$ for the fixed approximate eigenvalue, the identity
\[
Sf=2F+2\mu\wOmega^\ast-2\Delta_{\R^5} f-(3-2\widetilde{\lambda}_{app})f-2Kf-2\widetilde{Q}f,
\]
and the bounds \eqref{robu27.5} give
\begin{equation}\label{robu28}
\begin{split}
\|Sf\|_{(Z')^{even}}
&\leq 2\|F\|_{(Z')^{even}}+2|\mu|\|\wOmega^\ast\|_{(Z')^{even}}
+(5+2M)\|f\|_{Z^{even}}\\
&\leq 2\Big(1+\frac{\|\wOmega^\ast\|_{Z^{even}}}
{\||\nabla|^{-1}\wOmega^\ast\|_{L^2}}\Big)\|F\|_{(Z')^{even}}\\
&+\Big[5+2M+
\frac{(2+2M)\|\wOmega^\ast\|_{Z^{even}}+\|S\wOmega^\ast\|_{(Z')^{even}}}
{\||\nabla|^{-1}\wOmega^\ast\|_{L^2}}\Big]\|f\|_{Z^{even}}.
\end{split}
\end{equation}
Combining this with \eqref{robu26.1}, $\gamma\leq1/2$, and the first condition in
\eqref{robu37} gives
\[
\begin{split}
\|f\|_{Z_0^{even}}
&=\|f\|_{Z^{even}}+\delta\|Sf\|_{(Z')^{even}}\leq(0.99\gamma)^{-1}\|F\|_{(Z')^{even}},
\end{split}
\]
which proves \eqref{robu30}.

{\textbf{Step 3.}} We prove that $T$ is an isomorphism and establish the bounds
\eqref{robu38}--\eqref{robu39}. Set $\widetilde{E}g:=\sum_{j=1}^{\widetilde{n}}\widetilde{\ell}_j(g)\widetilde{e}_j$.
By \eqref{robu2.5} and the second condition in \eqref{robu37},
\[
\|\widetilde{E}g\|_{(Z')^{even}}
\leq\sum_{j=1}^{\widetilde{n}}|\widetilde{\ell}_j(g)|\,
\|\widetilde{e}_j\|_{(Z')^{even}}
\leq\sum_{j=1}^{\widetilde{n}}\|\widetilde{q}_j\|_{(Z')^{even}}\|\widetilde{e}_j\|_{(Z')^{even}}\|g\|_{Z^{even}}
\leq0.01\gamma\|g\|_{Z_0^{even}}.
\]
Consider
\[
T_+(\mu_0,g)+\widetilde{E}g=F,\qquad g\in Z_{0,\ast}^{even}.
\]
Applying \eqref{robu30} to $T_+(\mu_0,g)=F-\widetilde{E}g$, we obtain
\[
0.99\gamma\|g\|_{Z_0^{even}}
\leq\|F-\widetilde{E}g\|_{(Z')^{even}}
\leq\|F\|_{(Z')^{even}}+0.01\gamma\|g\|_{Z_0^{even}},
\]
and therefore
\[
\|g\|_{Z_0^{even}}\leq(0.98\gamma)^{-1}\|F\|_{(Z')^{even}}.
\]
In particular, $F=0$ implies $g=0$ and then $\mu_0=0$.
Thus $(\mu_0,g)\mapsto T_+(\mu_0,g)+\widetilde{E}g$ is injective and, by Step 1,
an isomorphism. For every $F\in(Z')^{even}$ there is consequently a unique
pair $(\mu_0,g)$ solving this equation.

Define
\[
\mu=\mu_0+\sum_{j=1}^{\widetilde{n}}\widetilde{\ell}_j(g)\widetilde{\mu}_j,
\qquad
f=g+\sum_{j=1}^{\widetilde{n}}\widetilde{\ell}_j(g)\widetilde{\Phi}_j.
\]
Since each $\widetilde{\Phi}_j$ lies in $Z_{0,\ast}^{even}$, so does $f$.
The response identities \eqref{robu34} give
\begin{equation}\label{robu36}
\begin{split}
T(\mu,f)
&=T(\mu_0,g)+\sum_{j=1}^{\widetilde{n}}
\widetilde{\ell}_j(g)(\widetilde{q}_j+\widetilde{e}_j)=T_+(\mu_0,g)+\widetilde{E}g=F.
\end{split}
\end{equation}
Hence $T$ is surjective, and Step 1 implies that it is an isomorphism.

Finally, $|\widetilde{\ell}_j(g)|\leq
\|\widetilde{q}_j\|_{(Z')^{even}}\|g\|_{Z^{even}}$ gives
\[
\|f\|_{Z_0^{even}}
\leq\Big(1+\sum_{j=1}^{\widetilde{n}}
\|\widetilde{\Phi}_j\|_{Z_0^{even}}\|\widetilde{q}_j\|_{(Z')^{even}}\Big)
\|g\|_{Z_0^{even}}.
\]
Together with the bound for $g$, this proves \eqref{robu38}.
Applying $m_0$ to $T(\mu,f)=F$, we have
\[
\mu=m_0(\mathcal{S}_{\widetilde{\lambda}_{app}}f)-m_0(F).
\]
Since $\langle\wOmega^\ast,\Delta^{-1}_{\R^5}\wOmega^\ast\rangle
=-\||\nabla|^{-1}\wOmega^\ast\|_{L^2}^2$,
identity \eqref{robu27} and the duality estimate \eqref{robu2.5} give both bounds in
\eqref{robu39}.
\end{proof}

\subsection{Linear and bilinear estimates} Recall the bilinear operator defined in \eqref{cl3}
\begin{equation}\label{nl2}
\mathcal{N}_5(g,h):=x_j\partial_jg\cdot \partial_zh-\partial_zg\cdot x_j\partial_jh+2g\cdot\partial_zh.
\end{equation}

\begin{lemma}\label{nl:bilinear}
Let $C_b:=0.08$. Then for any $f\in Z$ and axi-symmetric $g\in L^2(\R^5)$, the distribution
$\mathcal N_5(\Delta^{-1}_{\R^5}f,g)$ belongs to $Z'$ and satisfies
\begin{equation}\label{nl4}
\|\mathcal N_5(\Delta^{-1}_{\R^5}f,g)\|_{Z'}
\leq C_b\|f\|_Z\|g\|_{L^2(\R^5)}.
\end{equation}
In particular,
\begin{equation}\label{nl4.1}
\|\mathcal N_5(\Delta^{-1}_{\R^5}f,g)\|_{Z'}
\leq C_b\|f\|_Z\|g\|_Z
\qquad\text{for all }f,g\in Z.
\end{equation}
\end{lemma}

\begin{proof}
Assume first that $f,g\in C^\infty_c(\R^5)$ are axi-symmetric, set $\psi=\Delta^{-1}_{\R^5}f$, and let $\phi\in\mathcal S(\R^5)$ be axi-symmetric.
We regard these functions as profiles in $(r,z)$ on both $\R^5$ and $\R^3$.
For any integrable axi-symmetric $H$,
\begin{equation}\label{nl8}
\int_{\R^5}H\,dx\,dz
=2\pi^2\int_0^\infty\int_{\R}H(r,z)r^3\,dz\,dr
=\pi\int_{\R^3}r^2H\,dy.
\end{equation}
Writing $e_r,e_\theta,e_z$ for the cylindrical frame on $\R^3$, define
\begin{equation}\label{nl9}
\begin{split}
u&:=\operatorname{curl}(r\psi e_\theta)
=-r\partial_z\psi\,e_r+(r\partial_r\psi+2\psi)e_z,\\
v&:=\operatorname{curl}(r\phi e_\theta)
=-r\partial_z\phi\,e_r+(r\partial_r\phi+2\phi)e_z.
\end{split}
\end{equation}
Direct calculation gives $\operatorname{div}u=0$ and
$\operatorname{curl}u=-rf\,e_\theta$. Moreover,
\[
r|u|^2=r^3\big((\partial_r\psi)^2+(\partial_z\psi)^2\big)
+\partial_r(2r^2\psi^2).
\]
The boundary term vanishes at $r=0$ and at infinity, since $\psi$ is smooth and
$|\partial^\alpha\psi|\lesssim\langle(x,z)\rangle^{-3-|\alpha|}$.
Thus \eqref{nl8}, together with
$\|\nabla u\|_{L^2(\R^3)}=\|\operatorname{curl}u\|_{L^2(\R^3)}$, yields
\begin{equation}\label{nl10}
\|u\|_{L^2(\R^3)}^2=\frac1\pi\|\nabla\psi\|_{L^2(\R^5)}^2,
\qquad
\|\nabla u\|_{L^2(\R^3)}^2=\frac1\pi\|f\|_{L^2(\R^5)}^2.
\end{equation}
Applying the same identities to $\phi$, we obtain
\begin{equation}\label{nl11}
\begin{split}
\|u\|_{H^1(\R^3)}^2=\frac1\pi\|f\|_Z^2,\qquad \|v\|_{H^1(\R^3)}^2=\frac1\pi
\big(\|\nabla\phi\|_{L^2(\R^5)}^2+\|\Delta_{\R^5}\phi\|_{L^2(\R^5)}^2\big),
\end{split}
\end{equation}
where $\|V\|_{H^1(\R^3)}^2:=\|V\|_{L^2(\R^3)}^2+\|\nabla V\|_{L^2(\R^3)}^2$.

Set $w:=-rg\,e_\theta$ on $\R^3$. Then
\[
u\times w=rg\,u^z e_r-rg\,u^r e_z,
\]
and, since $\partial_r(ru^r)+r\partial_zu^z=0$,
\[
\operatorname{curl}(u\times w)
=\big(\partial_z(rg\,u^z)+\partial_r(rg\,u^r)\big)e_\theta
=r(u^r\partial_rg+u^z\partial_zg)e_\theta
=r\mathcal N_5(\psi,g)e_\theta.
\]
Using \eqref{nl8} and integrating the curl by parts therefore gives
\begin{equation}\label{nl12}
\begin{split}
\langle\mathcal N_5(\psi,g),\phi\rangle_{\R^5}
&=\pi\int_{\R^3}\operatorname{curl}(u\times w)\cdot(r\phi e_\theta)\,dy=\pi\int_{\R^3}(u\times w)\cdot v\,dy.
\end{split}
\end{equation}
Also, by \eqref{nl8},
\begin{equation}\label{nl12.1}
\|w\|_{L^2(\R^3)}=\pi^{-1/2}\|g\|_{L^2(\R^5)}.
\end{equation}

We now obtain an explicit bound for the embedding $H^1(\R^3)\hookrightarrow L^4(\R^3)$. The sharp fractional Sobolev inequality \cite{Li,CT} gives
\begin{equation}\label{nl16}
\|a\|_{L^4(\R^3)}^2
\leq C_4^2\big\||\nabla|^{3/4}a\big\|_{L^2(\R^3)}^2,
\qquad C_4^2:=\frac{8\Gamma(3/4)^2}{5\pi^2}.
\end{equation}
Since $\sup_{s\geq0}s^{3/2}/(1+s^2)=3^{3/4}/4$, for any vector field
$V\in H^1(\R^3;\R^3)$ with Cartesian components $V_k$,
\begin{equation}\label{nl17}
\begin{split}
\|V\|_{L^4(\R^3)}^2
&=\big\||V|^2\big\|_{L^2(\R^3)}
\leq\sum_{k=1}^3\|V_k\|_{L^4(\R^3)}^2\leq C_4^2\sum_{k=1}^3\big\||\nabla|^{3/4}V_k\big\|_{L^2(\R^3)}^2
\leq\frac{3^{3/4}C_4^2}{4}\|V\|_{H^1(\R^3)}^2.
\end{split}
\end{equation}

For axi-symmetric $F$, the definition of $Z'$ gives
\begin{equation}\label{nl13}
\|F\|_{Z'}
=\sup_{\substack{\phi\in\mathcal S(\R^5)\ \text{axi-symmetric}\\
\|\nabla\phi\|_{L^2(\R^5)}^2+\|\Delta_{\R^5}\phi\|_{L^2(\R^5)}^2\leq1}}
|\langle F,\phi\rangle_{\R^5}|.
\end{equation}
For a test function in this supremum, \eqref{nl11} gives
$\|v\|_{H^1(\R^3)}\leq\pi^{-1/2}$.
Combining \eqref{nl12}, \eqref{nl17}, \eqref{nl11}, and \eqref{nl12.1}, we obtain
\begin{equation}\label{nl18}
\begin{split}
|\langle\mathcal N_5(\psi,g),\phi\rangle_{\R^5}|
&\leq\pi\|u\|_{L^4(\R^3)}\|w\|_{L^2(\R^3)}\|v\|_{L^4(\R^3)}\\
&\leq\frac{\pi\,3^{3/4}C_4^2}{4}
\|u\|_{H^1(\R^3)}\|w\|_{L^2(\R^3)}\|v\|_{H^1(\R^3)}\\
&\leq\frac{3^{3/4}C_4^2}{4\sqrt\pi}\|f\|_Z\|g\|_{L^2(\R^5)}
\leq C_b\|f\|_Z\|g\|_{L^2(\R^5)},
\end{split}
\end{equation}
since $3^{3/4}C_4^2/(4\sqrt\pi)=2\,3^{3/4}\Gamma(3/4)^2/(5\pi^{5/2})<0.08=C_b$.
Taking the supremum proves \eqref{nl4} for $f,g\in C^\infty_c(\R^5)$.

Axi-symmetric functions in $C^\infty_c(\R^5)$ are dense in $Z$ and in the
axi-symmetric subspace of $L^2(\R^5)$, so \eqref{nl4} extends the bilinear map
continuously to these spaces. To identify the extension, we write, for $\psi=\Delta^{-1}_{\R^5}f$,
\[
\mathcal N_5(\psi,g)
=\operatorname{div}_{\R^5}\!\left[g\big(-x\,\partial_z\psi,\,
x\cdot\nabla_x\psi+2\psi\big)\right]
+2(\partial_z\psi)g.
\]
Since $\nabla\psi\in L^2(\R^5)$ and $\psi\in L^{10/3}(\R^5)$, the products on the
right are locally integrable and depend continuously on $f\in Z$, $g\in L^2(\R^5)$
in the sense of distributions. Thus the extension agrees with the nonlinearity
in \eqref{nl2}. Finally, \eqref{nl4.1} follows from $\|g\|_{L^2(\R^5)}\leq\|g\|_Z$.
\end{proof}

Taking $g=f$ in \eqref{nl4.1} bounds the nonlinear term in \eqref{cl2} by
$C_b\|f\|_Z^2$. Applying the same estimate to the two ordered terms in \eqref{cl2.5}
bounds their sum by $2C_b\|f\|_Z\|\widetilde f\|_Z$.

\begin{corollary}\label{nl:linear}
For any $f\in Z^{odd}$, the functions $\wt{\mathcal L}_3(\Delta^{-1}_{\R^5}f)$ and
$\wt{\mathcal L}_2f$ belong to $(Z')^{even}$ and satisfy
\begin{equation}\label{nl:linear-bounds}
\begin{split}
\|\wt{\mathcal L}_3(\Delta^{-1}_{\R^5}f)\|_{(Z')^{even}}
&\leq C_b\|\wt\Omega^\ast\|_{L^2(\R^5)}\|f\|_{Z^{odd}}
\leq0.0808\|f\|_{Z^{odd}},\\
\|\wt{\mathcal L}_2f\|_{(Z')^{even}}
&\leq C_b\|\wt\Omega^\ast\|_{Z^{even}}\|f\|_{L^2(\R^5)}.
\end{split}
\end{equation}
\end{corollary}

\begin{proof}
By \eqref{cl3}--\eqref{cl3.3},
\[
\wt{\mathcal L}_3(\Delta^{-1}_{\R^5}f)=\mathcal N_5(\Delta^{-1}_{\R^5}f,\wt\Omega^\ast),
\qquad
\wt{\mathcal L}_2f=\mathcal N_5(\Delta^{-1}_{\R^5}\wt\Omega^\ast,f).
\]
Both expressions are even in $z$, and the bounds follow from \eqref{nl4} and
the bound $\|\wt\Omega^\ast\|_{L^2(\R^5)}\leq1.01$ of Proposition~\ref{PropCentral2}.
\end{proof}

\subsection{Quantitative bounds and the correction equations}\label{sec:closing}

We solve \eqref{cl2}--\eqref{cl2.5} with the parameters
\begin{equation}\label{close:parameters}
\gamma=\frac1{10},\qquad n=26,\qquad \widetilde n=32,\qquad
\delta=10^{-6}.
\end{equation}
Theorem~\ref{c:thm:main} gives the weak coercivity inequalities \eqref{robu7} and \eqref{robu31}. The remaining inputs are the profile bounds of Proposition~\ref{PropCentral2} and the functional and response bounds in Table~\ref{close:coercivity}, proved in Subsections~\ref{cap:coercivity-verification} and~\ref{sec:responses}.

The bounds on $\mathcal R$ and $\wt{\mathcal R}$ in Table~\ref{close:profiles}
follow from the more detailed residual bounds in Subsection~\ref{sec:main-residual}. The bound on
$\mathcal R_\Psi$ also controls the passage from $\Psi_5$ to $\Psi^\ast$
in the coercivity estimates. Indeed, this passage changes the coercivity constant by at most
$C_b\|\mathcal R_\Psi\|_Z<1.001\times10^{-10}$, and the
calculation with $\Psi_5$ gives $\gamma_0=0.100001$, so we may take $\gamma=0.1$.

Section~\ref{sec:coercivity} constructs the sources in \eqref{robu8} and \eqref{robu32} by localizing the quadratic forms. The odd problem gives $26$ vector functionals. In the even problem the vector and scalar forms have separate partitions and give $23$ and $9$ functionals, respectively.

The responses $\Phi_j\in Z_0^{odd}$ and $(\widetilde\mu_j,\widetilde\Phi_j)\in\R\times Z_{0,\ast}^{even}$ are constructed in Subsection~\ref{sec:responses}. Their errors are
\begin{equation}\label{close:response-identities}
e_j=\mathcal S_0\Phi_j-q_j,\qquad
\widetilde e_j=
\mathcal S_{\widetilde\lambda_{app}}\widetilde\Phi_j
-\widetilde\mu_j\wOmega^\ast-\widetilde q_j.
\end{equation}
In particular $m_0(\widetilde\Phi_j)=0$. The errors in \eqref{close:response-identities} use the exact sources and include the defects of the approximate potentials, the cutoff terms, and the exterior contributions.

\begin{table}[!htbp]
\centering
\setlength{\arraycolsep}{12pt}
\setlength{\belowcaptionskip}{4pt}
\caption{Functional and response bounds, with $\delta=10^{-6}$ in the $Z_0$ norm.}
\label{close:coercivity}\label{close:responses}
\[
\begin{array}{|c|c|}
\hline
\rule{0pt}{2.6ex}\text{odd problem}&\text{even problem}\\[4pt] \hline
\displaystyle\sum_{j=1}^{26}\|q_j\|_{(Z')^{odd}}\leq19.2
&\displaystyle\sum_{j=1}^{32}\|\widetilde q_j\|_{(Z')^{even}}\leq34\\[8pt]
\displaystyle\sum_{j=1}^{26}\|q_j\|_{(Z')^{odd}}^2\leq19.86
&\displaystyle\sum_{j=1}^{32}\|\widetilde q_j\|_{(Z')^{even}}^2\leq61\\[8pt]
\displaystyle\sum_{j=1}^{26}\|q_j\|_{(Z')^{odd}}\|e_j\|_{(Z')^{odd}}\leq5\times10^{-4}
&\displaystyle\sum_{j=1}^{32}\|\widetilde q_j\|_{(Z')^{even}}\|\widetilde e_j\|_{(Z')^{even}}\leq10^{-3}\\[8pt]
\displaystyle\sum_{j=1}^{26}\|\Phi_j\|_{Z^{odd}}\|q_j\|_{(Z')^{odd}}\leq15.25
&\displaystyle\sum_{j=1}^{32}\|\widetilde\Phi_j\|_{Z_0^{even}}\|\widetilde q_j\|_{(Z')^{even}}\leq146\\[8pt]
\hline
\end{array}
\]
\end{table}

\begin{proposition}\label{close:proposition}
Equation~\eqref{cl2} has a
unique solution in the ball $\{f\in Z^{odd}:\|f\|_{Z^{odd}}\leq10^{-9}\}$,
and this solution satisfies $\|f\|_{Z_0^{odd}}\leq10^{-9}$. For this $f$, equation
\eqref{cl2.5} has a unique solution in the ball
\begin{equation}\label{close:radii}
\widetilde f\in Z_{0,\ast}^{even},\qquad
\|\widetilde f\|_{Z_0^{even}}\leq3\times10^{-7},\qquad
|\delta\lambda|\leq6\times10^{-6}.
\end{equation}
The equations hold in the sense of distributions. Moreover,
$\widetilde\lambda_{app}+\delta\lambda\in[0.23,0.24]$ and
$\|\wOmega^\ast+\widetilde f\|_{L^2}\geq1-3\times10^{-7}>0$.
\end{proposition}

\begin{proof}
\textbf{Step 1.} We first bound $K$, $\wOmega^\ast$, and the scalar $\mu$. For an axi-symmetric potential $g\in\dot H^1(\R^5)$, integration by parts
in $r$ gives
\begin{equation}\label{close:shifted-gradient}
\big\|(\partial_rg+2g/r,\partial_zg)\big\|_{L^2}
=\|\nabla g\|_{L^2}.
\end{equation}
Pairing the two terms of $K$ with a test function and using this identity,
the drift bound in Table~\ref{close:profiles}, and \eqref{nl4}, we obtain
\[
\|Kh\|_{Z'}\leq
\big(22.45+C_b\|\mathcal R_\Psi\|_Z\big)\|h\|_{L^2}
+13\||\nabla|^{-1}h\|_{L^2}.
\]
Here the first term uses $\Psi_5$ in place of $\Psi^\ast$ and includes
the difference $-\mathcal N_5(\Delta^{-1}_{\R^5}\mathcal R_\Psi,h)$. Thus
\begin{equation}\label{close:derived-profile-bounds}
\begin{split}
&\|K\|_{Z\to Z'}\leq
\big[(22.45+0.08\cdot1.251\times10^{-9})^2+13^2\big]^{1/2}<26,\\
&\|\wOmega^\ast\|_{Z^{even}}
\leq(1.01^2+0.89^2)^{1/2}<1.35,\qquad
\|S\wOmega^\ast\|_{(Z')^{even}}\leq22.55+5\cdot0.89=27.
\end{split}
\end{equation}
For the last inequality use
$S\wOmega^\ast=\operatorname{div}_{\R^5}((x,z)\wOmega^\ast)-5\wOmega^\ast$
and $\|\wOmega^\ast\|_{Z'}\leq\|\nabla\wPsi^\ast\|_{L^2}$.
Corollary~\ref{nl:linear} now gives
\begin{equation}\label{close:linear-forcing}
\|\wt{\mathcal L}_2 f+\wt{\mathcal L}_3(\Delta^{-1}_{\R^5}f)\|_{(Z')^{even}}
\leq0.08(1.35+1.01)\|f\|_{Z^{odd}}
=0.1888\|f\|_{Z^{odd}}.
\end{equation}

If $h\in Z_{0,\ast}^{even}$ and $g=\Delta^{-1}_{\R^5}h$, the normalization
implies $\langle h,\wPsi^\ast\rangle=\langle g,\wOmega^\ast\rangle=0$.
Integration by parts then yields
\begin{equation}\label{close:scalar-identity}
\left\langle\mathcal S_{\widetilde\lambda_{app}}h,\wPsi^\ast\right\rangle
=\langle h,a_\ast\rangle
+\left\langle(\partial_rg+2g/r,\partial_zg),b_\ast\right\rangle.
\end{equation}
For the scaling term we use $S\Delta_{\R^5} g=\Delta_{\R^5}(Sg-2g)$ together with $\langle g,\wOmega^\ast\rangle=0$, while the terms of $K$ follow directly from integration by parts. The identity extends by density using the displayed $L^2$ bounds.
Combining \eqref{close:shifted-gradient}--\eqref{close:scalar-identity}
with \eqref{robu2.5}, we find that $T(\mu,h)=F$ implies
\begin{equation}\label{close:scalar-bound}
\begin{split}
|\mu|&\leq
\frac{1.35}{0.88^2}\|F\|_{(Z')^{even}}
+\frac{14.4}{0.88^2}\|h\|_{Z^{even}}\leq1.75\|F\|_{(Z')^{even}}+18.6\|h\|_{Z^{even}}.
\end{split}
\end{equation}

\textbf{Step 2.} We prove the required inverse bounds.
By Table~\ref{close:responses}, the second conditions in \eqref{robu11.1}
and \eqref{robu37} hold, since $0.01\gamma=10^{-3}$. The choice $\delta=10^{-6}$ satisfies
the first conditions in \eqref{robu11.1} and \eqref{robu37}. Indeed, in the
more restrictive even condition the right-hand side is at least
\[
\frac{0.01}{6+2\cdot26+2\cdot61}
\left(1+\frac{1.35+27}{0.88}\right)^{-1}
>1.67\times10^{-6}.
\]
Thus Lemmas~\ref{weakCoercive1} and~\ref{weakCoercive2} give the two
isomorphisms. We use a $Z$ estimate for the odd inverse. In the notation
of the proof of Lemma~\ref{weakCoercive1}, pairing
$(\mathcal S_++E)g=F$ with $\Delta^{-1}_{\R^5}g$ gives
\[
(0.1-5\times10^{-4})\|g\|_{Z^{odd}}\leq\|F\|_{(Z')^{odd}}.
\]
Since $\mathcal S_0^{-1}F=(I+G)g$, Table~\ref{close:responses} yields
\begin{equation}\label{close:inverse-bounds}
\begin{split}
&\|\mathcal S_0^{-1}F\|_{Z^{odd}}
\leq\frac{1+15.25}{0.1-5\times10^{-4}}\|F\|_{(Z')^{odd}}
\leq164\|F\|_{(Z')^{odd}},\\
&\|h\|_{Z_0^{even}}\leq\frac{1+146}{0.98\cdot0.1}\|F\|_{(Z')^{even}}
=1500\|F\|_{(Z')^{even}}\qquad\text{if }T(\mu,h)=F.
\end{split}
\end{equation}
The first estimate uses the $Z$ norms of the odd responses, whereas the
second uses the full $Z_0$ norms of the even responses. For
$f=\mathcal S_0^{-1}F$, the equation gives
$Sf=2F-2\Delta_{\R^5}f-3f-2Kf$, hence
\[
\|Sf\|_{(Z')^{odd}}\leq2\|F\|_{(Z')^{odd}}+57\|f\|_{Z^{odd}}.
\]
Retaining the unrounded constant in \eqref{close:inverse-bounds}, we obtain
\[
\|\mathcal S_0^{-1}F\|_{Z_0^{odd}}
\leq\left[(1+57\delta)\frac{16.25}{0.0995}+2\delta\right]
\|F\|_{(Z')^{odd}}\leq163.326\|F\|_{(Z')^{odd}},
\qquad \delta=10^{-6}.
\]
In particular,
\eqref{close:scalar-bound} and \eqref{close:inverse-bounds} imply
$|\mu|\leq28000\|F\|_{(Z')^{even}}$.

\textbf{Step 3.} We solve \eqref{cl2} in the $Z_0^{odd}$ ball of radius $10^{-9}$.
Since $\|f\|_{Z^{odd}}\leq\|f\|_{Z_0^{odd}}$, the map
\[
f\longmapsto\mathcal S_0^{-1}
\left[\mathcal R-\mathcal N_5(\Delta^{-1}_{\R^5}f,f)\right]
\]
has norm at most
\[
164\big(6\times10^{-12}+0.08(10^{-9})^2\big)=9.8400001312\times10^{-10}<10^{-9}
\]
on this ball, and Lipschitz constant at most
$2\cdot0.08\cdot164\cdot10^{-9}=2.624\times10^{-8}$.
The contraction mapping theorem gives a unique fixed point in this ball. Any distributional
solution in the $Z^{odd}$ ball of radius $10^{-9}$ also has $Sf\in(Z')^{odd}$ by \eqref{cl2},
so it belongs to $Z_0^{odd}$. By the preceding inverse estimate, it lies in
the same $Z_0^{odd}$ ball where it is a fixed point of the same map.

\textbf{Step 4.} It remains to solve \eqref{cl2.5} in the ball \eqref{close:radii}.
Fix the solution $f$ of \eqref{cl2} constructed above. By \eqref{close:linear-forcing},
\begin{equation}\label{close:effective-residual}
\big\|\wt{\mathcal R}-\wt{\mathcal L}_2f
-\wt{\mathcal L}_3(\Delta^{-1}_{\R^5}f)\big\|_{(Z')^{even}}
\leq1.5\times10^{-12}+0.1888\cdot10^{-9}<1.91\times10^{-10}.
\end{equation}
Given $(\mu,h)\in\R\times Z_{0,\ast}^{even}$ in this ball, define
$(\mu^+,h^+)=T^{-1}F(\mu,h)$, where
\[
F(\mu,h):=\wt{\mathcal R}-\wt{\mathcal L}_2f
-\wt{\mathcal L}_3(\Delta^{-1}_{\R^5}f)
-\mathcal N_5(\Delta^{-1}_{\R^5}f,h)-\mathcal N_5(\Delta^{-1}_{\R^5}h,f)+\mu h.
\]
Here $\|h\|_{Z'}\leq\|h\|_Z$. The parity of each term is even, and the
normalization of $h^+$ is imposed by $T^{-1}$. Hence
\[
\|F(\mu,h)\|_{Z'}
\leq1.91\times10^{-10}
+(0.16\cdot10^{-9}+6\times10^{-6})\,3\times10^{-7}
=1.92800048\times10^{-10}.
\]
Equations \eqref{close:scalar-bound} and \eqref{close:inverse-bounds} give
\[
\|h^+\|_{Z_0^{even}}\leq2.89200072\times10^{-7}<3\times10^{-7},
\qquad |\mu^+|<5.4\times10^{-6}<6\times10^{-6}.
\]
For $(\mu_1,h_1),(\mu_2,h_2)$ in the ball with distance
\[
D:=\max\left\{\frac{\|h_1-h_2\|_{Z_0^{even}}}{3\times10^{-7}},\,
\frac{|\mu_1-\mu_2|}{6\times10^{-6}}\right\},
\]
the forcing difference is bounded by
\[
\|F(\mu_1,h_1)-F(\mu_2,h_2)\|_{Z'}
\leq3\times10^{-7}(0.16\cdot10^{-9}+12\times10^{-6})D
=3.600048\times10^{-12}D.
\]
By the inverse bounds, the distance between the images is at most
\[
3.600048\times10^{-12}
\max\left\{\frac{1500}{3\times10^{-7}},\,
\frac{28000}{6\times10^{-6}}\right\}D
=0.01800024D<0.019D.
\]
By the contraction mapping theorem, \eqref{cl2.5} has a unique solution in the ball \eqref{close:radii}. Finally, the eigenvalue interval claimed in the proposition follows from the fixed value \eqref{lambda:app} and $|\delta\lambda|\leq6\times10^{-6}$, and the
nontriviality follows from $\|\widetilde f\|_{L^2}\leq3\times10^{-7}$.
\end{proof}

\subsection{Completion of the proofs of Proposition~\ref{PropCentral} and Theorem~\ref{Main3D}}\label{finish}
We reconstruct the profiles using \eqref{cl4}, and set $R=(r^2+z^2)^{1/2}$. In this subsection $y\in\R^3$, and $U^r,U^z$ denote cylindrical velocity components. Their three-dimensional velocities and homogeneous limiting field are
\begin{equation}\label{fin1}
\begin{split}
U&=r\partial_z\Psi e_r-(r\partial_r\Psi+2\Psi)e_z,\qquad
\widetilde U=r\partial_z\wPsi e_r-(r\partial_r\wPsi+2\wPsi)e_z,\\
U_0&=\sigma R^{-4}\big[r(r^2-z^2)e_r-2z^3e_z\big].
\end{split}
\end{equation}
Define $U_{app}$ by the same formula as $U$, with $\Psi$ replaced by $\Psi_{app}$.

\begin{lemma}\label{finlem}
The reconstructed profiles are smooth. For every multi-index $\alpha$,
\begin{equation}\label{fin2}
|\partial^\alpha(U-U_0)(y)|\lesssim_\alpha |y|^{-3-|\alpha|}
\quad(|y|\geq1),\qquad
|\partial^\alpha\widetilde U(y)|\lesssim_\alpha\langle y\rangle^{-4-|\alpha|}.
\end{equation}
\end{lemma}
\begin{proof}
The identities \eqref{nl11} and the explicit asymptotics give
$H:=U-e^{\Delta_{\R^3}} U_0\in H^1(\R^3)$ and $\widetilde U\in H^1(\R^3)$. Indeed,
\begin{align*}
    \|U-U_{app}\|_{H^1(\R^3)}^2=\pi^{-1}\|f-\mathcal R_\Psi\|_Z^2,\qquad \|\widetilde U\|_{H^1(\R^3)}^2=\pi^{-1}\|\wOmega^\ast+\widetilde f\|_Z^2.
\end{align*}
In particular, $U\in L^6(\R^3)$, $\nabla U\in L^2(\R^3)$. The scalar equations imply
\begin{equation}\label{fin3}
\begin{split}
\Delta_{\R^3} U+\tfrac12y\cdot\nabla U+\tfrac12U
&=\Pi\operatorname{div}_{\R^3}(U\otimes U),\\
\Delta_{\R^3}\widetilde U+\tfrac12y\cdot\nabla\widetilde U+
(\tfrac12-\widetilde\lambda)\widetilde U
&=\Pi\operatorname{div}_{\R^3}(U\otimes\widetilde U+\widetilde U\otimes U),
\end{split}
\end{equation}
where $\Pi$ is the Leray projection on $\R^3$. Here the $L^6(\R^3)$ bounds exclude an additional harmonic polynomial. Local Stokes estimates first give $U\in W^{2,3/2}_{loc}(\R^3)$ and then, by bootstrapping, $U,\widetilde U\in C^\infty$. The axis-regular identity
$\Psi(r,z)=-\int_0^1sU^z(sr,z)\,ds$, and its tilded version, give smooth scalar profiles.

The self-similar solution $u(t)=e^{t\Delta_{\R^3}}U_0+t^{-1/2}H(\cdot/\sqrt t)$ has initial trace $U_0$ in $L^2_{loc}(\R^3)$, since the second term has $L^2(\R^3)$ norm $t^{1/4}\|H\|_{L^2(\R^3)}$. Also $\|\nabla u(t)\|_{L^2(\R^3)}=t^{-1/4}\|\nabla U\|_{L^2(\R^3)}$. Together with
$|e^{t\Delta_{\R^3}}U_0(y)|\lesssim(\sqrt t+|y|)^{-1}$, these bounds give the uniform local energy bounds and spatial tail condition of \cite[Definition 3.1]{JiSv1}. Smoothness for $t>0$ gives the local energy identity. Thus \cite[Theorem 4.1]{JiSv1} applies and yields
\begin{align*}
    |\partial^\alpha(U-e^{\Delta_{\R^3}} U_0)|\lesssim_\alpha\langle y\rangle^{-3-|\alpha|}.
\end{align*}
The corresponding heat-kernel estimate for $e^{\Delta_{\R^3}} U_0-U_0$ proves the first bound in \eqref{fin2}.

Next, put $F=U\otimes\widetilde U+\widetilde U\otimes U\in L^2(\R^3)$.
The field $t^{\widetilde\lambda-1/2}\widetilde U(\cdot/\sqrt t)$ solves the forced heat equation with tensor forcing $-t^{\widetilde\lambda-1}F(\cdot/\sqrt t)$ and has $L^2(\R^3)$ norm $t^{\widetilde\lambda+1/4}\|\widetilde U\|_{L^2(\R^3)}\to0$. Hence
\begin{equation}\label{fin4}
\widetilde U=-\mathcal T_{\widetilde\lambda}F,\qquad
\mathcal T_\lambda F:=\int_0^1s^{\lambda-1}
e^{(1-s)\Delta_{\R^3}}\Pi\operatorname{div}_{\R^3}[F(\cdot/\sqrt s)]\,ds.
\end{equation}
For $\lambda>0$, $\mathcal T_\lambda$ is bounded on $L^2(\R^3)$, with norm bounded by
$C\int_0^1s^{\lambda-1/4}(1-s)^{-1/2}\,ds$. It is also bounded on the weighted spaces
\begin{equation}\label{fin5}
\|h\|_{X_m}:=\sum_{|\alpha|\leq m}
\sup_{y\in\R^3}\langle y\rangle^{4+|\alpha|}|\partial^\alpha h(y)|,\qquad
\|\mathcal T_\lambda F\|_{X_m}\lesssim_{\lambda,m}\|F\|_{X_m}.
\end{equation}
Indeed, the kernel $K_t$ of $e^{t\Delta_{\R^3}}\Pi$ satisfies
$|\partial^\alpha\nabla K_t(y)|\lesssim_\alpha(|y|+\sqrt t)^{-4-|\alpha|}$ and
$\|\nabla K_t\|_{L^1(\R^3)}\lesssim t^{-1/2}$. Fix $|y|=R\geq2$ and split the convolution near $Y=y$, using a cutoff at scale $R/4$. In the near part move derivatives onto $F(Y/\sqrt s)$ and the cutoff; in the far part use $\|F(\cdot/\sqrt s)\|_{L^1(\R^3)}=s^{3/2}\|F\|_{L^1(\R^3)}$. The resulting bound is
$C_m\|F\|_{X_m}R^{-4-|\alpha|}[s^2(1-s)^{-1/2}+s^{3/2}]$, integrable against $s^{\lambda-1}ds$.
For bounded $y$, split at $s=1/2$, placing derivatives on the kernel for small $s$ and on $F$ otherwise. This proves \eqref{fin5}.

Finally, we write $U=U_{\rm in}+U_{\rm out}$ using a smooth cutoff at radius $L$. From above we have
\[
\max_{|\alpha|\leq m}\sup_{y\in\R^3}\langle y\rangle^{|\alpha|}
|\partial^\alpha U_{\rm out}(y)|\lesssim_m L^{-1}.
\]
Thus $\mathcal B_Lh:=-\mathcal T_{\widetilde\lambda}(U_{\rm out}\otimes h+h\otimes U_{\rm out})$ has norm $O_m(L^{-1})$ on both $L^2(\R^3)$ and $X_m$. Choose $L$ so both norms are below $1/2$. Equation \eqref{fin4} becomes
$\widetilde U=G_L+\mathcal B_L\widetilde U$, where
$G_L=-\mathcal T_{\widetilde\lambda}(U_{\rm in}\otimes\widetilde U+\widetilde U\otimes U_{\rm in})\in X_m$,
since its tensor source is smooth and compactly supported.
The Neumann series gives a solution in $X_m\subset L^2(\R^3)$; uniqueness in $L^2(\R^3)$ identifies it with $\widetilde U$. Since $m$ is arbitrary, this proves the second bound in \eqref{fin2}.
\end{proof}

To finish Proposition~\ref{PropCentral}, we use $\partial_z\Psi=U^r/r$ and integrate from $z$ to the infinity of the same sign. Near the axis use $U^r/r=\int_0^1\partial_rU^r(sr,z)\,ds$, and away from it divide directly by $r$. Applied to $U-U_0$ and $\widetilde U$, \eqref{fin2} gives
$|\partial^\beta(\Psi-\sigma z/R^2)|\lesssim_\beta R^{-3-|\beta|}$ for $R\geq1$ and
$|\partial^\beta\wPsi|\lesssim_\beta\langle R\rangle^{-4-|\beta|}$.
Applying $L_0$ yields \eqref{sf2}. Reality, parity, the eigenvalue interval, and nontriviality follow from Proposition~\ref{close:proposition}.

Finally, $U_0$ is divergence-free, smooth off the origin, and homogeneous of degree $-1$.
Equations \eqref{fin3} and \eqref{fin2} verify the hypotheses of \cite[Theorem 1.2]{IoJiPa}; the velocity eigenfunction is nonzero because
$\operatorname{curl}\widetilde U=r\wOmega e_\theta\not\equiv0$.
That theorem gives all the conclusions of Theorem~\ref{Main3D}, including the energy identity and preservation of axial symmetry and absence of swirl.

\section{The approximate solutions $(\Psi_{app},\Omega_{app})$ and $(\widetilde{\lambda}_{app},\widetilde{\Psi}_{app},\widetilde{\Omega}_{app})$}

We now describe in detail the construction of the approximate solutions in Proposition~\ref{PropCentral2}.

\subsection{Asymptotic expansions}\label{fixed_asymp} Recall the nonlinear system we are looking to solve
\begin{equation}\label{gf1}
L_0\Psi-\Omega=0,\qquad L_0\Omega+\frac{1}{2}\big(r\partial_r+z\partial_z+3\big)\Omega-r\partial_z\Psi\partial_r\Omega+r\partial_r\Psi\partial_z\Omega+2\Psi\partial_z\Omega=0
\end{equation}
in the space of functions that are odd in $z$, where $L_0=\partial_r^2+\partial_z^2+(3/r)\partial_r$.
Throughout this subsection we write
\begin{equation}\label{gf0}
R:=\big(r^2+z^2\big)^{1/2},\qquad \ell:=\log R.
\end{equation}

We start with the top order terms, homogeneous of orders $-1$ for $\Psi$ and $-3$ for $\Omega$,
\begin{equation}\label{gf2}
\Psi_{-1}(r,z):=\frac{\sigma z}{R^{2}},\qquad \Omega_{-3}(r,z):=\frac{-6\sigma z}{R^{4}}.
\end{equation}
The first equation in \eqref{gf1} is satisfied identically for these functions. The second equation can be used to determine the homogeneous term $\Omega_{-5}$, by the formula
\begin{equation}\label{gf3}
L_0\Omega_{-3}+\frac{1}{2}\big(r\partial_r+z\partial_z+3\big)\Omega_{-5}-r\partial_z\Psi_{-1}\partial_r\Omega_{-3}+(r\partial_r\Psi_{-1}+2\Psi_{-1})\partial_z\Omega_{-3}=0,
\end{equation}
which leads to
\begin{equation}\label{gf4}
\Omega_{-5}(r,z):=\frac{12\sigma^{2}z}{R^{8}}\Big[(-2r^{2}+3z^{2})+2\sigma^{-1}R^{2}\Big].
\end{equation}
We can then use the first equation in \eqref{gf1} to find the term of order $-3$ in $\Psi$,
\begin{equation}\label{gf5}
\Psi_{-3}(r,z):=\frac{\sigma^{2}z}{3R^{6}}\Big[(3r^{2}-7z^{2})-18\sigma^{-1}R^{2}\Big]
+\frac{\sigma'}{R^{3}},
\end{equation}
where $\sigma'\in\R$. Since our base solution is odd in $z$ we must have $\sigma'=0$, therefore
\begin{equation}\label{gf6}
\begin{split}
&\Psi_{\geq -3}(r,z)=\frac{\sigma z}{R^{2}}
+\frac{\sigma^{2}z}{3R^{6}}\Big[(3r^{2}-7z^{2})-18\sigma^{-1}R^{2}\Big],\\
&\Omega_{\geq -5}(r,z)=\frac{-6\sigma z}{R^{4}}
+\frac{12\sigma^{2}z}{R^{8}}\Big[(-2r^{2}+3z^{2})+2\sigma^{-1}R^{2}\Big].
\end{split}
\end{equation}
This gives us the boundary data we use to construct our solutions in \cite{IoJiPa}.

To get good residuals in the asymptotic region we need a more precise expansion. Some care is needed, because the equation $L_0(\Psi)=0$ admits general homogeneous solutions of degree $-4$ at infinity of the form
\begin{equation}\label{gf7}
\Psi_{-4}(r,z):=\frac{\sigma_2 z}{R^{5}},\qquad\sigma_2\in\R.
\end{equation}
These solutions are odd in $z$, thus acceptable. The choice of the parameter $\sigma_2$ is determined by global smoothness considerations and cannot be found just by asymptotic analysis.

To continue the expansion, the next step is to calculate explicitly the correction $\Omega_{-7}$. The relevant equation is analogous to equation \eqref{gf3}
\begin{equation}\label{gf8}
\begin{split}
L_0\Omega_{-5}+\frac{1}{2}\big(r\partial_r+z\partial_z+3\big)\Omega_{-7}&-r\partial_z\Psi_{-1}\partial_r\Omega_{-5}+(r\partial_r\Psi_{-1}+2\Psi_{-1})\partial_z\Omega_{-5}\\
&-r\partial_z\Psi_{-3}\partial_r\Omega_{-3}+(r\partial_r\Psi_{-3}+2\Psi_{-3})\partial_z\Omega_{-3}=0.
\end{split}
\end{equation}
This gives
\begin{equation}\label{gf9}
\Omega_{-7}(r,z):=\frac{2\sigma^{3}z}{R^{12}}\Big[(-39r^{4}+190r^{2}z^{2}-111z^{4})
+36\sigma^{-1}R^{2}(3r^{2}-5z^{2})+36\sigma^{-2}R^{4}\Big],
\end{equation}
and then
\begin{equation}\label{gf10}
\Psi_{-5}(r,z):=\frac{\sigma^{3}z}{3R^{10}}\Big[(-8r^{4}-47r^{2}z^{2}+29z^{4})
-108\sigma^{-1}R^{2}(r^{2}-z^{2})+36\sigma^{-2}R^{4}\Big].
\end{equation}

We can now calculate the next correction term $\Omega_{-8}$, similarly, leading to
\begin{equation}\label{gf12}
\Omega_{-8}(r,z):=\frac{12\sigma\sigma_2 z}{5R^{11}}\Big[-r^{2}+6z^{2}\Big].
\end{equation}
Then we calculate the $\Psi_{-6}$ correction term,
\begin{equation}\label{gf13}
\Psi_{-6}(r,z):=\frac{\sigma\sigma_2 z}{105R^{9}}
\Big[\big(28\ell+105\,\sigma^{-1}\sigma_2^{-1}\sigma_3\big)(3r^{2}-4z^{2})+36R^{2}\Big],
\end{equation}
where $\sigma_3$ is another independent constant.

Then we calculate the correction $\Omega_{-9}$, as in \eqref{gf8},
\begin{equation}\label{gf15}
\begin{split}
\Omega_{-9}(r,z)&:=\frac{8\sigma^{4}z}{R^{16}}\big[
(-30r^{6}+282r^{4}z^{2}-531r^{2}z^{4}+179z^{6})\\
&+2\sigma^{-1}R^{2}(78r^{4}-403r^{2}z^{2}+215z^{4})+12\sigma^{-2}R^{4}(6r^{2}-11z^{2})+72\sigma^{-3}R^{6}\big],
\end{split}
\end{equation}
and the correction term $\Psi_{-7}$ by solving $L_0\Psi_{-7}=\Omega_{-9}$, leading to
\begin{equation}\label{gf17}
\begin{split}
\Psi_{-7}&(r,z):=\frac{\sigma^{4}z}{15R^{14}}\big[
(-451r^{6}-149r^{4}z^{2}+2685r^{2}z^{4}-537z^{6})\\
&+6\sigma^{-1}R^{2}(-471r^{4}+1462r^{2}z^{2}-387z^{4})+36\sigma^{-2}R^{4}(27r^{2}-41z^{2})+360\sigma^{-3}R^{6}\big].
\end{split}
\end{equation}

To continue the expansion, we calculate the correction $\Omega_{-10}$ in the same way as before,
\begin{equation}\label{gf19}
\begin{split}
\Omega_{-10}(r,z):=\frac{4\sigma^{2}\sigma_2 z}{245R^{15}}\big[
&\big(28\ell+105\,\sigma^{-1}\sigma_2^{-1}\sigma_3\big)(3r^{4}+56r^{2}z^{2}-24z^{4})\\
&-6(279r^{4}-1001r^{2}z^{2}+1702z^{4})
+504\sigma^{-1}R^{2}(3r^{2}+5z^{2})\big].
\end{split}
\end{equation}
Then we calculate the $\Psi_{-8}$ term as a solution of \(L_0\Psi_{-8}=\Omega_{-10}\). The result is
\begin{equation}\label{gf20}
\begin{split}
\Psi_{-8}(r,z)&:=\frac{\sigma^{2}\sigma_2 z}{R^{13}}\Big\{
R^{4}\Big[\frac{512}{6615}\ell-\frac{352204}{416745}
+\frac{216}{245}\sigma^{-1}+\frac{128}{441}\sigma^{-1}\sigma_2^{-1}\sigma_3\Big]\\
&+R^{2}(3r^{2}-4z^{2})\Big[\frac{32}{539}\ell+\frac{159332}{207515}
-\frac{288}{2695}\sigma^{-1}+\frac{120}{539}\sigma^{-1}\sigma_2^{-1}\sigma_3\Big]\\
&+(5r^{4}-20r^{2}z^{2}+8z^{4})\Big[\frac{8}{195}\ell^{2}+\frac{134144}{195195}\ell
+\frac{4}{13}\sigma^{-1}\sigma_2^{-1}\sigma_3\,\ell
+\sigma^{-2}\sigma_2^{-1}\sigma_4\Big]\Big\},
\end{split}
\end{equation}
where \(\sigma_4\in\mathbb R\) is a new constant that multiplies the degree $-8$ homogeneous solution $z(5r^4-20r^2z^2+8z^4)R^{-13}$.  

Finally, we calculate \(\Omega_{-11}\) using the degree \(-11\) equation in the same way as before,
\begin{equation}\label{gf22}
\begin{split}
\Omega_{-11}&(r,z):=\frac{2\sigma^{5}z}{5R^{20}}\big[
(-1938r^{8}+24458r^{6}z^{2}-99933r^{4}z^{4}+104448r^{2}z^{6}-23623z^{8})\\
&+4\sigma^{-1}R^{2}
(3403r^{6}-39831r^{4}z^{2}+62949r^{2}z^{4}-17177z^{6})\\
&+12\sigma^{-2}R^{4}
(293r^{4}-2764r^{2}z^{2}+2083z^{4})+14400\sigma^{-3}R^{6}(r^{2}-2z^{2})
+18000\sigma^{-4}R^{8}\big].
\end{split}
\end{equation}

Solving \(L_0\Psi_{-9}=\Omega_{-11}\), we obtain
\begin{equation}\label{gf23}
\begin{split}
\Psi_{-9}&(r,z):=\frac{\sigma^{5}z}{7875R^{18}}\big[
(-1413624r^{8}+5108254r^{6}z^{2}+9506631r^{4}z^{4}
-11069496r^{2}z^{6}+785751z^{8})\\
&+36\sigma^{-1}R^{2}
(-202799r^{6}+1280103r^{4}z^{2}-1106397r^{2}z^{4}+109201z^{6})+36\sigma^{-2}R^{4}\\
&\times(43578r^{4}-224719r^{2}z^{2}+117203z^{4})+151200\sigma^{-3}R^{6}(9r^{2}-16z^{2})
+1134000\sigma^{-4}R^{8}\big].
\end{split}
\end{equation}
No new independent constants occur in \(\Psi_{-9}\), since the axially
symmetric homogeneous solution of \(L_0f=0\) of degree \(-9\) is even
in \(z\), and is therefore excluded by the odd symmetry. 

We also need the first terms in the asymptotic expansion of the eigenfunction. The eigenpair $(\wt\lambda,\wPsi,\wOmega)$ solves the linear system \eqref{sf1.2}, which we write in the form
\begin{equation}\label{gf24}
\begin{split}
&L_0\wPsi-\wOmega=0,\qquad L_0\wOmega+\frac{1}{2}\big(r\partial_r+z\partial_z+3-2\wt\lambda\big)\wOmega-r\partial_z\Psi\partial_r\wOmega+(r\partial_r\Psi+2\Psi)\partial_z\wOmega\\
&\qquad\qquad-r\partial_z\wPsi\partial_r\Omega+(r\partial_r\wPsi+2\wPsi)\partial_z\Omega=0,
\end{split}
\end{equation}
in the space of functions that are even in $z$. The decaying axi-symmetric homogeneous solutions of $L_0f=0$ that are even in $z$ are the functions $h_{2k}R^{-4k-3}$, $k\geq0$, where $h_{2k}$ is the homogeneous axi-symmetric harmonic polynomial of degree $2k$ in $\R^5$ that is even in $z$. We will use
\begin{equation}\label{gf25}
h_2(r,z):=4z^2-r^2,\qquad h_4(r,z):=r^4-12r^2z^2+8z^4.
\end{equation}

The slowest decaying term allowed in the expansion of $\wPsi$ is the homogeneous solution $\tau_1R^{-3}$ ($k=0$), in which case $\wOmega=O(R^{-7})$ and $\tau_1=-\frac14\int_{(0,\infty)\times\R}\wOmega\,d\mu$, where $d\mu:=r^3\,dr\,dz$ on $(0,\infty)\times\R$. This coefficient vanishes. Indeed, we integrate the second equation in \eqref{gf24} with respect to $d\mu$. All the boundary terms vanish because of the decay, and $\int_{(0,\infty)\times\R}(r\partial_r+z\partial_z)g\,d\mu=-5\int_{(0,\infty)\times\R} g\,d\mu$, $\int_{(0,\infty)\times\R}\big(-r\partial_zf\,\partial_rg+(r\partial_rf+2f)\partial_zg\big)\,d\mu=2\int_{(0,\infty)\times\R} g\,\partial_zf\,d\mu$. Therefore
\[
(1+\wt\lambda)\int_{(0,\infty)\times\R}\wOmega\,d\mu=2\int_{(0,\infty)\times\R}\big(L_0\wPsi\,\partial_z\Psi+L_0\Psi\,\partial_z\wPsi\big)\,d\mu=0,
\]
since $L_0$ is symmetric with respect to $d\mu$ and commutes with $\partial_z$, while $\partial_z$ is antisymmetric. Thus $\tau_1=0$. The term $\wOmega_{-7}$ is a multiple of $\sigma\tau_1$, so $\wOmega_{-7}=0$ and the top order term of $\wPsi$ is the homogeneous solution
\begin{equation}\label{gf26}
\wPsi_{-5}(r,z):=\frac{\tau_2h_2}{R^{7}}=\frac{\tau_2(4z^2-r^2)}{R^{7}},\qquad\tau_2\in\R.
\end{equation}
The correction $\wOmega_{-9}$ is determined by the terms of degree $-9$ in the second equation in \eqref{gf24}, which reduce to
\begin{equation}\label{gf27}
\frac{1}{2}\big(r\partial_r+z\partial_z+3-2\wt\lambda\big)\wOmega_{-9}-r\partial_z\wPsi_{-5}\partial_r\Omega_{-3}+(r\partial_r\wPsi_{-5}+2\wPsi_{-5})\partial_z\Omega_{-3}=0.
\end{equation}
Since $\frac12(r\partial_r+z\partial_z+3-2\wt\lambda)=-(\wt\lambda+3)$ on functions homogeneous of degree $-9$, we can take
\begin{equation}\label{gf28}
\wOmega_{-9}(r,z):=-\frac{18\sigma\tau_2}{(\wt\lambda_a+3)R^{13}}\big(r^4+8r^2z^2-8z^4\big),
\end{equation}
where $\wt\lambda_a$ is a good approximation of the eigenvalue. Then we solve $L_0\wPsi_{-7}=\wOmega_{-9}$. The source $\wOmega_{-9}$ has a nontrivial component along $h_4$, which corresponds to the homogeneous solution $h_4R^{-11}$ of the same degree $-7$, so a logarithmic term appears:
\begin{equation}\label{gf29}
\wPsi_{-7}(r,z):=-\frac{\sigma\tau_2}{(\wt\lambda_a+3)R^{11}}\Big[\frac{4}{49}\big(17r^4-15r^2z^2-32z^4\big)+\frac{90}{77}\,\ell\,h_4\Big]+\frac{\tau_3h_4}{R^{11}},
\end{equation}
where $\tau_3\in\R$ is a new constant. As for $\sigma_2,\sigma_3,\sigma_4$, the constants $\tau_2,\tau_3$ are not determined by the asymptotic analysis.

\subsubsection{The full asymptotic solution}\label{ss:trunc}
This procedure can be continued to find more terms in the asymptotic expansion of the solution, but for our purpose we can stop here. Since $L_0\Psi_{-1}=\Omega_{-3}$, $L_0\Psi_{-3}=\Omega_{-5}$, $L_0\Psi_{-4}=0$, and
$L_0\Psi_{-k}=\Omega_{-(k+2)}$ for $k=5,\ldots,9$, the pair
\begin{equation}\label{e:trnl}
\Psi_{\geq-9}:=\sum_{k\in\{1,3,4,5,6,7,8,9\}}\Psi_{-k},
\qquad
\Omega_{\geq-11}:=\sum_{k\in\{3,5,7,8,9,10,11\}}\Omega_{-k}
\end{equation}
satisfies $L_0\Psi_{\geq-9}-\Omega_{\geq-11}=0$ identically, and leaves in
the second equation of \eqref{gf1} a residual $O(R^{-12}\ell^{2})$. We fix the constants
\begin{equation}\label{constants}
\begin{split}
\sigma:=130, \qquad\sigma_2:=288016.2,\qquad\sigma_3:= 90012367,\qquad\sigma_4:= 3.2103\times 10^{10}.
\end{split}
\end{equation}

For the eigenfunction we stop at the pair
\begin{equation}\label{e:trnl2}
\wPsi_{\geq-7}:=\wPsi_{-5}+\wPsi_{-7},\qquad\wOmega_{\geq-9}:=\wOmega_{-9}.
\end{equation}
Since $L_0\wPsi_{-5}=0$ and $L_0\wPsi_{-7}=\wOmega_{-9}$, it satisfies $L_0\wPsi_{\geq-7}-\wOmega_{\geq-9}=0$ identically. We set $(\wt\Psi_{app},\wt\Omega_{app})=(\wt\Psi_{\geq-7},\wt\Omega_{\geq-9})$ in the asymptotic region, and set the constants
\begin{equation}\label{constants2}
\wt\lambda_a:=0.235059686,\qquad\tau_2:=227.94186,\qquad\tau_3:=-152371.97.
\end{equation}

The asymptotic recursion leaves $\sigma_2,\sigma_3,\sigma_4,\tau_2,\tau_3$ free. To choose them we apply Green's formula on large balls to $L_0\Psi=\Omega$ and $L_0\wPsi=\wOmega$, testing against homogeneous axi-symmetric harmonic polynomials of degrees $1,3,5$ and $2,4$, respectively. We evaluate the interior moments from the splines and the exterior contributions from the expansions, subtracting the known divergent terms, including logarithms, before taking finite parts. For the eigenfunction we use the normalization of Proposition~\ref{PropCentral2}. We fix the resulting rounded coefficients in \eqref{constants} and \eqref{constants2}. The residual bounds account for the mismatch where the fields are glued. The constant $\wt\lambda_a$ in \eqref{constants2} approximates $\widetilde{\lambda}_{app}$ only in the exterior expansion.

\subsection{Analysis in the computational domain}\label{sec:computational-domain} 
Following \cite{IoJiPa} we make the change of variables
\begin{equation}\label{u:cov}
u=f(r^2),\qquad v=f(z^2),\qquad f(x)=(1+12x)^{1/4}-1,
\end{equation}
or equivalently
\begin{equation}\label{u:U0}
r^2=U_0(u),\qquad z^2=U_0(v),\qquad
U_0(x):=\frac{(x+1)^4-1}{12}.
\end{equation}

For any functions $g\in C^2_{e,o}$ and $h\in C^2_{e,e}$ we define the functions $g^{(c)},h^{(c)}:[0,\infty)^2\to\C$ by
\begin{equation}\label{transported}
g(r,z)=zg^{(c)}(u,v),\qquad h(r,z)=h^{(c)}(u,v).
\end{equation}
Notice the different definitions depending on the $z$-parity of the function. In particular we define our main variables $\Psi_{app}^{(c)}$, $\Omega_{app}^{(c)}$, $\widetilde{\Psi}_{app}^{(c)}$, $\widetilde{\Omega}_{app}^{(c)}$ in the computational domain by the formulas
\begin{equation}\label{u:norm}
\begin{split}
&\Psi_{app}(r,z)=z\,\Psi_{app}^{(c)}(u,v),\qquad \Omega_{app}(r,z)=z\,\Omega_{app}^{(c)}(u,v),\\
&\wt\Psi_{app}(r,z)=\wt\Psi_{app}^{(c)}(u,v),\qquad \wt\Omega_{app}(r,z)=\wt\Omega_{app}^{(c)}(u,v).
\end{split}
\end{equation}
As explained in \cite{IoJiPa}, these changes of variables serve two main purposes: (1) they encode the desired symmetries of the functions $(\Psi,\Omega, \wPsi,\wOmega)$, and (2) they compress very large physical domains of the form $[0,L']\times[-L',L']$ in $(r,z)$, $L'\approx 3000$, to computational domains of moderate size.

In terms of the new variables the system \eqref{sf3} becomes
\begin{equation}\label{u:nl}
\begin{split}
&\Luv\Psi_{app}^{(c)}-\Omega_{app}^{(c)}=E_{\Psi}^{(c)},\\
&\Luv\Omega_{app}^{(c)}+\big[K(u)\partial_u+K(v)\partial_v+2\big]\Omega_{app}^{(c)}
+2\,\mathcal{N}[\Psi_{app}^{(c)},\Omega_{app}^{(c)}]=E_{\Omega}^{(c)},
\end{split}
\end{equation}
where
\begin{equation}\label{u:K}
K(u)=\frac{(u+1)^4-1}{4(u+1)^3},\qquad K(v)=\frac{(v+1)^4-1}{4(v+1)^3},
\end{equation}
\begin{equation}\label{u:Luv}
\Luv:=\frac{3((u+1)^4-1)}{(u+1)^6}\,\partial_u^2
+\frac{3(5(u+1)^4+3)}{(u+1)^7}\,\partial_u
+\frac{3((v+1)^4-1)}{(v+1)^6}\,\partial_v^2
+\frac{9((v+1)^4+1)}{(v+1)^7}\,\partial_v,
\end{equation}
\begin{equation}\label{u:N}
\mathcal{N}[P,O]:=2K(u)K(v)\big(\partial_uP\,\partial_vO-\partial_vP\,\partial_uO\big)
-K(u)\big(P\partial_uO-\partial_uP\,O\big)+2K(v)P\partial_vO+PO,
\end{equation}
\begin{equation}\label{Emod}
E_\Psi(r,z)=z E_\Psi^{(c)}(u,v),\qquad E_\Omega(r,z)=zE_\Omega^{(c)}(u,v).
\end{equation}
The system \eqref{sf3.2} becomes
\begin{equation}\label{u:eig}
\begin{split}
&\wt{L}_{uv}\wt\Psi_{app}^{(c)}-\wt\Omega_{app}^{(c)}=E_{\wt\Psi}^{(c)},\qquad \wt{L}_{uv}\wt\Omega_{app}^{(c)}+\big[K(u)\partial_u+K(v)\partial_v+\tfrac32-\wt\lambda_{app}\big]\wt\Omega_{app}^{(c)}\\
&\qquad\qquad\qquad+2\,\wt{\mathcal{N}}_1(\Psi_{app}^{(c)},\wt\Omega_{app}^{(c)})+2\,\wt{\mathcal{N}}_2(\wt\Psi_{app}^{(c)},\Omega_{app}^{(c)})=E_{\wt\Omega}^{(c)},
\end{split}
\end{equation}
where
\begin{equation}\label{u:Lpuv}
\wt{L}_{uv}:=\frac{3((u+1)^4-1)}{(u+1)^6}\,\partial_u^2
+\frac{3(5(u+1)^4+3)}{(u+1)^7}\,\partial_u
+\frac{3((v+1)^4-1)}{(v+1)^6}\,\partial_v^2
+\frac{3(3-(v+1)^4)}{(v+1)^7}\,\partial_v,
\end{equation}
\begin{equation}\label{u:halves}
\begin{split}
\wt{\mathcal{N}}_1(P,\wt{O})&:=2K(u)K(v)\big(\partial_uP\partial_v\wt{O}
-\partial_vP\partial_u\wt{O}\big)
-K(u)P\partial_u\wt{O}+2K(v)P\partial_v\wt{O},\\
\wt{\mathcal{N}}_2(\wt P,O)&:=2K(u)K(v)\big(\partial_u\wt{P}\partial_vO
-\partial_v\wt{P}\partial_uO\big)
+K(u)\partial_u\wt{P}O+2K(v)\wt{P}\partial_vO
+\wt{P}O.
\end{split}
\end{equation}
\begin{equation}\label{Emod2}
E_{\wt\Psi}(r,z)=E_{\wt\Psi}^{(c)}(u,v),\qquad E_{\wt\Omega}(r,z)=E_{\wt\Omega}^{(c)}(u,v).
\end{equation}
Notice the slight differences between the operators $L_{uv}$ and $\wt{L}_{uv}$, caused by the parity change and difference in the changes of variables \eqref{transported}. For the same reason, the quadratic forms $\wt{\mathcal{N}}_1$ and $\wt{\mathcal{N}}_2$ are not exactly the formal linearizations of the original quadratic form $\mathcal{N}$.

\subsubsection{Construction of the approximate solutions} We use the computational fields directly and define 
\begin{equation}\label{sec3eq1}
\begin{split}
&\Psi_{ext}(u,v):=\Psi_{\geq -9}(r,z)/z,\qquad \Omega_{ext}(u,v):=\Omega_{\geq -11}(r,z)/z,\\
&\wt\Psi_{ext}(u,v):=\wt\Psi_{\geq -7}(r,z),\qquad \wt\Omega_{ext}(u,v):=\wt\Omega_{\geq -9}(r,z).
\end{split}
\end{equation}
The main step in the construction is to find numerical solutions $(\Psi_{in},\Omega_{in})$ and $(\wt\Psi_{in},\wt{\Omega}_{in})$ of the systems \eqref{u:nl} and \eqref{u:eig} in the interior part of the computational domain, the square $S_L:=[0,L]^2$ with $L=100$. We first find a good numerical solution $(\Psi_{in},\Omega_{in})$ of the nonlinear system \eqref{u:nl} in $S_L$, i.e.
\begin{equation}\label{sec3eq2}
\begin{split}
&L_{uv}\Psi_{in}-\Omega_{in}=E_{\Psi,in},\\
&L_{uv}\Omega_{in}+K(u)\partial_u\Omega_{in}+K(v)\partial_v\Omega_{in}+2\Omega_{in}+2\,\mathcal{N}[\Psi_{in},\Omega_{in}]=E_{\Omega,in},
\end{split}
\end{equation}
with small residuals $E_{\Psi,in}$, $E_{\Omega,in}$ in $S_L$. The prescribed boundary conditions are
\begin{equation}\label{sec3eq3}
\Psi_{in}=\Psi_{ext},\qquad \Omega_{in}=\Omega_{ext}\qquad\text{ on }\partial_{ext}S_L:=\{(u,v)\in S_L:\,u=L\text{ or }v=L\}.
\end{equation}
These data are imposed at the boundary collocation nodes up to the numerical residual; no exact matching along the whole edges is required for the gluing below. No boundary conditions are imposed on the interior boundary of the square $\partial_{in}S_L:=\{(u,v)\in S_L:\,u=0\text{ or }v=0\}$, since the symmetries there are already encoded in the change of variables \eqref{u:cov}--\eqref{transported}. Then, keeping $(\Psi_{in},\Omega_{in})$ fixed, we find a numerical eigenpair of the linear system \eqref{u:eig} in $S_L$, i.e.
\begin{equation}\label{sec3eq4}
\begin{split}
&\wt{L}_{uv}\wt\Psi_{in}-\wt\Omega_{in}=E_{\wt\Psi,in},\qquad \wt{L}_{uv}\wt\Omega_{in}+K(u)\partial_u\wt\Omega_{in}+K(v)\partial_v\wt\Omega_{in}+\big(\tfrac32-\wt\lambda_{app}\big)\wt\Omega_{in}\\
&\qquad+2\wt{\mathcal N}_1\big(\Psi_{in},\wt\Omega_{in}\big)+2\wt{\mathcal N}_2\big(\wt\Psi_{in},\Omega_{in}\big)=E_{\wt\Omega,in},
\end{split}
\end{equation}
with small residuals $E_{\wt\Psi,in}$, $E_{\wt\Omega,in}$ in $S_L$. The numerical eigenproblem uses boundary data from the asymptotic expansion and determines the eigenvalue along with the interior eigenfunction. The residuals in \eqref{sec3eq4} are evaluated at the fixed value $\wt\lambda_{app}$ in \eqref{lambda:app}. The approximate eigenfunction is specified by the normalized spline coefficients used in Section~\ref{sec:numerical-residual}; the exterior fields use the exact constants $\tau_2,\tau_3$ fixed in \eqref{constants2}, in the same normalization. Any boundary mismatch, including that due to rounding these constants, is included in the gluing residuals below. 

We then glue the interior solutions to the exterior solutions, which are $\big(\Psi_{ext},\Omega_{ext}\big)$ for the profile and $(\wPsi_{ext},\wOmega_{ext})$ for the eigenpair. For this we need to define good cutoff functions. We define first the function $\vartheta:\R\to[0,1]$ by
\begin{equation}\label{sec3eq5}
\begin{split}
&\vartheta(t):=0\,\,\,\text{ if }t\in(-\infty, 0], \qquad \vartheta(t):=1\,\,\,\text{ if }t\in [1,\infty),\\
&\vartheta(t):=126t^5-420t^6+540t^7-315t^8+70t^9\,\,\,\text{ if }t\in[0,1].
\end{split}
\end{equation}
Since $\vartheta'(t)=630t^4(1-t)^4$ for $t\in[0,1]$, it follows that $\vartheta$ is an increasing $C^4$ cutoff function on $\R$. Then we define 
\begin{equation}\label{res:cutoff}
\chi_0(s):=(1-\vartheta)((s-L+d)/d),\qquad \chi(u,v)=\chi_0(u)\chi_0(v).
\end{equation}
We set $d:=16$ and recall $L=100$, thus $\chi\in C^4([0,\infty)^2)$, $0\leq\chi\leq1$, $\chi=1$ in $S_{84}$, and $\chi=0$ in $[0,\infty)^2\setminus S_{100}$. We define, for $(u,v)\in[0,\infty)^2$,
\begin{equation}\label{sec3eq5.2}
\begin{split}
&\Psi^{(c)}_{app}:=\chi\Psi_{in}+(1-\chi)\Psi_{ext},\qquad
\Omega^{(c)}_{app}:=\chi\Omega_{in}+(1-\chi)\Omega_{ext},\\
&\wt\Psi^{(c)}_{app}:=\chi\,\wt\Psi_{in}+(1-\chi)\wt\Psi_{ext},\qquad
\wt\Omega^{(c)}_{app}:=\chi\wt\Omega_{in}+(1-\chi)\wt\Omega_{ext}.
\end{split}
\end{equation}
Since the interior solutions are $C^{12}$ splines in $S_L$ and the asymptotic solutions are smooth away from the origin, the functions in \eqref{sec3eq5.2} are $C^4$ in $[0,\infty)^2$. We transfer the functions $\Psi^{(c)}_{app}, \Omega^{(c)}_{app}, \wt\Psi^{(c)}_{app}, \wt\Omega^{(c)}_{app}$ to the physical space using the formulas \eqref{u:norm}, with $u=f(r^2)$ and $v=f(z^2)$. Since $f$ is smooth on $[0,\infty)$, this defines $(\Psi_{app},\Omega_{app})\in(C^4_{e,o})^2$ and $(\wt\Psi_{app},\wt\Omega_{app})\in(C^4_{e,e})^2$. The identities \eqref{sf3.3} clearly hold, by definition. 

We estimate the residuals in \eqref{sf3} and \eqref{sf3.2} using \eqref{Emod} and \eqref{Emod2}. Clearly, in $S_{84}$,
\begin{equation}\label{sec3eq7.5}
E_\Psi^{(c)}=E_{\Psi,in},\qquad E_{\Omega}^{(c)}=E_{\Omega,in},\qquad E_{\wt\Psi}^{(c)}=E_{\wt{\Psi},in},\qquad E_{\wt\Omega}^{(c)}=E_{\wt{\Omega}, in}.
\end{equation}
Moreover, in $S_{ext}:=[0,\infty)^2\setminus S_{100}$,
\begin{equation}\label{sec3eq7.6}
\begin{split}
&E_\Psi^{(c)}=0,\qquad E_{\wt\Psi}^{(c)}=0,\qquad E_{\Omega}^{(c)}=E_{\Omega,ext},\qquad E_{\wt\Omega}^{(c)}=E_{\wt\Omega,ext},\\
&E_{\Omega,ext}:=L_{uv}\Omega_{ext}+[K(u)\partial_u+K(v)\partial_v+2]\Omega_{ext}+2\,\mathcal N[\Psi_{ext},\Omega_{ext}],\\
&E_{\wt\Omega,ext}:=\wt L_{uv}\wt\Omega_{ext}+\big[K(u)\partial_u+K(v)\partial_v+\tfrac32-\wt\lambda_{app}\big]\wt\Omega_{ext}\\
&\qquad\qquad+2\wt{\mathcal N}_1\big(\Psi_{ext},\wt\Omega_{ext}\big)+2\wt{\mathcal N}_2\big(\wt\Psi_{ext},\Omega_{ext}\big).
\end{split}
\end{equation}
To find the residuals in the gluing region we define
\begin{equation}\label{sec3eq8}
\Psi_{diff}:=\Psi_{in}-\Psi_{ext},\quad \Omega_{diff}:=\Omega_{in}-\Omega_{ext},\quad \wt\Psi_{diff}:=\wt\Psi_{in}-\wt\Psi_{ext},\quad \wt\Omega_{diff}:=\wt\Omega_{in}-\wt\Omega_{ext}
\end{equation}
in $S_L$, so that, in $S_L$, $\Psi^{(c)}_{app}=\Psi_{ext}+\chi\,\Psi_{diff}$, $\Omega^{(c)}_{app}=\Omega_{ext}+\chi\,\Omega_{diff}$, $\wt\Psi^{(c)}_{app}=\wt\Psi_{ext}+\chi\,\wt\Psi_{diff}$, and $\wt\Omega^{(c)}_{app}=\wt\Omega_{ext}+\chi\,\wt\Omega_{diff}$. The cutoff function enters through the commutators
\begin{equation}\label{sec3eq10}
\begin{split}
[L_{uv},\chi]g&:=\frac{6((u+1)^4-1)}{(u+1)^6}\chi_u\partial_u g+\frac{6((v+1)^4-1)}{(v+1)^6}\chi_v\partial_v g+(L_{uv}\chi)g,\\
[\wt L_{uv},\chi]g&:=\frac{6((u+1)^4-1)}{(u+1)^6}\chi_u\partial_u g+\frac{6((v+1)^4-1)}{(v+1)^6}\chi_v\partial_v g+(\wt L_{uv}\chi)g,
\end{split}
\end{equation}
where $\chi_u:=\partial_u\chi$ and $\chi_v:=\partial_v\chi$, and through the first-order expressions
\begin{equation}\label{sec3eq11}
\begin{split}
\mathcal A_\chi O&:=2K(u)K(v)\big(\chi_u\partial_v O-\chi_v\partial_u O\big)+K(u)\chi_u O,\\
\mathcal B_\chi P&:=2K(u)K(v)\big(\chi_v\partial_u P-\chi_u\partial_v P\big)-K(u)\chi_u P+2K(v)\chi_v P,\\
\wt{\mathcal A}_\chi\wt O&:=2K(u)K(v)\big(\chi_u\partial_v\wt O-\chi_v\partial_u\wt O\big),\\
\wt{\mathcal B}_\chi\wt P&:=2K(u)K(v)\big(\chi_v\partial_u\wt P-\chi_u\partial_v\wt P\big)+2K(v)\chi_v\wt P.
\end{split}
\end{equation}
Using \eqref{sec3eq10}--\eqref{sec3eq11} and the definitions \eqref{sec3eq5.2}, in $S_{band}:=S_{100}\setminus S_{84}$ we have
\begin{equation}\label{sec3eq13}
\begin{split}
E^{(c)}_{\Psi}&=\chi E_{\Psi,in}+[L_{uv},\chi]\Psi_{diff},\\
E^{(c)}_{\Omega}&=\chi E_{\Omega,in}+(1-\chi)E_{\Omega,ext}+[L_{uv},\chi]\Omega_{diff}+\big(K(u)\chi_u+K(v)\chi_v\big)\Omega_{diff}\\
&+2\big(\Psi_{diff}\,\mathcal A_\chi\Omega^{(c)}_{app}+\Omega_{diff}\,\mathcal B_\chi\Psi^{(c)}_{app}\big)
-2\chi(1-\chi)\,\mathcal N[\Psi_{diff},\Omega_{diff}],\\
E^{(c)}_{\wt\Psi}&=\chi E_{\wt\Psi,in}+[\wt L_{uv},\chi]\wt\Psi_{diff},\\
E^{(c)}_{\wt\Omega}&=\chi E_{\wt\Omega,in}+(1-\chi)E_{\wt\Omega,ext}+[\wt L_{uv},\chi]\wt\Omega_{diff}+\big(K(u)\chi_u+K(v)\chi_v\big)\wt\Omega_{diff}\\
&+2\big(\wt\Psi_{diff}\,\mathcal A_\chi\Omega^{(c)}_{app}+\wt\Omega_{diff}\,\mathcal B_\chi\Psi^{(c)}_{app}\big)+2\big(\Psi_{diff}\,\wt{\mathcal A}_\chi\wt\Omega^{(c)}_{app}+\Omega_{diff}\,\wt{\mathcal B}_\chi\wt\Psi^{(c)}_{app}\big)\\
&-2\chi(1-\chi)\big(\wt{\mathcal N}_1(\Psi_{diff},\wt\Omega_{diff})+\wt{\mathcal N}_2(\wt\Psi_{diff},\Omega_{diff})\big).
\end{split}
\end{equation}

The four interior fields are degree-$13$ tensor-product splines computed by collocation in double-double arithmetic. The nonlinear predictor comes from \cite{IoJiPa}. The internal knots are simple, so the fields are $C^{12}$ on $S_{100}$. Applying \eqref{sec3eq2} and \eqref{sec3eq4} to the spline pieces gives the interior residuals. In the transition region, \eqref{sec3eq8} gives the differences between the interior and exterior fields. The exterior residuals in \eqref{sec3eq7.6} enter with the factor $1-\chi$ and hence only where $\max(r,|z|)>2085$. Section~\ref{sec:numerical-residual} bounds these terms on the spline cells and the infinite exterior.

\subsection{Estimating the residuals in $Z'$}\label{sec:main-residual}

With the notation in Subsection~\ref{L0L1}, we estimate both sources in \eqref{cl3.3},
\begin{equation}\label{res:starting}
\mathcal R =-\mathcal R_\Omega+ \mathcal L_3(\Delta^{-1}_{\R^5}\mathcal R_\Psi),\qquad \widetilde{\mathcal R} =-\mathcal R_{\widetilde\Omega} +\mathcal L_3(\Delta^{-1}_{\R^5}\mathcal R_{\widetilde\Psi}) +\widetilde{\mathcal L}_3(\Delta^{-1}_{\R^5}\mathcal R_\Psi).
\end{equation}
Unless another domain is indicated, norms and duality pairings are taken on
$\R^5$, under the identification described in Subsection~\ref{L0L1}. We use real-valued fields and real test functions throughout, and the spaces defined in \eqref{robu1}. In particular,
\begin{equation}\label{res:negative-norm}
\|F\|_{Z'}\leq \|F\|_{\dot H^{-1}(\R^5)}:=\||\nabla|^{-1}F\|_{L^2}.
\end{equation}
Our basic source estimates are ordinary $L^2$ estimates for divergence
fields and weighted $L^2$ estimates for scalar terms. For odd scalar
sources we use the weight $w$ below; for even scalar sources we use the
ordinary cylindrical Hardy multiplier $r$. 

\subsubsection{General estimates}
We first record the cylindrical Hardy inequality, which does not
require any symmetry in $z$.

\begin{lemma}[General Hardy inequality]\label{res:radial-hardy}
For any $h\in \dot{H}^1(\R^5)$,
\begin{equation}\label{res:radial-primal}
\|h/r\|_{L^2}\leq\|\nabla_xh\|_{L^2}\leq\|\nabla h\|_{L^2}.
\end{equation}
Consequently, for any scalar source $F$ and vector-field $A$ on $\R^5$ 
\begin{equation}\label{res:radial-source}
\|F\|_{\dot H^{-1}(\R^5)}\leq\|rF\|_{L^2},\qquad
\|\operatorname{div}_{\R^5}A\|_{\dot H^{-1}(\R^5)}
\leq\|A\|_{L^2}.
\end{equation}
\end{lemma}

\begin{proof}
For smooth compactly supported $h$, integration by parts in $\R^4$ gives
\[
0\leq\|\nabla_xh+xh/r^2\|_{L^2}^2
=\|\nabla_xh\|_{L^2}^2-\|h/r\|_{L^2}^2.
\]
This gives  \eqref{res:radial-primal}. Moreover, for any test function $h$,
\[
|\langle F,h\rangle|\leq\|rF\|_{L^2}\|h/r\|_{L^2},\qquad
|\langle\operatorname{div}_{\R^5}A,h\rangle|
\leq\|A\|_{L^2}\|\nabla h\|_{L^2},
\]
and the desired bounds \eqref{res:radial-source} follow by duality.
\end{proof}

For functions that are odd in $z$ we can improve the scalar multiplier near $z=0$.
Let $k$ be even, positive for $z\neq0$, and assume that the half-line inequality
\begin{equation}\label{res:axial-hypothesis}
\int_0^\infty\frac{|a(z)|^2}{k(z)^2}\,dz
\leq\int_0^\infty|a'(z)|^2\,dz
\end{equation}
holds for any function $a$ with $a(0)=0$. Define
\begin{equation}\label{res:w}
w(r,z):=\frac{r^2k(z)^2}{r^2+k(z)^2}\leq r^2,\qquad
\frac1w=\frac1{r^2}+\frac1{k(z)^2}.
\end{equation}

\begin{lemma}[Odd Hardy inequality]\label{res:hardy-lemma}
(i) For every $h\in \dot{H}^1(\R^5)$ that is odd in $z$ we have
\begin{equation}\label{res:hardy-primal}
\|h/\sqrt w\|_{L^2}\leq\|\nabla h\|_{L^2}.
\end{equation}
Therefore, for any odd scalar source $F$,
\begin{equation}\label{res:hardy-dual}
\|F\|_{\dot H^{-1}(\R^5)}\leq\|\sqrt w\,F\|_{L^2}.
\end{equation}

(ii) The bound \eqref{res:axial-hypothesis} holds if $k:\R\to[0,\infty)$ is an even function defined by
\begin{equation}\label{res:k}
k(z)^2=\frac{|z|(1+|z|)^2(2+|z|)}6.
\end{equation}
\end{lemma}

\begin{proof}
(i) We apply \eqref{res:axial-hypothesis} for each $x$, on both half-lines, and
add \eqref{res:radial-primal}, so
\[
\|h/\sqrt w\|_{L^2}^2
=\|h/r\|_{L^2}^2+\|h/k\|_{L^2}^2
\leq\|\nabla_xh\|_{L^2}^2+\|\partial_zh\|_{L^2}^2.
\]
This proves \eqref{res:hardy-primal}, and \eqref{res:hardy-dual} then follows by duality.

(ii) To prove \eqref{res:axial-hypothesis}, define $h_0:[0,\infty)\to[0,\infty)$ by $h_0(z)=1-(1+z)^{-2}$. Then
\[
h_0(z)\geq 0,\qquad
h_0'(z)=2(1+z)^{-3},\qquad h_0''(z)=-6(1+z)^{-4},\qquad -h_0''/h_0=1/k^2.
\]
For compactly supported smooth $a$ with $a(0)=0$, the substitution $a=h_0b$ and integration by parts yields the identity
\begin{equation}\label{res:axial-square}
\int_0^\infty |a'|^2\,dz-\int_0^\infty |a|^2/k^2\,dz
=\int_0^\infty h_0^2|(a/h_0)'|^2\,dz\geq0.
\end{equation}
This proves the required half-line inequality. 
\end{proof}

We examine now the terms in \eqref{res:starting} that contain the operators $\mathcal L_3$ and
$\widetilde{\mathcal L}_3$. Notice that $\mathcal L_3=\mathcal L_{\Omega^\ast}$ and
$\widetilde{\mathcal L}_3=\mathcal L_{\wt{\Omega}^{\ast}}$, where for any real $C^1$ axi-symmetric coefficient $H$ we set
\[
\mathcal L_Hg:=H_z(rg_r+2g)-rH_rg_z=\mathcal N_5(g,H),
\qquad
\mathcal Tg:=(g_r+2g/r,g_z).
\]

\begin{lemma}[Bounds on $\|\mathcal{L}_Hg\|_{\dot H^{-1}(\R^5)}$]\label{res:coupling-lemma}
 For axi-symmetric $g\in\dot H^1(\R^5)$,
\begin{equation}\label{res:T-isometry}
\|\mathcal Tg\|_{L^2}=\|\nabla g\|_{L^2}.
\end{equation}
If $rH\,\mathcal Tg\in L^2$, then
\begin{equation}\label{res:L3local}
\|\mathcal L_Hg\|_{\dot H^{-1}(\R^5)}\leq\|rH\,\mathcal Tg\|_{L^2}\leq\|rH\|_{L^\infty}\|\nabla g\|_{L^2}.
\end{equation}
\end{lemma}

\begin{proof}
For compactly supported smooth $g$, use the radial measure $r^3\,dr$:
\[
\int_0^\infty r^3|g_r+2g/r|^2\,dr
=\int_0^\infty r^3|g_r|^2\,dr+
4\int_0^\infty r^2gg_r\,dr+4\int_0^\infty rg^2\,dr
=\int_0^\infty r^3|g_r|^2\,dr.
\]
After integration in $z$ and in the angular variables this proves
\eqref{res:T-isometry}. Moreover, integration by parts gives, for any axi-symmetric test function $h$,
\begin{equation}\label{res:coupling-pairing}
\langle\mathcal L_Hg,h\rangle
=\int_{\R^5}rH\left[
g_z(h_r+2h/r)-h_z(g_r+2g/r)\right]\,dx\,dz.
\end{equation}
Therefore
\[
|\langle\mathcal L_Hg,h\rangle|
\leq\|rH\,\mathcal Tg\|_{L^2}\|\mathcal Th\|_{L^2}
=\|rH\,\mathcal Tg\|_{L^2}\|\nabla h\|_{L^2}.
\]
Duality proves the first inequality in \eqref{res:L3local}, and the second inequality follows from \eqref{res:T-isometry}.
\end{proof}

For $g=\Delta^{-1}_{\R^5}F$ we have $\|\nabla g\|_{L^2}=\|F\|_{\dot H^{-1}(\R^5)}$. Therefore \eqref{res:L3local}, \eqref{res:radial-source}, and \eqref{res:hardy-dual} give, for any axi-symmetric coefficient $H$,
\begin{equation}\label{res:near}
\|\mathcal L_H(\Delta^{-1}_{\R^5}F)\|_{\dot H^{-1}(\R^5)}\leq\|rH\|_{L^\infty}\|F\|_{\dot H^{-1}(\R^5)}\leq
\begin{cases}
\|rH\|_{L^\infty}\|\sqrt w\,F\|_{L^2}&\text{ if $F$ is odd in $z$},\\
\|rH\|_{L^\infty}\|rF\|_{L^2}&\text{ if $F$ is even in $z$}.
\end{cases}
\end{equation}

These bounds are efficient for sources concentrated near the origin. For sources supported far from the origin the weights $w$ and $r^2$ are large, and we use instead the Newton kernel \eqref{sf10} in a large ball. In the rest of this subsection $R:=|y|$ for $y=(x,z)\in\R^5$, and $|\mathbb S^4|=8\pi^2/3$.

\begin{lemma}[Distant sources]\label{res:far-lemma}
Assume that $0<R_0<R_1$ and that $F\in L^2(\R^5)$ is compactly supported, axi-symmetric, odd in $z$, and supported in $\{R\geq R_1\}$. Let $B_{R_0}:=\{y\in\R^5:\,|y|\leq R_0\}$,
\begin{equation}\label{res:Phi}
\Phi(R)^2:=\frac{3}{8\pi^2}\sum_{j=0}^4\binom4j\frac{R^{4-j}}{(7-j)(R_1-R)^{7-j}},\quad I_o(H,R_0):=\|(r+2|z|)\,\Phi(|y|)H\|_{L^2(B_{R_0})}.
\end{equation}
Then, for any axi-symmetric coefficient $H$,
\begin{equation}\label{res:far-bound}
\|\mathcal L_H(\Delta^{-1}_{\R^5}F)\|_{\dot H^{-1}(\R^5)}^2\leq I_o(H,R_0)^2\,\|F\|_{L^2}^2+\|rH\|_{L^\infty(\R^5\setminus B_{R_0})}^2\,\|\sqrt w\,F\|_{L^2}^2.
\end{equation}
\end{lemma}

\begin{proof}
Let $g=\Delta^{-1}_{\R^5}F$. Notice that $\Phi$ is increasing on $[0,R_1)$. If $|y|<R_1$ then we use \eqref{sf10}, the identity $|\nabla_y|y-Y|^{-3}|=3|y-Y|^{-4}$, and the Cauchy inequality to estimate
\begin{equation}\label{sec3eq40}
|\nabla g(y)|^2\leq\frac{9\|F\|_{L^2}^2}{64\pi^4}\int_{\{Y\in\R^5:\,|Y|\geq R_1\}}\frac{dY}{(|Y|-|y|)^8}
=\frac{3\|F\|_{L^2}^2}{8\pi^2}\int_{R_1}^\infty\frac{s^4\,ds}{(s-|y|)^8}=\Phi(|y|)^2\|F\|_{L^2}^2.
\end{equation}
Since $g$ is odd in $z$, integration along vertical segments gives $|g(y)|\leq|z|\Phi(|y|)\|F\|_{L^2}$ for $|y|<R_1$. As $|r\mathcal Tg|\leq r|\nabla g|+2|g|$, it follows that
\[
\|rH\mathcal Tg\|_{L^2(B_{R_0})}\leq I_o(H,R_0)\,\|F\|_{L^2}.
\]
In $\R^5\setminus B_{R_0}$, \eqref{res:T-isometry} and \eqref{res:hardy-dual} give
\[
\begin{split}
\|rH\mathcal Tg\|_{L^2(\R^5\setminus B_{R_0})}\leq\|rH\|_{L^\infty(\R^5\setminus B_{R_0})}\|\nabla g\|_{L^2}\leq\|rH\|_{L^\infty(\R^5\setminus B_{R_0})}\|\sqrt w\,F\|_{L^2}.
\end{split}
\]
Adding the squares and using \eqref{res:L3local} gives \eqref{res:far-bound}.
\end{proof}

\subsubsection{Bounds on the residuals}\label{sec:residual-bounds}
We identify functions of $(u,v)$ with functions of $(r,z)$ via \eqref{u:cov}, and with their lifts to $\R^5$. Let $\mathbf 1_{\rm in}$, $\mathbf 1_{\rm band}$, $\mathbf 1_{\rm ext}$ denote the characteristic functions of $S_{84}$, $S_{band}$, $S_{ext}$, and for $\alpha\in\{\Psi,\Omega,\wPsi,\wOmega\}$ let
\[
\mathcal R_\alpha^{\rm in}:=\mathcal R_\alpha\mathbf 1_{\rm in},\qquad
\mathcal R_\alpha^{\rm band}:=\mathcal R_\alpha\mathbf 1_{\rm band},\qquad
\mathcal R_\alpha^{\rm ext}:=\mathcal R_\alpha\mathbf 1_{\rm ext}.
\]
These pieces have the same parity in $z$ as $\mathcal R_\alpha$. By \eqref{sec3eq7.6},
\begin{equation}\label{res:exterior-identities}
\mathcal R_\Psi^{\rm ext}=\mathcal R_{\wPsi}^{\rm ext}=0.
\end{equation}
Since $\max(r,|z|)>U_0(84)^{1/2}$ in $S_{band}$, the band pieces are supported in $\{R\geq r_\ast\}$, where
\begin{equation}\label{res:source-radius}
r_\ast:=U_0(84)^{1/2}=2085.677\ldots.
\end{equation}
All the norms below are integrals in the variables $(u,v)$: if $h_o$ is odd in $z$,  $h_e$ is even in $z$, and we define $h_o(r,z)=zh_o^{(c)}(u,v)$ and $h_e(r,z)=h_e^{(c)}(u,v)$, then
\begin{equation}\label{res:measure}
\begin{split}
&\|h_o\|_{L^2(\R^5)}^2=4\pi^2\int_0^\infty\!\!\int_0^\infty|h_o|^2\,r^3\,dr\,dz
=\frac{\pi^2}{9}\int_0^\infty\!\!\int_0^\infty|h_o^{(c)}|^2U_0(u)(u+1)^3(v+1)^3U_0(v)^{1/2}\,du\,dv,\\
&\|h_e\|_{L^2(\R^5)}^2=4\pi^2\int_0^\infty\!\!\int_0^\infty|h_e|^2\,r^3\,dr\,dz
=\frac{\pi^2}{9}\int_0^\infty\!\!\int_0^\infty|h_e^{(c)}|^2\,\frac{U_0(u)(u+1)^3(v+1)^3}{U_0(v)^{1/2}}\,du\,dv.
\end{split}
\end{equation}

Let $I_o$ be defined by \eqref{res:Phi} with $R_1=r_\ast$. Table~\ref{certi1} gives the coefficient and residual bounds proved in Subsection~\ref{cap:profile-verification}.

\begin{table}[!htbp]
\centering
\setlength{\arraycolsep}{12pt}
\setlength{\belowcaptionskip}{0pt}
\caption{Bounds on residuals and coefficients.}\label{certi1}\label{certi2}
\[
\begin{array}{|c|c|}
\hline
\rule{0pt}{2.6ex}\text{nonlinear solution}&\text{eigenpair}\\[4pt] \hline
\|r\Omega^\ast\|_{L^\infty}\leq13
&\|r\wt\Omega^\ast\|_{L^\infty}\leq0.025\\[8pt]
\|r\Omega^\ast\|_{L^\infty(\R^5\setminus B_{900})}\leq5\times10^{-4}
&\|r\wt\Omega^\ast\|_{L^\infty(\R^5\setminus B_{900})}\leq10^{-18}\\[8pt]
I_o(\Omega^\ast,900)\leq0.38
&I_o(\wt\Omega^\ast,900)\leq3.2\times10^{-5}\\[8pt]
&\|\wt\Omega^\ast\|_{L^2}\in[1,1.01]\\[8pt]
\hline
\rule{0pt}{3.2ex}\|\sqrt w\,\mathcal R_\Omega\|_{L^2}\leq4.3\times10^{-12}
&\|r\mathcal R_{\wOmega}\|_{L^2}\leq1.25\times10^{-12}\\[8pt]
\|\sqrt w\,\mathcal R_\Psi^{\rm in}\|_{L^2}\leq6.8\times10^{-14}
&\|r\mathcal R_{\wPsi}\|_{L^2}\leq1.6\times10^{-14}\\[8pt]
\|\mathcal R_\Psi^{\rm band}\|_{L^2}\leq5.6\times10^{-13}
&\\[8pt]
\|\sqrt w\,\mathcal R_\Psi^{\rm band}\|_{L^2}\leq1.25\times10^{-9}
&\\[8pt]
\hline
\end{array}
\]
\end{table}

\begin{lemma}\label{res:bounds-lemma}
We have
\begin{equation}\label{res:main-bound}
\begin{split}
\|\mathcal R\|_{\dot H^{-1}(\R^5)}\leq{}&\|\sqrt w\,\mathcal R_\Omega\|_{L^2}+\|r\Omega^\ast\|_{L^\infty}\|\sqrt w\,\mathcal R_\Psi^{\rm in}\|_{L^2}\\
&+\Big(I_o(\Omega^\ast,900)^2\,\|\mathcal R_\Psi^{\rm band}\|_{L^2}^2+\|r\Omega^\ast\|_{L^\infty(\R^5\setminus B_{900})}^2\,\|\sqrt w\,\mathcal R_\Psi^{\rm band}\|_{L^2}^2\Big)^{1/2}
\end{split}
\end{equation}
and
\begin{equation}\label{res:even-total}
\begin{split}
\|\wt{\mathcal R}\|_{\dot H^{-1}(\R^5)}\leq{}&\|r\mathcal R_{\wOmega}\|_{L^2}+\|r\Omega^\ast\|_{L^\infty}\|r\mathcal R_{\wPsi}\|_{L^2}+\|r\wt\Omega^\ast\|_{L^\infty}\|\sqrt w\,\mathcal R_\Psi^{\rm in}\|_{L^2}\\
&+\Big(I_o(\wt\Omega^\ast,900)^2\,\|\mathcal R_\Psi^{\rm band}\|_{L^2}^2+\|r\wt\Omega^\ast\|_{L^\infty(\R^5\setminus B_{900})}^2\,\|\sqrt w\,\mathcal R_\Psi^{\rm band}\|_{L^2}^2\Big)^{1/2}.
\end{split}
\end{equation}
\end{lemma}

\begin{proof}
By \eqref{res:starting} and \eqref{res:exterior-identities},
\[
\begin{split}
\mathcal R&=-\mathcal R_\Omega+\mathcal L_3(\Delta^{-1}_{\R^5}\mathcal R_\Psi^{\rm in})+\mathcal L_3(\Delta^{-1}_{\R^5}\mathcal R_\Psi^{\rm band}),\\
\wt{\mathcal R}&=-\mathcal R_{\wOmega}+\mathcal L_3(\Delta^{-1}_{\R^5}\mathcal R_{\wPsi})+\wt{\mathcal L}_3(\Delta^{-1}_{\R^5}\mathcal R_\Psi^{\rm in})+\wt{\mathcal L}_3(\Delta^{-1}_{\R^5}\mathcal R_\Psi^{\rm band}).
\end{split}
\]
We use the triangle inequality, \eqref{res:negative-norm}, and the identities $\mathcal L_3=\mathcal L_{\Omega^\ast}$, $\wt{\mathcal L}_3=\mathcal L_{\wt\Omega^\ast}$. The odd source $\mathcal R_\Omega$ is estimated by \eqref{res:hardy-dual}, the even source $\mathcal R_{\wOmega}$ by \eqref{res:radial-source}, and the three terms with the sources $\mathcal R_\Psi^{\rm in}$ and $\mathcal R_{\wPsi}$ by \eqref{res:near}. The two terms with $\mathcal R_\Psi^{\rm band}$ are estimated by Lemma~\ref{res:far-lemma} with $R_1=r_\ast$, $R_0=900$, and $H=\Omega^\ast$, respectively $H=\wt\Omega^\ast$; the support condition holds by \eqref{res:source-radius}.
\end{proof}

Applying Lemma~\ref{res:bounds-lemma} with Table~\ref{certi1} gives
\begin{equation}\label{res:final}
\begin{split}
\|\mathcal R\|_{Z'}&\leq\Big(4.3+13\cdot0.068+\big((0.38\cdot0.56)^2+(5\cdot10^{-4}\cdot1250)^2\big)^{1/2}\Big)\cdot10^{-12}<5.85\times10^{-12},\\
\|\wt{\mathcal R}\|_{Z'}&\leq\Big(1.25+13\cdot0.016+0.025\cdot0.068\\
&\qquad+\big((3.2\cdot10^{-5}\cdot0.56)^2+(10^{-18}\cdot1250)^2\big)^{1/2}\Big)\cdot10^{-12}<1.5\times10^{-12}.
\end{split}
\end{equation}

\section{Stabilized coercivity}\label{sec:coercivity}

In this section we prove Theorem~\ref{c:thm:main}, which gives \eqref{robu7} and \eqref{robu31} with $\gamma=1/10$, $n=26$ and $\widetilde n=32$. Both parity classes are treated by the same scheme: the negative part of the quadratic form is confined to finitely many cells, and on each cell it is controlled by finitely many modes of a weighted Neumann Laplacian. The computer-assisted inputs are collected in Proposition~\ref{c:prop:CAP} at the end of the section.

Throughout this section, $\Delta=\Delta_{\R^5}$ and $\R^2_+=(0,\infty)\times\R$. We prove the coercive inequalities with the constant
\begin{equation}\label{c:parameters}
\gamma_0=0.100001
\end{equation}
in place of $\gamma=1/10$ to absorb the error \eqref{c:profile-loss} from the difference between $\Psi_5$ and $\Psi^\ast$.

\begin{theorem}\label{c:thm:main}
The functions defined in \eqref{c:qfinal} belong to the respective parity
classes, $q_j\in(Z')^{odd}$ for $1\leq j\leq26$ and
$\widetilde q_j\in(Z')^{even}$ for $1\leq j\leq32$, and satisfy
\begin{equation}\label{c:qbounds}
\sum_{j=1}^{26}\|q_j\|_{Z'} \leq19.2,\qquad \sum_{j=1}^{26}\|q_j\|_{Z'}^2\leq19.86,\qquad \sum_{j=1}^{32}\|\widetilde q_j\|_{Z'} \leq34,\qquad \sum_{j=1}^{32}\|\widetilde q_j\|_{Z'}^2\leq61.
\end{equation}
For every $f\in Z_0^{odd}$ and $h\in Z_0^{even}$,
\begin{equation}\label{c:main1}
\langle\mathcal S_0f,\Delta^{-1}f\rangle +\sum_{j=1}^{26}\langle q_j,\Delta^{-1}f\rangle^2 \geq\tfrac1{10}\|f\|_Z^2,\qquad \langle\mathcal S_{\widetilde\lambda_{app}}h,\Delta^{-1}h\rangle +\sum_{j=1}^{32}\langle\widetilde q_j,\Delta^{-1}h\rangle^2 \geq\tfrac1{10}\|h\|_Z^2.
\end{equation}
\end{theorem}

The second inequality holds on the full even space, hence on the normalized
subspace used in Lemma~\ref{weakCoercive2}.

\subsection{The quadratic form in terms of $g$}\label{c:ss:start}

We compute the coefficients of the quadratic form with the glued spline potential $\Psi_5$ in place of $\Psi^\ast$. Define
\begin{equation}\label{c:profile-operator}
\mathcal S^{(5)}_\lambda f:=\mathcal O_{\R^5}f-\lambda f
 +\mathcal N_5(\Psi_5,f)+\mathcal N_5(\Delta^{-1}f,\Omega^\ast).
\end{equation}
Only the potential in the first bilinear term differs from
$\mathcal S_\lambda$. By \eqref{cl3.3} and Lemma~\ref{nl:bilinear},
\begin{equation}\label{c:profile-loss}
\begin{split}
(\mathcal S_\lambda-\mathcal S^{(5)}_\lambda)f
 &=-\mathcal N_5(\Delta^{-1}\mathcal R_\Psi,f),\\
\big|\langle(\mathcal S_\lambda-\mathcal S^{(5)}_\lambda)f,
                         \Delta^{-1}f\rangle\big|
 &\leq C_b\|\mathcal R_\Psi\|_Z\|f\|_Z^2
 <1.001\times10^{-10}\|f\|_Z^2.
\end{split}
\end{equation}

In what follows, we write $g_r=\partial_rg$, $g_{rz}=\partial_r\partial_zg$, and
$\nabla g=(g_r,g_z)$ when integrating on $\R^2_+$.

\begin{lemma}\label{c:lem:Z}
Let $f\in Z$ be axially symmetric and let $g:=\Delta^{-1}f$. Then
\begin{equation}\label{c:Z1}
\frac{1}{2\pi^2}\|f\|_Z^2=\int_{\R^2_+}\big(g_r^2+g_z^2\big)\,r^3drdz
+\int_{\R^2_+}\Big(g_{rr}^2+g_{zz}^2+2g_{rz}^2+\frac{3}{r^2}g_r^2\Big)\,r^3drdz .
\end{equation}
\end{lemma}

\begin{proof}
From the definition \eqref{robu1} and Plancherel's theorem, $\|f\|_Z^2=\|\nabla_{\R^5}g\|^2_{L^2}+\|D^2g\|^2_{L^2}$, where $D^2g$ is the Hessian of $g$ on $\R^5$. Expanding in coordinates and integrating over the angular variables of $\R^4$, which gives the factor $|\mathbb S^3|=2\pi^2$, we obtain \eqref{c:Z1}.
\end{proof}

\begin{proposition}\label{c:prop:ident}
Let $\lambda\geq0$, let $f\in C^\infty_c(\R^5)$ be axially symmetric, and let
$g:=\Delta^{-1}f$. Then
\begin{equation}\label{mos1}
\begin{split}
\frac{1}{2\pi^2}\langle \mathcal S^{(5)}_\lambda f,\Delta^{-1}f\rangle
&=\int_{\R^2_+}\big(g_{rr}^2+g_{zz}^2+2g_{rz}^2\big)\,r^3drdz
+\int_{\R^2_+}\Big(\frac14+\lambda+\frac{3}{r^2}+e_{rr}\Big)g_r^2\,r^3drdz\\
&\quad+\int_{\R^2_+}\Big(\frac14+\lambda+e_{zz}\Big)g_z^2\,r^3drdz
+\int_{\R^2_+}2a_{12}\,g_rg_z\,r^3drdz+\int_{\R^2_+} e_0\, g^2\,r^3drdz,
\end{split}
\end{equation}
where
\begin{equation}\label{mos2}
\begin{split}
&e_{rr} :=-r\partial_r\partial_z\Psi_5-2\partial_z\Psi_5,\qquad e_{zz} :=r\partial_r\partial_z\Psi_5+\partial_z\Psi_5,\\
&a_{12} :=\tfrac12\big[r(\partial_r^2\Psi_5-\partial_z^2\Psi_5)+3\partial_r\Psi_5\big],\qquad e_0 :=\partial_z\Delta\Psi_5.
\end{split}
\end{equation}
\end{proposition}

\begin{proof}
Here we write $\nabla$ for the gradient in $\R^5$, whose indices we denote by $j,k\in\{1,\ldots,4\}$, $\al,\be\in\{1,\ldots,5\}$. Using \eqref{cl3} and $f=\Delta g$ we write
$\langle\mathcal S^{(5)}_\lambda f,\Delta^{-1}f\rangle=J_1+J_2+J_3$, where
\begin{equation}\label{pin51}
\begin{split}
&J_1:=\int_{\R^5}\Big(\Delta f+\tfrac12 x_\al\partial_\al f+(\tfrac32-\lambda)f\Big)g\,dxdz,\\
&J_2:=\int_{\R^5}\Big(-\partial_z\Psi_5\cdot x_j\partial_jf+(x_j\partial_j\Psi_5+2\Psi_5)\partial_zf\Big)g\,dxdz,\\
&J_3:=\int_{\R^5}\Big(-(x_j\partial_j\Omega^\ast)\partial_zg+\partial_z\Omega^\ast\cdot x_j\partial_jg+2\partial_z\Omega^\ast\,g\Big)g\,dxdz .
\end{split}
\end{equation}
Since $f\in C^\infty_c(\R^5)$, we have $g=O(R^{-3})$ and $\nabla g=O(R^{-4})$ as $R\to\infty$, so all the integrals below converge absolutely and all the integrations by parts are justified.

First consider the term $J_1$. Since $[\Delta,x_\al\partial_\al]=2\Delta$, we have $x_\al\partial_\al\Delta g=\Delta(x_\al\partial_\al g)-2\Delta g$. Therefore,
\[
\tfrac12\int(x_\al\partial_\al\Delta g)\,g=\tfrac12\int(\Delta (x_\al\partial_\al g))g-\int(\Delta g)g
=\tfrac12\int(x_\al\partial_\al g)\Delta g+\int|\nabla g|^2 .
\]
Integrating by parts again and using $\partial_\al x_\al=5$ we have
\[
\int(x_\al\partial_\al g)\Delta g=-\int|\nabla g|^2-\tfrac12\int x_\al\partial_\al\big(|\nabla g|^2\big)
=-\int|\nabla g|^2+\tfrac{5}{2}\int|\nabla g|^2=\tfrac32\int|\nabla g|^2 .
\]
Summing these identities, we have
\begin{equation}\label{pin52}
J_1=\int_{\R^5}(\Delta g)^2\,dxdz+\Big(\frac14+\lambda\Big)\int_{\R^5}|\nabla g|^2\,dxdz .
\end{equation}

Next, $J_3$ vanishes identically:
\[
J_3=\frac12\int x_j\partial_j\partial_z\Omega^\ast\,g^2
 -\frac12\int(4\partial_z\Omega^\ast+x_j\partial_j\partial_z\Omega^\ast)g^2
 +2\int\partial_z\Omega^\ast g^2=0.
\]

Finally, consider the term $J_2$. Let $V_j:=-x_j\partial_z\Psi_5$ for $j\in\{1,\ldots,4\}$ and
$V_5:=x_j\partial_j\Psi_5+2\Psi_5$, so $J_2=\int V_\al\partial_\al\Delta g\cdot g$.
We integrate by parts in $\partial_\be$ to find
\begin{equation}\label{pin54}
J_2=\int_{\R^5}\Big\{\frac{\partial_\al V_\al}{2}|\nabla g|^2
+\partial_\be V_\al\,\partial_\al g\,\partial_\be g
-\frac{\Delta(\partial_\al V_\al)}{2}\,g^2\Big\}\,dxdz .
\end{equation}
Since $\partial_j(-x_j\partial_z\Psi_5)=-4\partial_z\Psi_5-x_j\partial_j\partial_z\Psi_5$ and
$\partial_zV_5=x_j\partial_j\partial_z\Psi_5+2\partial_z\Psi_5$, we have
\begin{equation}\label{pin55}
\partial_\al V_\al=-2\partial_z\Psi_5,\qquad
-\tfrac12\Delta(\partial_\al V_\al)=\Delta\partial_z\Psi_5=e_0 .
\end{equation}
To calculate the middle term in \eqref{pin54} we use the axial symmetry. Since $\partial_jg=(x_j/r)g_r$ and $\partial_j\Psi_5=(x_j/r)\partial_r\Psi_5$, we have
\[
\partial_kV_j=-\delta_{jk}\partial_z\Psi_5-\frac{x_jx_k}{r}\partial_r\partial_z\Psi_5,
\qquad
\partial_zV_j+\partial_jV_5=x_j\Big(\partial_r^2\Psi_5-\partial_z^2\Psi_5+\frac3r\partial_r\Psi_5\Big).
\]
Therefore,
\[
\partial_\be V_\al\,\partial_\al g\,\partial_\be g
=-\big(\partial_z\Psi_5+r\partial_r\partial_z\Psi_5\big)g_r^2
+\big(r\partial_r\partial_z\Psi_5+2\partial_z\Psi_5\big)g_z^2
+2a_{12}\,g_rg_z .
\]
Since $\frac12(\partial_\al V_\al)|\nabla g|^2=-\partial_z\Psi_5(g_r^2+g_z^2)$, it follows that
\begin{equation}\label{pin56}
\begin{split}
\frac{1}{2\pi^2}J_2=&\int_{\R^2_+}\big(-2\partial_z\Psi_5-r\partial_r\partial_z\Psi_5\big)g_r^2\,r^3drdz
+\int_{\R^2_+}\big(\partial_z\Psi_5+r\partial_r\partial_z\Psi_5\big)g_z^2\,r^3drdz\\
&+\int_{\R^2_+} 2a_{12}\,g_rg_z\,r^3drdz+\int_{\R^2_+} e_0\, g^2\,r^3drdz .
\end{split}
\end{equation}
Collecting the terms yields \eqref{mos1}.
\end{proof}

By Proposition~\ref{c:prop:ident} and Lemma~\ref{c:lem:Z},
\begin{equation}\label{c:Q}
\frac{1}{2\pi^2}\big(\langle\mathcal S^{(5)}_\lambda f,\Delta^{-1}f\rangle
      -\gamma_0\|f\|_Z^2\big)=\mathcal Q_{\gamma_0,\lambda}(g),
\end{equation}
where
\[
\begin{split}
\mathcal Q_{\gamma_0,\lambda}(g)
={}&(1-\gamma_0)\int_{\R^2_+}(|\nabla g_r|^2+|\nabla g_z|^2)r^3\,drdz\\
&+\int_{\R^2_+}\big((a_{11}+\lambda)g_r^2+(a_{22}+\lambda)g_z^2+2a_{12}g_rg_z+e_0g^2\big)r^3\,drdz,
\end{split}
\]
and
\begin{equation}\label{c:coef}
a_{11}:=\frac14-\gamma_0+\frac{3(1-\gamma_0)}{r^2}+e_{rr},
\qquad a_{22}:=\frac14-\gamma_0+e_{zz}.
\end{equation}

\subsection{The odd quadratic form}\label{c:ss:Hardy}

The coefficient $e_0=\partial_z\Delta\Psi_5$ of the zeroth order term in \eqref{c:Q} has a nontrivial negative part $(e_0)_-:=\max(-e_0,0)$, which we absorb into the first order term using the following Poincar\'e inequality in $z$. Note that the oddness in $z$ is being used in a crucial way.

\begin{lemma}\label{c:lem:odd}
Let $L>0$ and let $\varphi\in C^1([-L,L])$ be an odd function. Then
\begin{equation}\label{c:odd1}
\int_{-L}^L\varphi^2\,dz\leq\frac{4L^2}{\pi^2}\int_{-L}^L(\varphi')^2\,dz .
\end{equation}
\end{lemma}

\begin{proof}
Since $\varphi$ is odd, it suffices to prove the inequality on $[0,L]$, for functions with $\varphi(0)=0$. Then \eqref{c:odd1} is the Poincar\'e inequality on $[0,L]$ with a Dirichlet condition at $0$ and a Neumann condition at $L$, whose first eigenvalue is $\pi^2/(4L^2)$.
\end{proof}

To apply Lemma~\ref{c:lem:odd}, we bound $(e_0)_-(r,\cdot)$ at every fixed $r$ by a finite sum of step functions of $z$ of the form $h_p(r)\mathbf 1_{\{|z|<L_p(r)\}}$, which we call boxes.

\begin{corollary}\label{c:cor:zero}
Assume that for every $r>0$ we are given finitely many pairs of numbers $h_p(r)\geq0$, $L_p(r)\geq0$ such that
\begin{equation}\label{c:box}
(e_0)_-(r,z)\leq\sum_ph_p(r)\,\mathbf 1_{\{|z|<L_p(r)\}}
\qquad\text{for all }(r,z)\in\R^2_+ ,
\end{equation}
and let
\begin{equation}\label{c:zeta}
\zeta(r,z):=\frac{4}{\pi^2}\sum_ph_p(r)L_p(r)^2\,\mathbf 1_{\{|z|<L_p(r)\}} .
\end{equation}
Then, for any smooth function $g$ that is odd in $z$,
\begin{equation}\label{c:zero1}
\int_{\R^2_+}(e_0)_-\,g^2\,r^3drdz\leq\int_{\R^2_+}\zeta\,g_z^2\,r^3drdz .
\end{equation}
\end{corollary}

\begin{proof}
Using \eqref{c:box} and Fubini's theorem,
\[\int_{\R^2_+}(e_0)_-g^2\,r^3drdz\leq\int_0^\infty r^3\sum_ph_p(r)\big(\int_{-L_p}^{L_p}g(r,z)^2dz\big)dr.\]
We apply Lemma~\ref{c:lem:odd} to $\varphi=g(r,\cdot)$ on $[-L_p,L_p]$, so
\begin{align*}
    \int_{-L_p}^{L_p}g(r,z)^2dz\leq \frac{4L_p^2}{\pi^2}\int_{-L_p}^{L_p}g_z(r,z)^2dz.
\end{align*}
The claim then follows after multiplying by $h_p(r)$, summing over $p$, and integrating in $r$.
\end{proof}

We defer the choice of the boxes to Proposition~\ref{c:prop:boxes} below. Having eliminated the zeroth order contribution, we can now bound $\mathcal Q_{\gamma_0,0}(g)$ from below by a quadratic form in $\nabla g$.

\begin{proposition}\label{c:prop:vec}
Let $\zeta$ be defined by \eqref{c:zeta} for boxes satisfying \eqref{c:box}. For a vector field $u=(u_1,u_2)$ on $\R^2_+$ define
\begin{equation}\label{c:Mdef}
\mathcal E(u):=(1-\gamma_0)\int_{\R^2_+}|\nabla u|^2\,r^3drdz+\int_{\R^2_+} u^{T}M^ou\,r^3drdz,
\qquad
M^o:=\begin{pmatrix}a_{11}&a_{12}\\ a_{12}&a_{22}-\zeta\end{pmatrix}.
\end{equation}
Let $f\in C^\infty_c(\R^5)$ be axially symmetric and odd in $z$, and let $g=\Delta^{-1}f$. Then $\mathcal Q_{\gamma_0,0}(g)\geq\mathcal E(\nabla g)$. Consequently, if $q_1,\ldots,q_n\in Z'$ satisfy
\begin{equation}\label{c:vec2}
\mathcal E(\nabla g)+\frac{1}{2\pi^2}\sum_{j=1}^n\langle q_j,g\rangle^2\geq0,
\end{equation}
then
\begin{equation}\label{c:vec3}
\langle\mathcal S^{(5)}_0f,\Delta^{-1}f\rangle+\sum_{j=1}^n\langle q_j,\Delta^{-1}f\rangle^2\geq\gamma_0\|f\|_Z^2.
\end{equation}
\end{proposition}

\begin{proof}
By Corollary~\ref{c:cor:zero},
\begin{align*}
    \int_{\R^2_+} e_0g^2\,r^3drdz\geq-\int_{\R^2_+}(e_0)_-g^2\,r^3drdz\geq-\int_{\R^2_+}\zeta g_z^2\,r^3drdz.
\end{align*}
We insert this in the definition of $\mathcal Q_{\gamma_0,0}$, thus $\mathcal Q_{\gamma_0,0}(g)\geq\mathcal E(\nabla g)$. The second claim then follows from \eqref{c:Q}.
\end{proof}

\subsection{Confinement of the odd negative part}\label{c:ss:far}

Next we fix the boxes and show that $M^o\geq0$ outside a bounded rectangle. We start with the far field, recalling that $R=(r^2+z^2)^{1/2}$.

\begin{lemma}\label{c:lem:lead}
Let $e^{(1)}_{rr}$, $e^{(1)}_{zz}$, $a^{(1)}_{12}$, $e^{(1)}_0$ denote functions defined as in
\eqref{mos2}, but with $\Psi_5$ and $\Delta\Psi_5$ replaced by
\[
\Psi_{-1} = \sigma z R^{-2}
\qquad\text{and}\qquad
\Omega_{-3} = L_0 \Psi_{-1} = -6\sigma z R^{-4}.
\]
Then we have the exact formulas
\begin{equation}\label{c:lead1}
\begin{aligned}
e^{(1)}_{rr} &= \frac{2\sigma z^2(z^2-3r^2)}{R^6}, &\qquad
e^{(1)}_{zz} &= \frac{\sigma(-r^4+6r^2z^2-z^4)}{R^6}, \\[1ex]
a^{(1)}_{12} &= -\frac{\sigma rz(5z^2-3r^2)}{R^6}, &\qquad
e^{(1)}_0 &= -\frac{6\sigma(r^2-3z^2)}{R^6},
\end{aligned}
\end{equation}
which admit the estimates
\begin{equation}\label{c:lead2}
|e^{(1)}_{rr}| \le \frac{2\sigma}{R^2}, \qquad |e^{(1)}_{zz}| \le \frac{\sigma}{R^2}, \qquad |a^{(1)}_{12}| \le \frac{1.37\,\sigma}{R^2},\qquad (e^{(1)}_0)_- \le \frac{6\sigma}{r^4}.
\end{equation}
\end{lemma}

\begin{proof}
If $\Psi = \sigma z R^{-2}$, then
\begin{align*}
\partial_z \Psi &= \sigma(r^2-z^2)R^{-4}, &
\partial_r \Psi &= -2\sigma rz R^{-4}, \\
\partial_r\partial_z \Psi &= 2\sigma r(3z^2-r^2)R^{-6}, &
\partial_r^2 \Psi &= -2\sigma z(z^2-3r^2)R^{-6}, \\
\partial_z^2 \Psi &= 2\sigma z(z^2-3r^2)R^{-6}, &
\partial_z \Omega_{-3} &= -6\sigma(R^2-4z^2)R^{-6}.
\end{align*}
Substituting these formulas into \eqref{mos2} gives \eqref{c:lead1}.

To prove \eqref{c:lead2}, let $t := z^2/R^2 \in [0,1]$. Then
\begin{align*}
&\frac{R^2\,|e^{(1)}_{rr}|}{\sigma} = 2t\,|4t-3| \le 2,\qquad \frac{R^2\,|e^{(1)}_{zz}|}{\sigma} = |-8t^2+8t-1| \le 1,\\
&\frac{R^2\,|a^{(1)}_{12}|}{\sigma} = \sqrt{t(1-t)}\,|8t-3| \le 1.368,
\end{align*}
where the last maximum is attained at $t \approx 0.8237$. Finally,
\[
(e^{(1)}_0)_- = 6\sigma(r^2-3z^2)_+ R^{-6} \le 6\sigma r^2 R^{-6} \le 6\sigma r^{-4}.
\qedhere
\]
\end{proof}

\begin{lemma}\label{c:lem:far}
We have the following in the far field:
\begin{enumerate}[(i)]
\item $(e_0)_-\equiv0$ where $R\ge 200$ and $|z|\ge 0.65\,r$;
\item $\sup_{z}\,(e_0)_-(r,z)\le7\sigma r^{-4}$ where $r\ge 200$.
\end{enumerate}
\end{lemma}

\begin{proof}
(i) Assume $R\ge 200$ and $|z|\ge 0.65\,r$, and let $t:=z^2/R^2$. Then
\[
r^2 = (1-t)R^2 \le \frac{z^2}{0.4225} = \frac{tR^2}{0.4225},
\]
so $t\ge 0.4225/1.4225>0.297$, and therefore
\[
3z^2-r^2 = (4t-1)R^2 \ge 0.188\,R^2.
\]
Using \eqref{c:lead1} and the bounds \eqref{c:Hyp}, we obtain
\[
e_0 \ge e^{(1)}_0-\sigma R^{-4}
\ge (6\cdot 0.188-1)\sigma R^{-4} > 0.
\]
If $r\ge 200$, then $R\ge r\ge200$, so by
\eqref{c:Hyp} and \eqref{c:lead2},
\[
(e_0)_- \le (e^{(1)}_0)_- + \sigma R^{-4} \le 7\sigma r^{-4},
\]
which proves (ii).
\end{proof}

\begin{proposition}\label{c:prop:boxes}
Define the boxes as follows:
\begin{enumerate}[(a)]
\item If $r\geq 200$, there is a single box, with $L_1(r):=0.65\,r$ and $h_1(r):=7\sigma r^{-4}$.
\item If $r\in I_k:=[r_{k-1},r_k)$, $1\leq k\leq K_b$, there are two boxes, with $L_p(r):=L_{p,k}$ and $h_p(r):=h_{p,k}$, $p=1,2$, where $r_k$, $L_{p,k}$, $h_{p,k}$ are the data listed before Proposition~\ref{c:prop:CAP}.
\end{enumerate}
Then \eqref{c:box} holds.
\end{proposition}

\begin{proof}
If $r\geq 200$, then $(e_0)_-(r,z)\leq7\sigma r^{-4}$ by Lemma~\ref{c:lem:far}(ii), and $(e_0)_-(r,z)=0$ for $|z|\geq0.65\,r$ by Lemma~\ref{c:lem:far}(i), which gives \eqref{c:box}.

Let now $r\in I_k$. The right-hand side of \eqref{c:box} is $h_{1,k}+h_{2,k}$ if $|z|<L_{1,k}$, $h_{2,k}$ if $L_{1,k}\leq|z|<L_{2,k}$, and $0$ if $|z|\geq L_{2,k}$. If $|z|\leq 200$, \eqref{c:box} therefore follows from the bounds \eqref{c:boxes}. If $|z|\geq 200$, then $|z|\geq L_{2,k}$, and $(e_0)_-(r,z)=0$ by Lemma~\ref{c:lem:far}(i), since $|z|\geq 200>0.65\,r$.
\end{proof}

From now on $\zeta$ is defined with these boxes.

\begin{proposition}\label{c:prop:supp}
The matrix $M^o$ is nonnegative outside the rectangle
\begin{equation}\label{c:R0}
\mathcal{R}_0:=[1,42]\times[-20,20].
\end{equation}
\end{proposition}

\begin{proof}
For $\max(r,|z|)\leq200$, this is Proposition~\ref{c:prop:CAP}(iii). Assume now that $\max(r,|z|)\ge 200$, so that $R\ge 200$, and set
\[
x := \sigma R^{-2} \le \sigma\cdot200^{-2} <0.0033.
\]
If $r\ge 200$, then by Proposition~\ref{c:prop:boxes}(a) and \eqref{c:zeta},
\[
\zeta(r,z) \le \frac{4}{\pi^2}\cdot7\sigma r^{-4}\,(0.65\,r)^2\,
\mathbf 1_{\{|z|<0.65\,r\}}.
\]
Since $R^2<1.4225\,r^2$ when $|z|<0.65\,r$, it follows that
\[
\zeta \le \frac{4\cdot 0.4225\cdot 1.4225}{\pi^2}\cdot7\,\frac{\sigma}{R^2}
\le 0.2436\cdot 7\,x \le 1.71x.
\]
If $r<200\le|z|$, then $\zeta(r,z)=0$, because $|z|\ge 200>L_{2,k}\geq L_{1,k}$ for the two boxes at this $r$, by Proposition~\ref{c:prop:CAP}(ii). Thus $\zeta\le 1.71\,x$ in both cases.

Using \eqref{c:coef}, \eqref{mos2}, \eqref{c:lead2}, and \eqref{c:Hyp}, we have
\[
a_{11} \ge \tfrac14-\gamma_0-2.05\,x>0.14, \qquad
a_{22}-\zeta \ge \tfrac14-\gamma_0-2.76\,x>0.14, \qquad
|a_{12}| \le 1.42\,x<0.005,
\]
so $a_{11}(a_{22}-\zeta)>0.14^2>a_{12}^2$, and $M^o\ge0$.
\end{proof}

\subsection{Cells and the partition of unity}\label{c:ss:cells}

We cover the rectangle $\mathcal{R}_0$ by products of one-dimensional bump functions. These are constructed using the cutoff
\begin{equation}\label{c:S}
\begin{split}
&\chi(t):=126t^5-420t^6+540t^7-315t^8+70t^9\qquad (0\leq t\leq1),\\
&\chi(t):=0\quad (t<0),\qquad \chi(t):=1\quad (t>1).
\end{split}
\end{equation}
For $a<n_1\leq n_2<b$, we define the plateau bump
\begin{equation}\label{c:plateau}
P[a,n_1,n_2,b](x):=
\begin{cases}
\chi\big(\frac{x-a}{n_1-a}\big),&a\leq x\leq n_1,\\[2pt]
1,&n_1\leq x\leq n_2,\\[2pt]
\chi\big(\frac{b-x}{b-n_2}\big),&n_2\leq x\leq b
\end{cases}
\end{equation}
and $P[a,n_1,n_2,b](x)=0$ otherwise. We note that
\begin{equation}\label{c:S1}
\chi'(t)=630\,t^4(1-t)^4\quad (0\leq t\leq1),\qquad
\chi(t)+\chi(1-t)=1\quad (t\in\R).
\end{equation}
A cell is a product $\varrho(r,z)=T(r)U(z)$ of plateau bumps $T=P[r_a,\cdot,\cdot,r_b]$, with $r_a>0$, and $U=P[z_a,\cdot,\cdot,z_b]$. Its support is the rectangle $Q:=[r_a,r_b]\times[z_a,z_b]$, and the cell is symmetric if $U$ is even. We use the radial bumps
\begin{equation}\label{c:cells1}
T_0:=P[0.5,\,1,\,15,\,17],\qquad
T_1:=P[15,\,17,\,29,\,31],\qquad
T_2:=P[29,\,31,\,42,\,45],
\end{equation}
the vertical bumps
\begin{equation}\label{c:cells2}
U_0:=P[-7,\,-5,\,5,\,7],\qquad U_1:=P[5,\,7,\,20,\,23],\qquad U_1^-(z):=U_1(-z)
\end{equation}
to construct the nine cells $\varrho=TU$, $T\in\{T_0,T_1,T_2\}$,
$U\in\{U_0,U_1,U_1^-\}$. The three cells $T_kU_0$ are symmetric, and the cells $T_kU_1$ and $T_kU_1^-$ form a reflected pair. By Lemmas~\ref{c:lem:q}
and~\ref{c:lem:mirror} below, the two cells of a reflected pair lead to the same functionals, so only six cells need to be analyzed.

\begin{lemma}\label{c:lem:pu}
Let $P=P[a,n_1,n_2,b]$ be a plateau bump.
\begin{enumerate}[(i)]
\item $P\in C^4(\R)$, $0\leq P\leq1$, $\{P>0\}=(a,b)$, and $P=1$ on $[n_1,n_2]$.
\item The bumps in \eqref{c:cells1}--\eqref{c:cells2} satisfy $T_0+T_1+T_2=1$ on $[1,42]$, $U_1^-+U_0+U_1=1$ on $[-20,20]$, and both sums are $\leq1$ on $\R$. Therefore the nine cells form a partition of unity on $\mathcal{R}_0$ and satisfy $\sum\varrho\leq1$ on $\R^2_+$.
\end{enumerate}
\end{lemma}

\begin{proof}
(i) follows from \eqref{c:S1}, since $\chi^{(k)}(0)=\chi^{(k)}(1)=0$ for $1\leq k\leq4$, and (ii) from $\chi(t)+\chi(1-t)=1$ on the overlaps of consecutive bumps.
\end{proof}

\subsection{The even quadratic form}\label{c:ss:even}

For the even problem, let $h\in Z_0^{even}$ and $g=\Delta^{-1}h$, and take $\lambda=\widetilde\lambda_{app}$ in \eqref{c:Q}. The function $g$ is now even in $z$, so Corollary~\ref{c:cor:zero} is not available. Instead, the gradient term will absorb the negative part of $e_0$ through a Hardy inequality in $r$. Set
\begin{equation}\label{c:even-split-constants}
\alpha:=0.224,\qquad \beta:=0.189,
\end{equation}
and define
\begin{align}\label{c:even-vector-matrix}
M^v&:=\begin{pmatrix}
\frac14+\widetilde\lambda_{app}-\gamma_0-\alpha
 +\frac{3(1-\gamma_0)}{r^2}+e_{rr}&a_{12}\\
a_{12}&\frac14+\widetilde\lambda_{app}-\gamma_0-\alpha+e_{zz}
\end{pmatrix},\\
V&:=e_0+\frac{\alpha-\beta}{r^2}.\notag
\end{align}
The two forms used in the even construction are
\begin{equation}\label{c:even-forms}
\begin{split}
\widetilde{\mathcal E}_{v}(u)
&:=(1-\gamma_0)\int_{\R^2_+}|\nabla u|^2r^3\,drdz
   +\int_{\R^2_+}u^TM^v u\,r^3\,drdz,\\
\widetilde{\mathcal E}_{s}(g)
&:=\beta\int_{\R^2_+}|\nabla g|^2r^3\,drdz
   +\int_{\R^2_+}Vg^2r^3\,drdz.
\end{split}
\end{equation}
Here $u=\nabla g=(g_r,g_z)$ has respective parities even and odd in $z$, whereas the scalar argument $g$ is even.

The transverse Hardy inequality is
\[
\int_{\R^2_+}\frac{g^2}{r^2}\,r^3\,drdz
\leq\int_{\R^2_+}g_r^2r^3\,drdz.
\]
It follows by integrating $\partial_r(r^2)=2r$ against $g^2$ and applying Cauchy--Schwarz. Consequently
\begin{equation}\label{c:even-Hardy-split}
\alpha\int|\nabla g|^2r^3\,drdz+\int e_0g^2r^3\,drdz
\geq\beta\int|\nabla g|^2r^3\,drdz
 +\int Vg^2r^3\,drdz.
\end{equation}
Thus \eqref{c:Q} gives
\[
\mathcal Q_{\gamma_0,\widetilde\lambda_{app}}(g)
\geq\widetilde{\mathcal E}_v(\nabla g)+\widetilde{\mathcal E}_s(g).
\]

The matrix $M^v$ is related to the odd matrix $M^o$ of \eqref{c:Mdef} by
\begin{equation}\label{c:even-vector-confinement}
M^v=M^o+(\widetilde\lambda_{app}-\alpha)I_2+\begin{pmatrix}0&0\\0&\zeta\end{pmatrix},
\qquad\widetilde\lambda_{app}-\alpha>0.011.
\end{equation}
Hence $M^v\geq0$ outside $\mathcal R_0$ by Proposition~\ref{c:prop:supp}, and the vector form can be localized with the same nine cells as in the odd case.

For the scalar form, we show that
\begin{equation}\label{c:scalar-confinement}
V\geq0\quad\hbox{outside}\quad
\mathcal R_s:=[0.2,152]\times[-35,35].
\end{equation}
Where $R\geq200$, \eqref{c:lead1} gives $(e^{(1)}_0)_-\leq6\sigma r^2R^{-6}\leq6\sigma R^{-4}$, so $(e_0)_-\leq7\sigma R^{-4}$ by \eqref{c:Hyp}. Since $r\leq R$ and $\sigma=130$,
\begin{equation}\label{c:scalar-far}
r^2V\geq0.035-7\sigma\frac{r^2}{R^4}\geq0.035-\frac{910}{200^2}>0.
\end{equation}
The remaining region $R<200$ lies in $(0,200]\times[-200,200]$, where \eqref{c:scalar-confinement} follows from Proposition~\ref{c:prop:CAP}(iii).

Using the plateau function \eqref{c:plateau}, set
\begin{equation}\label{c:scalar-cells}
\begin{split}
&S_0:=P[0.1,0.2,15,17],\quad
 S_1:=P[15,17,29,31],\quad
 S_2:=P[29,31,152,160],\\
&V_0:=P[-40,-35,35,40],\qquad
 V_1:=P[-200,-180,180,200],\\
&\varrho_{s,0}:=S_0V_0,\qquad
 \varrho_{s,1}:=S_1V_0,\qquad
 \varrho_{s,2}:=S_2V_1.
\end{split}
\end{equation}
All three scalar cells are symmetric in $z$. Their sum equals one on $\mathcal R_s$ and is at most one everywhere, since the radial bumps have sum at most one and the vertical bumps lie in $[0,1]$.

\subsection{The weighted Neumann Laplacian of a cell}\label{c:ss:spec}

Let $\varrho=TU$ be a cell with support $Q=[r_a,r_b]\times[z_a,z_b]$. We define
\begin{equation}\label{c:cellnorm}
\langle u,v\rangle_\varrho:=\int_Q u\cdot v\,\varrho\,r^3drdz,\qquad
\|u\|_\varrho^2:=\langle u,u\rangle_\varrho,\qquad L^2(\varrho):=L^2(Q,\varrho\,r^3drdz),
\end{equation}
for scalar functions or vector fields $u,v$ on $Q$. Notice that $\varrho\,r^3=\rho_1(r)\rho_2(z)$, where
\begin{equation}\label{c:rho12}
\rho_1(r):=T(r)r^3\ \text{ on }I:=[r_a,r_b],\qquad \rho_2(z):=U(z)\ \text{ on }J:=[z_a,z_b].
\end{equation}
Therefore the quadratic form $u\mapsto\int_Q\varrho|\nabla u|^2\,r^3drdz$ is the sum of two one-dimensional forms, which we consider first.

\begin{definition}\label{c:def:weight}
A weight on $[a,b]$ is a function $\rho\in C^1([a,b])$, positive on $(a,b)$, such that at each endpoint either $\rho>0$, or $\rho$ vanishes to fifth order, in the sense that $\rho(a+\tau)=\tau^5p(\tau)$ (respectively $\rho(b-\tau)=\tau^5p(\tau)$) for small $\tau\geq0$, where $p$ is a polynomial with $p(0)>0$. We define $L^2(\rho):=L^2([a,b],\rho\,dx)$, with inner product
$\langle\cdot,\cdot\rangle_\rho$, and
\[
H_\rho:=\{\varphi\in L^2(\rho)\cap H^1_{loc}(a,b):\ \varphi'\in L^2(\rho)\},
\qquad
\mathfrak a_\rho(\varphi):=\int_a^b\rho\,(\varphi')^2dx .
\]
\end{definition}

Since $\chi$ is a polynomial with $\chi(t)=126t^5+O(t^6)$, every plateau bump vanishes to fifth order at its endpoints. Therefore the functions $\rho_1$ and $\rho_2$ defined in \eqref{c:rho12} are weights in the sense of Definition~\ref{c:def:weight}.

\begin{lemma}\label{c:lem:1dspec}
Let $\rho$ be a weight on $[a,b]$.
\begin{enumerate}[(i)]
\item The form $\mathfrak a_\rho$ is closed on $H_\rho$, and the embedding $H_\rho\hookrightarrow L^2(\rho)$ is compact.
\item There is an orthonormal basis $(\varphi_m)_{m\geq0}$ of $L^2(\rho)$, and there are numbers $0=\lambda_0<\lambda_1\leq\lambda_2\leq\ldots\to\infty$, such that $\varphi_m\in H_\rho$, $\varphi_0$ is constant, and
\begin{equation}\label{c:1dexp}
\mathfrak a_\rho(\varphi)=\sum_{m\geq0}\lambda_m\langle\varphi,\varphi_m\rangle_\rho^2
\qquad\text{for any }\varphi\in H_\rho .
\end{equation}
\item The functions $\varphi_m$ belong to $C^2([a,b])$ and satisfy
$-(\rho\varphi_m')'=\lambda_m\rho\varphi_m$ on $(a,b)$ and $\varphi_m'(a)=\varphi_m'(b)=0$.
\item The eigenvalues are simple, and $\varphi_m$ has exactly $m$ zeros in $(a,b)$.
\item If $[a,b]=[-l,l]$ and $\rho$ is even, then $\varphi_m$ is even if $m$ is even and odd if $m$ is odd.
\end{enumerate}
\end{lemma}

\begin{proof}
(i) We show first that the form $\mathfrak a_\rho$ is closed. If $(\varphi_n)$ is a Cauchy sequence with respect to the norm
$\|\varphi\|^2_{H_\rho}:=\|\varphi\|^2_\rho+\mathfrak a_\rho(\varphi)$, then
$\varphi_n\to\varphi$ and $\varphi_n'\to\psi$ in $L^2(\rho)$. The weight is bounded from below on every compact
subinterval of $(a,b)$, so $\varphi_n\to\varphi$ in
$H^1_{loc}$ and $\varphi'=\psi$. Therefore $\varphi\in H_\rho$ and
$\varphi_n\to\varphi$ in $H_\rho$.

To prove that the embedding is compact we show that there are constants $C,\eps_0>0$, depending only on $\rho$, such that
\begin{equation}\label{c:tail}
\int_a^{a+\eps}\rho\varphi^2dx+\int_{b-\eps}^b\rho\varphi^2dx\leq C\eps\,\|\varphi\|^2_{H_\rho}
\qquad\text{for any }\eps\in(0,\eps_0],\ \varphi\in H_\rho.
\end{equation}
Compactness follows from \eqref{c:tail} by a diagonal argument. Indeed, the weight is bounded from below on $[a+\eps,b-\eps]$, so by Rellich's theorem any sequence $(\varphi_n)$ with $\|\varphi_n\|_{H_\rho}\leq1$ has a subsequence that converges in $L^2([a+\eps,b-\eps],\rho\,dx)$. We take
$\eps=1/k$, $k\to\infty$, and a diagonal subsequence, which is a Cauchy sequence in $L^2(\rho)$, since by \eqref{c:tail} the contribution of the two ends to $\|\varphi_n-\varphi_{n'}\|^2_\rho$ is bounded by $4C\eps$.

We prove now \eqref{c:tail} at the endpoint $a$. If $\rho(a)>0$ then $\rho\geq c>0$ on
$[a,a+\eps_0]$, so $\varphi\in H^1(a,a+\eps_0)$ with
$\|\varphi\|_{H^1}\leq c^{-1/2}\|\varphi\|_{H_\rho}$. Therefore, using the one-dimensional Sobolev embedding,
\[
\int_a^{a+\eps}\rho\varphi^2\leq\eps\|\rho\|_\infty\|\varphi\|^2_{L^\infty(a,a+\eps_0)}
\leq C\eps\|\varphi\|^2_{H_\rho}.
\]
Assume now that $\rho$ vanishes to fifth order at $a$. Then there are $c_1,c_2,\tau_0>0$ such that $c_1\tau^5\leq\rho(a+\tau)\leq c_2\tau^5$ for $0\leq\tau\leq\tau_0$. Let
$\eps\leq\eps_0:=\tau_0/2$, $x_0:=a+\eps$, and let $\varphi\in H_\rho$, which is continuous on $(a,b)$. For $x\in(a,x_0]$ we have $\varphi(x)=\varphi(x_0)-\int_x^{x_0}\varphi'$, so, using the Cauchy--Schwarz inequality,
\[
\varphi(x)^2\leq2\varphi(x_0)^2+2\Big(\int_x^{x_0}\frac{dy}{\rho(y)}\Big)\int_x^{x_0}\rho(\varphi')^2
\leq2\varphi(x_0)^2+\frac{\mathfrak a_\rho(\varphi)}{2c_1(x-a)^4}.
\]
We multiply by $\rho(x)\leq c_2(x-a)^5$ and integrate over $(a,x_0)$, thus
\[
\int_a^{a+\eps}\rho\varphi^2\leq\frac{c_2\eps^6}{3}\varphi(x_0)^2+\frac{c_2\eps^2}{4c_1}\mathfrak a_\rho(\varphi).
\]
It remains to bound $\varphi(x_0)^2$. Since $\rho\geq c_1\eps^5$ on $[x_0,x_0+\tau_0/2]\subset(a,a+\tau_0]$, there is a point $x_1$ in this interval such that
$\varphi(x_1)^2\leq\frac{2}{\tau_0}\int_{x_0}^{x_0+\tau_0/2}\varphi^2\leq\frac{2}{\tau_0c_1\eps^5}\|\varphi\|^2_\rho$. Therefore
\begin{equation}\label{c:point}
\varphi(x_0)^2\leq2\varphi(x_1)^2+2(x_1-x_0)\int_{x_0}^{x_1}(\varphi')^2
\leq\frac{4}{\tau_0c_1\eps^5}\|\varphi\|^2_\rho+\frac{\tau_0}{c_1\eps^5}\mathfrak a_\rho(\varphi),
\end{equation}
so $\frac{c_2\eps^6}{3}\varphi(x_0)^2\leq C\eps\|\varphi\|^2_{H_\rho}$, and
\eqref{c:tail} follows. The endpoint $b$ is similar.

(ii) By (i), the nonnegative self-adjoint operator $L_\rho$ associated to the closed form $\mathfrak a_\rho$ has compact resolvent. Therefore it has an orthonormal basis of eigenfunctions $\varphi_m\in D(L_\rho)\subset H_\rho$, with eigenvalues
$\lambda_m\to\infty$, and
\[
\mathfrak a_\rho(\varphi)=\|L_\rho^{1/2}\varphi\|^2_\rho=\sum_m\lambda_m\langle\varphi,\varphi_m\rangle^2_\rho
\qquad\text{for }\varphi\in H_\rho=D(L_\rho^{1/2}).
\]
The kernel of $\mathfrak a_\rho$ consists of the functions $\varphi\in H_\rho$ with $\varphi'=0$ a.e. in $(a,b)$, i.e. the constant functions, which belong to $L^2(\rho)$ since $\rho$ is bounded. Thus $\lambda_0=0$ is a simple eigenvalue and $\lambda_1>0$.

(iii) The weak form of the eigenvalue equation gives $-(\rho\varphi_m')'=\lambda_m\rho\varphi_m$ in the sense of distributions on $(a,b)$. Hence $\rho\varphi_m'\in H^1_{loc}$ is continuous, so $\varphi_m\in C^1$, and then $\rho\varphi_m'\in C^1$, so $\varphi_m\in C^2((a,b))$, since $\rho\in C^1$ is positive on $(a,b)$. At an endpoint where $\rho>0$, the eigenfunctions are $C^2$ up to the endpoint and satisfy the Neumann condition $\varphi_m'=0$, since $\int\rho\varphi_m'\chi'=\lambda_m\int\rho\varphi_m\chi$ for all $\chi\in H_\rho$. At an endpoint where $\rho$ vanishes to fifth order, say $a$, the equation reads, in the variable $\tau=x-a$,
\[
\varphi''+\Big(\frac5\tau+\frac{p'(\tau)}{p(\tau)}\Big)\varphi'+\lambda\varphi=0,
\]
which has a regular singular point at $\tau=0$ with indicial exponents $0$ and $-4$. The solutions with exponent $-4$ are not in $L^2(\rho)$ near $a$, since $\int_0\tau^5\tau^{-8}d\tau=\infty$, so $\varphi_m$ is the solution with exponent $0$, which is given by a convergent power series in $\tau$ (as in \eqref{eq:decomp}); in particular $\varphi_m\in C^2([a,b])$ and $\varphi_m'(a)=0$.

(iv) At an endpoint where $\rho$ vanishes to fifth order, the endpoint is of limit-point type by (iii), and it is non-oscillatory since the indicial exponents are real. Therefore $L_\rho$ is the self-adjoint realization of $-\rho^{-1}(\rho\varphi')'$ with Neumann conditions at the endpoints where $\rho>0$, and (iv) follows from the oscillation theorem for Sturm--Liouville operators with regular or non-oscillatory limit-point endpoints, see \cite[Chapter 10]{Zettl}.

(v) The reflection $\varphi\mapsto\varphi(-\cdot)$ is a unitary operator on $L^2(\rho)$ that preserves $H_\rho$ and $\mathfrak a_\rho$, so it commutes with $L_\rho$. By (iv) every eigenspace is one-dimensional, so each $\varphi_m$ is even or odd. If $\varphi_m$ is even then $\varphi_m(0)\neq0$, since otherwise $\varphi_m(0)=\varphi_m'(0)=0$ and $\varphi_m\equiv0$, because $\rho(0)>0$. Thus the zeros of an even eigenfunction come in pairs $\pm x$, while the zeros of an odd eigenfunction are $0$ and pairs $\pm x$. The claim follows from (iv).
\end{proof}

Let $(\lambda_m,\varphi_m)_{m\geq0}$ and $(\mu_k,\psi_k)_{k\geq0}$ denote the eigenpairs given by Lemma~\ref{c:lem:1dspec} for the weights $\rho_1$ on $I$ and $\rho_2$ on $J$, with the signs fixed by $\varphi_m(r_a)>0$ and $\psi_k(z_a)>0$. Note that this is possible since $\varphi_m(r_a)\neq0$; indeed, by the proof of Lemma~\ref{c:lem:1dspec}(iii), $\varphi_m$ is the solution with exponent $0$ at $r_a$, which is determined by its value at $r_a$. Let
\begin{equation}\label{c:modes}
e_{mk}(r,z):=\varphi_m(r)\psi_k(z),\qquad \Lambda_{mk}:=\lambda_m+\mu_k,
\qquad c_{mk}(u):=\langle u,e_{mk}\rangle_\varrho .
\end{equation}
For a symmetric cell the function $\psi_k$ has the parity of $k$, by Lemma~\ref{c:lem:1dspec}(v), and we call the mode $(m,k)$ even or odd according to the parity of $k$.

\begin{proposition}\label{c:prop:spec}
Let $\varrho=TU$ be a cell.
\begin{enumerate}[(i)]
\item The family $\{e_{mk}\}_{m,k\geq0}$ is an orthonormal basis of $L^2(\varrho)$, and
\begin{equation}\label{c:sep2}
\int_Q\varrho|\nabla u|^2\,r^3drdz=\sum_{m,k\geq0}\Lambda_{mk}\,c_{mk}(u)^2\qquad\text{for any }u\in C^1(\bar Q).
\end{equation}
\item Let $\Sigma$ be a finite set of modes, and let
\[
P_\Sigma u:=\sum_{(m,k)\in\Sigma}c_{mk}(u)e_{mk},\qquad
\Lambda_*(\Sigma):=\min\{\Lambda_{mk}:(m,k)\notin\Sigma\}.
\]
Then
\begin{equation}\label{c:sep3}
\int_Q\varrho|\nabla u|^2\,r^3drdz\geq\sum_{(m,k)\in\Sigma}\Lambda_{mk}\,c_{mk}(u)^2
+\Lambda_*(\Sigma)\,\|u-P_\Sigma u\|^2_\varrho .
\end{equation}
\item Assume that the cell is symmetric. If $u$ is odd (respectively even) in $z$, then $c_{mk}(u)=0$ for every even (respectively odd) mode. If, in addition, $\Sigma$ consists of odd (respectively even) modes, then \eqref{c:sep3} holds with $\Lambda_*(\Sigma)$ replaced by the minimum of $\Lambda_{mk}$ over the odd (respectively even) modes not in $\Sigma$.
\end{enumerate}
\end{proposition}

\begin{proof}
(i) Since $\varrho\,r^3drdz=\rho_1(r)\rho_2(z)\,drdz$, we have
$L^2(\varrho)=L^2(I,\rho_1dr)\otimes L^2(J,\rho_2dz)$, and the products of the
elements of two orthonormal bases form an orthonormal basis of the tensor
product. Let $u\in C^1(\bar Q)$. For any fixed $z$ the function $u(\cdot,z)$ belongs to
$C^1(I)\subset H_{\rho_1}$, so, using \eqref{c:1dexp},
\[
\int_I\rho_1|\partial_ru|^2dr=\sum_m\lambda_m\,\al_m(z)^2,\qquad
\al_m(z):=\langle u(\cdot,z),\varphi_m\rangle_{\rho_1}.
\]
We integrate in $z$ against $\rho_2dz$ and use the monotone convergence theorem and Parseval's identity in $L^2(J,\rho_2dz)$, thus
\[
\int_Q\varrho|\partial_ru|^2\,r^3drdz=\sum_m\lambda_m\int_J\rho_2\al_m^2dz
=\sum_m\lambda_m\sum_k\langle\al_m,\psi_k\rangle^2_{\rho_2}=\sum_{m,k}\lambda_mc_{mk}(u)^2,
\]
where we used Fubini's theorem in the last identity. An analogous statement holds for $\partial_zu$, and \eqref{c:sep2} follows.

(ii) In \eqref{c:sep2} we bound $\Lambda_{mk}\geq\Lambda_*(\Sigma)$ for $(m,k)\notin\Sigma$ and use Parseval's identity.

(iii) We have $c_{mk}(u)=\int_I\rho_1\varphi_m\big(\int_J\rho_2u\psi_kdz\big)dr$, and the inner integral vanishes if $u(r,\cdot)$ and $\psi_k$ have opposite parities. Thus the sum in \eqref{c:sep2} only contains the modes with the same parity as $u$, and the claim follows as in (ii).
\end{proof}

\subsection{A lower bound on a cell}\label{c:ss:cellthm}

The cell bound below is applied separately to three forms: the vector forms $\mathcal E$ and $\widetilde{\mathcal E}_v$ of the odd and even problems, and the scalar form $\widetilde{\mathcal E}_s$ of the even problem. We therefore state it for $s\in\{1,2\}$ components and a diffusion coefficient $d>0$. Fix a parity $\varsigma\in\{-1,1\}$ and set
\[
J_\varsigma:=\begin{cases}
\varsigma,&s=1,\\
\operatorname{diag}(\varsigma,-\varsigma),&s=2.
\end{cases}
\]
If $g(r,-z)=\varsigma g(r,z)$, then $u=g$ (for $s=1$) and $u=\nabla g$ (for $s=2$) satisfy $u(r,-z)=J_\varsigma u(r,z)$, and we refer to $J_\varsigma$ as the component parities of $u$.

\begin{theorem}\label{c:thm:refined}
Let $\varrho=TU$ be a cell supported on $Q$, and let $M$ be a bounded
measurable function on $Q$ with values in the symmetric $s\times s$ matrices, such that $M\geq-WI_s$ for some $W\geq0$. For $1\leq a\leq s$, let $\Sigma_a$ be a finite set of modes. If the cell is symmetric, we assume further that all modes in $\Sigma_a$ have parity $(J_\varsigma)_{aa}$. Let $\Lambda_*>0$ be a lower bound for $\Lambda_*(\Sigma_a)$, $1\leq a\leq s$, or, if the cell is symmetric, for the corresponding minima in Proposition~\ref{c:prop:spec}(iii). For a mode $j=(m,k)$, write $e_j=e_{mk}$ and $\Lambda_j=\Lambda_{mk}$. Assume $d\Lambda_*>W$, and define, for $j\in\Sigma_a$ and $l\in\Sigma_b$,
\begin{equation}\label{c:NGG}
\begin{split}
&G_{(j,a),(l,b)} :=\int_Q\varrho M_{ab}e_je_l\,r^3drdz,\qquad (G_2)_{(j,a),(l,b)} :=\int_Q\varrho(M^2)_{ab}e_je_l\,r^3drdz,\\
&N :=d\Lambda+G-\frac{G_2-G^2}{d\Lambda_*-W},
\end{split}
\end{equation}
where $\Lambda$ is diagonal with $\Lambda_{(j,a),(j,a)}=\Lambda_j$. For $u\in C^1(\bar Q)^s$, with the prescribed component parities on a symmetric cell, let $c_{(j,a)}:=\langle u_a,e_j\rangle_\varrho$. Then
\begin{equation}\label{c:ref1}
\mathcal E_\varrho(u):=\int_Q\varrho\big(d|\nabla u|^2+u^TMu\big)r^3drdz
\geq c^TNc.
\end{equation}
Moreover, $G_2-G^2\geq0$ and $|c|\leq\|u\|_\varrho$.
\end{theorem}

\begin{proof}
Let $P$ be the orthogonal projection in $L^2(\varrho)^s$ onto
$\mathcal V=\operatorname{span}\{e_j\mathbf e_a:j\in\Sigma_a\}$, and set
$v=(I-P)u$. Proposition~\ref{c:prop:spec}, applied componentwise, gives
\[
d\int_Q\varrho|\nabla u|^2r^3drdz
\geq d c^T\Lambda c+d\Lambda_*\|v\|_\varrho^2.
\]
Also $\|Pu\|_\varrho^2=|c|^2$, $\langle Pu,MPu\rangle_\varrho=c^TGc$, and
\[
\|(I-P)MPu\|_\varrho^2
=\|MPu\|_\varrho^2-\|PMPu\|_\varrho^2=c^T(G_2-G^2)c.
\]
This proves the asserted positivity. Writing $b=(I-P)MPu$, we obtain
\[
\mathcal E_\varrho(u) \geq d c^T\Lambda c+c^TGc +(d\Lambda_*-W)\|v\|_\varrho^2+2\langle v,b\rangle_\varrho \geq d c^T\Lambda c+c^TGc-\frac{\|b\|_\varrho^2}{d\Lambda_*-W},
\]
which is \eqref{c:ref1}.
\end{proof}

The following version of Gershgorin's theorem converts the finite matrix bound to explicit functionals.

\begin{lemma}\label{c:lem:gersh}
Let $N$ be a real symmetric matrix, and let $X$ have columns $\xi_1,\ldots,\xi_K$, with $\tau:=\|I-XX^T\|\leq1/3$. Define
\begin{equation}\label{c:gersh1}
D:=X^TNX,\qquad g_l:=D_{ll}-\sum_{k\ne l}|D_{lk}|.
\end{equation}
Then $|\xi_l|^2\leq1+\tau$, and
\begin{equation}\label{c:gersh2}
c^TNc\geq\sum_lg_l(\xi_l\cdot c)^2-3\tau\|N\|\,|c|^2
\geq-\sum_{l:g_l<0}|g_l|(\xi_l\cdot c)^2-3\tau\|N\|\,|c|^2.
\end{equation}
\end{lemma}

\begin{proof}
Put $b=X^Tc$ and $e=c-Xb=(I-XX^T)c$. Symmetry gives
$b^TDb\geq\sum_lg_lb_l^2$. Since $\|X\|^2\leq1+\tau$,
$|Xb|\leq(1+\tau)|c|$ and $|e|\leq\tau|c|$. Consequently,
\[
c^TNc=b^TDb+2(Xb)^TNe+e^TNe
\geq\sum_lg_lb_l^2-(2+3\tau)\tau\|N\|\,|c|^2,
\]
and $2+3\tau\leq3$. The bound on $|\xi_l|^2$ follows from
$|\xi_l|\leq\|X\|$.
\end{proof}

\begin{corollary}\label{c:cor:cell}
Under the assumptions of Theorem~\ref{c:thm:refined} and Lemma~\ref{c:lem:gersh}, choose a set $L$ containing every index with
$g_l<0$, and let
\[
\nu_l\geq\max(0,-g_l)\quad(l\in L),\qquad
\varepsilon:=3\tau\|N\|,
\]
and define
\begin{equation}\label{c:Xi}
\Xi_l:=\sum_{a=1}^s\sum_{j\in\Sigma_a}
\xi_{l,(j,a)}e_j\mathbf e_a,\qquad l\in L.
\end{equation}
Then $\|\Xi_l\|_\varrho^2=|\xi_l|^2\leq1+\tau$, the functions $\Xi_l$ have the prescribed component parities on a symmetric cell, and
\begin{equation}\label{c:cell2}
\mathcal E_\varrho(u)
\geq-\sum_{l\in L}\nu_l\langle u,\Xi_l\rangle_\varrho^2
-\varepsilon\|u\|_\varrho^2.
\end{equation}
\end{corollary}

\begin{proof}
The modes are orthonormal, so $\xi_l\cdot c=\langle u,\Xi_l\rangle_\varrho$ and $\|\Xi_l\|_\varrho=|\xi_l|$. Apply \eqref{c:ref1} and \eqref{c:gersh2}, using $|c|\leq\|u\|_\varrho$.
\end{proof}

\subsection{From functionals to the functions $q\in Z'$}\label{c:ss:fun}

Write $\Pi_\varsigma F(r,z):=\tfrac12(F(r,z)+\varsigma F(r,-z))$.
This parity projection is a contraction on $Z'$ and on $L^2$ with any
weight even in $z$.

\begin{lemma}\label{c:lem:q}
Let $\varrho$ be a cell supported on $Q$, and let $\Xi$ be Lipschitz on $\bar Q$, with $s$ components. Extend $\varrho\Xi$ by zero. For $s=2$ set $V=\varrho(\Xi^1\hat e_r+\Xi^2\hat e_z)$ on $\R^5$ and define
\begin{equation}\label{c:q1}
q_{\varrho,\Xi}:=
\begin{cases}
\Pi_\varsigma(\varrho\Xi),&s=1,\\
-\Pi_\varsigma\!\left[\partial_r(\varrho\Xi^1)
 +\dfrac3r\varrho\Xi^1+\partial_z(\varrho\Xi^2)\right],&s=2.
\end{cases}
\end{equation}
These are bounded, compactly supported functions in $Z'$ of parity $\varsigma$. For $s=2$,
\begin{equation}\label{c:q2}
\|q_{\varrho,\Xi}\|_{Z'}
\leq\|\mathcal P_\varsigma V\|_{L^2(\R^5)}
\leq\|V\|_{L^2(\R^5)}
=\left(2\pi^2\int_Q\varrho^2|\Xi|^2r^3drdz\right)^{1/2},
\end{equation}
where $\mathcal P_\varsigma$ projects the radial and vertical components onto parities $\varsigma$ and $-\varsigma$, respectively. For $s=1$, put
\[
F(r,z):=r^{-3}\int_0^r t^3q_{\varrho,\Xi}(t,z)\,dt.
\]
Then $\|q_{\varrho,\Xi}\|_{Z'}\leq\|F\hat e_r\|_{L^2(\R^5)}$. In either case, for a smooth axially symmetric $g$ of parity $\varsigma$, with $u=g$ if $s=1$ and $u=\nabla g$ if $s=2$,
\begin{equation}\label{c:q3}
\langle q_{\varrho,\Xi},g\rangle
=2\pi^2\langle u,\Xi\rangle_\varrho.
\end{equation}
If $\varrho^-(r,z)=\varrho(r,-z)$ and
$\Xi^-(r,z)=J_\varsigma\Xi(r,-z)$, then
$q_{\varrho^-,\Xi^-}=q_{\varrho,\Xi}$ and
$\langle u,\Xi^-\rangle_{\varrho^-}=\langle u,\Xi\rangle_\varrho$.
\end{lemma}

\begin{proof}
The plateaux vanish at the boundary, and their radial supports avoid the axis, so the vector field $V$ is Lipschitz and compactly supported. Moreover, $\Pi_\varsigma\operatorname{div}V
=\operatorname{div}\mathcal P_\varsigma V$. The Fourier multiplier estimate
\[
\||\nabla|^{-1}\langle\nabla\rangle^{-1}\operatorname{div}H\|_2
\leq\|H\|_2
\]
gives \eqref{c:q2}. For $s=1$, $\operatorname{div}(F\hat e_r)=q_{\varrho,\Xi}$ in distributions. The field $F\hat e_r$ is in $L^2$: it vanishes for $r<r_a$, and for $r>r_b$ it is $r^{-3}$ times a compactly supported bounded function of $z$. The same multiplier estimate applies. Identity \eqref{c:q3} follows from the self-adjointness of $\Pi_\varsigma$ and, for vector fields, integration by parts. The reflection identities follow by $z\mapsto-z$ and the component parities of $u$.
\end{proof}

\begin{lemma}\label{c:lem:mirror}
Suppose that $M$ has the reflection symmetry
\begin{equation}\label{c:Msym}
M(r,-z)=J_\varsigma M(r,z)J_\varsigma,
\end{equation}
and let $\varrho^-(r,z)=\varrho(r,-z)$ be the reflected cell. Let $N$ be the
matrix \eqref{c:NGG} of $\varrho$, and $N^-$ the matrix \eqref{c:NGG} of
$\varrho^-$ computed with the reflected modes $e_j^-(r,z)=e_j(r,-z)$ and the same
$d,W,\Lambda_*$ and index sets. Let $\mathcal D$ be the
diagonal matrix with entry $(J_\varsigma)_{aa}$ at index $(j,a)$. Then, for
every $K\times K$ matrix $X$,
\[
N^-=\mathcal DN\mathcal D,\qquad
(\mathcal DX)^TN^-(\mathcal DX)=X^TNX,\qquad
\|I-(\mathcal DX)(\mathcal DX)^T\|=\|I-XX^T\|.
\]
Reflected cells have the same $g_l$, $\nu_l$ and $\varepsilon$, and the associated fields are $\Xi_l^-(r,z)=J_\varsigma\Xi_l(r,-z)$.
\end{lemma}

\begin{proof}
Changing $z$ to $-z$ in \eqref{c:NGG} gives $G^-=\mathcal DG\mathcal D$ and $G_2^-=\mathcal DG_2\mathcal D$. The diagonal matrix $\Lambda$ commutes with $\mathcal D$ and $\mathcal D^2=I$, proving the assertions.
\end{proof}

\subsection{Assembling the cells}\label{c:ss:assembly}

We now combine the cell bounds into a global inequality.

\begin{theorem}\label{c:thm:assembly}
Let $d>0$, $s\in\{1,2\}$, and let $M$ be a symmetric matrix field on $\R^2_+$ satisfying \eqref{c:Msym}. Let $\{\varrho_i:i\in\mathcal C\}$ be a finite family of cells invariant under reflection, such that
\[
0\leq\sum_i\varrho_i\leq1,\qquad
M\geq0\quad\hbox{on }\left\{\sum_i\varrho_i<1\right\}.
\]
Assume that for each $i\in\mathcal C$ the hypotheses of Corollary~\ref{c:cor:cell} hold on the cell $\varrho_i$, with $M_\eta=M-\eta I_s$ in place of $M$ and with component parities $J_\varsigma$, and denote by $N_i,\tau_i,L_i,\nu_{i,l},\Xi_{i,l}$ the corresponding data. Assume also that the data of the two cells of a reflected pair are related as in Lemma~\ref{c:lem:mirror}, and that
\begin{equation}\label{c:loss}
\varepsilon_i=3\tau_i\|N_i\|\leq\eta\qquad(i\in\mathcal C).
\end{equation}
Choose one representative from each reflected pair, and retain every symmetric cell. For $l\in L_i$, set
\begin{equation}\label{c:qref}
q_{(i,l)}:=\kappa_i\sqrt{\frac{\nu_{i,l}}{2\pi^2}}\,
q_{\varrho_i,\Xi_{i,l}},\qquad
\kappa_i:=\begin{cases}1,&\varrho_i\text{ symmetric},\\
\sqrt2,&\varrho_i\text{ in a reflected pair}.
\end{cases}
\end{equation}
Then the functions $q_{(i,l)}$ belong to $Z'$ and have parity $\varsigma$, and for $u=g$ when $s=1$, or
$u=\nabla g$ when $s=2$, with $g$ smooth and of parity $\varsigma$,
\begin{equation}\label{c:ass1}
\mathcal E_{d,M}(u):=\int_{\R^2_+}
\big(d|\nabla u|^2+u^TMu\big)r^3drdz
\geq-\frac1{2\pi^2}\sum_{(i,l)}\langle q_{(i,l)},g\rangle^2.
\end{equation}
\end{theorem}

\begin{proof}
Since $M\geq0$ where $\sum_i\varrho_i<1$, the weight $1-\sum_i\varrho_i\geq0$ contributes a nonnegative term to $\mathcal E_{d,M}(u)$. Writing $M=M_\eta+\eta I_s$ inside each cell therefore gives
\begin{equation}\label{c:loc}
\mathcal E_{d,M}(u)
\geq\sum_{i\in\mathcal C}\left[
\int\varrho_i\big(d|\nabla u|^2+u^TM_\eta u\big)r^3drdz
+\eta\|u\|_{\varrho_i}^2\right].
\end{equation}
By \eqref{c:cell2} and \eqref{c:loss}, the $i$th summand is at least $-\sum_{l\in L_i}\nu_{i,l}\langle u,\Xi_{i,l}\rangle_{\varrho_i}^2$. Further, by Lemmas~\ref{c:lem:q} and~\ref{c:lem:mirror}, the two members of a reflected pair have equal functionals. Thus \eqref{c:q3} and \eqref{c:qref} give, after combining each pair,
\[
\kappa_i^2\nu_{i,l}\langle u,\Xi_{i,l}\rangle_{\varrho_i}^2
=\frac1{2\pi^2}\langle q_{(i,l)},g\rangle^2.
\]
Summing proves \eqref{c:ass1}.
\end{proof}

We now bound the norms of the $q_{(i,l)}$. On a symmetric cell the fields have the prescribed parity. On a positive-$z$ cell and its disjoint reflection, the squared norm of the parity projection is one half of the unprojected squared norm, canceling $\kappa_i^2=2$. Consequently, for vector functionals,
\begin{equation}\label{c:vector-norm}
\|q_{(i,l)}\|_{Z'}^2
\leq\nu_{i,l}\int_{Q_i}\varrho_i^2|\Xi_{i,l}|^2r^3drdz
\leq\nu_{i,l}(1+\tau_i).
\end{equation}
For a scalar functional define $H_{i,l}(r,z):=\int_0^r t^3\varrho_i(t,z)\Xi_{i,l}(t,z)\,dt$.
The radial primitive in Lemma~\ref{c:lem:q} yields the sharper bound
\begin{equation}\label{c:scalar-norm}
\|q_{(i,l)}\|_{Z'}^2
\leq\nu_{i,l}\int_{z_a}^{z_b}
\left[\int_0^{r_b}\frac{|H_{i,l}(r,z)|^2}{r^3}\,dr
+\frac{|H_{i,l}(r_b,z)|^2}{2r_b^2}\right]dz.
\end{equation}
The last term is the contribution of $r>r_b$, where $H_{i,l}(r,z)=H_{i,l}(r_b,z)$.

\subsection{The functions and the two coercivity inequalities}\label{c:ss:q}

For each row $i$ of Table~\ref{c:input-cells} and $l\in L_i$, let $q_{(i,l)}$ be the function \eqref{c:qref}, with $\Xi_{i,l}$ given by \eqref{c:Xi}. Explicitly,
\begin{equation}\label{c:qfinal}
q_{(i,l)}=\kappa_i\sqrt{\frac{\nu_{i,l}}{2\pi^2}}
\begin{cases}
-\Pi_{\varsigma_i}\!\left[\partial_r(\varrho_i\Xi^1_{i,l})
 +\dfrac3r\varrho_i\Xi^1_{i,l}
 +\partial_z(\varrho_i\Xi^2_{i,l})\right],&s_i=2,\\[3pt]
\Pi_{\varsigma_i}(\varrho_i\Xi_{i,l}),&s_i=1.
\end{cases}
\end{equation}
Let $q_1,\ldots,q_{26}$ be an enumeration of the functions $q_{(i,l)}$ for the first block of Table~\ref{c:input-cells}, and $\widetilde q_1,\ldots,\widetilde q_{32}$ an enumeration of those for the second and third blocks. These are the key functions $q_j$ and $\widetilde q_j$ used in \eqref{robu8} and \eqref{robu32}.

\begin{proof}[Proof of Theorem~\ref{c:thm:main}]
For the odd problem, apply Theorem~\ref{c:thm:assembly} with $d=1-\gamma_0$, $M=M^o$, $\eta=\eta_v:=0.01$, the odd component parities, and the cells and data of the first block of Table~\ref{c:input-cells}. Lemma~\ref{c:lem:pu} and Proposition~\ref{c:prop:supp} give the partition and confinement assumptions. Since $\Psi_5$ is odd in $z$, the functions $e_{rr}$, $e_{zz}$, $e_0$ are even and $a_{12}$ is odd in $z$. Since $\zeta$ is even, $M^o$ satisfies \eqref{c:Msym}, and so do $M^v$ and $V$. Proposition~\ref{c:prop:CAP}(iv)--(v) gives the hypotheses of Theorem~\ref{c:thm:refined}, except $d\Lambda_{*,i}>W_i$, which can be read off from Table~\ref{c:input-cells}. Part (vi) of Proposition~\ref{c:prop:CAP} gives the Gershgorin inequalities of Corollary~\ref{c:cor:cell} and $3\tau_i\|N_i\|\leq9\times10^{-11}<\eta_i$ in \eqref{c:loss}. Therefore \eqref{c:Q}, Proposition~\ref{c:prop:vec}, and \eqref{c:ass1} imply
\[
\langle\mathcal S^{(5)}_0f,\Delta^{-1}f\rangle
+\sum_{j=1}^{26}\langle q_j,\Delta^{-1}f\rangle^2
\geq\gamma_0\|f\|_Z^2.
\]

For the even problem, apply Theorem~\ref{c:thm:assembly} twice: with $(d,M,\eta)=(1-\gamma_0,M^v,\eta_v)$ and the second block of Table~\ref{c:input-cells}, and with $(d,M,\eta)=(\beta,V,\eta_s)$, $\eta_s:=10^{-8}$, and the third block. The partition and confinement assumptions were verified in Subsection~\ref{c:ss:even}, and the remaining hypotheses follow from Proposition~\ref{c:prop:CAP}(iv)--(vi) as in the odd case. Adding the two inequalities and using \eqref{c:even-Hardy-split} yields
\[
\langle\mathcal S^{(5)}_{\widetilde\lambda_{app}}h,\Delta^{-1}h\rangle
+\sum_{j=1}^{32}\langle\widetilde q_j,\Delta^{-1}h\rangle^2
\geq\gamma_0\|h\|_Z^2.
\]

In both cases, \eqref{c:profile-loss} gives $\gamma_0-1.001\times10^{-10}>0.1$, proving \eqref{c:main1}. The calculation first applies to smooth compactly supported functions of the specified parity and extends to $Z_0$ by density. Finally the norm bounds \eqref{c:vector-norm}--\eqref{c:scalar-norm} and Proposition~\ref{c:prop:CAP}(vii) give \eqref{c:qbounds}.
\end{proof}

\subsection{Computer-assisted inputs}\label{c:ss:CAP}

The data in \cite{proof_code} specify the following objects for Proposition~\ref{c:prop:CAP}.
\begin{enumerate}
\item The radial partition $0=r_0<r_1<\cdots<r_{K_b}=200$, with $K_b=2482$, and the numbers $L_{1,k},L_{2,k},h_{1,k},h_{2,k}$ of \eqref{c:boxes}.
\item For each row $i$ of Table~\ref{c:input-cells}, the retained mode sets $\Sigma_{i,a}$, $1\leq a\leq s_i$, and a real square matrix $X_i$ of order $K_i=\sum_a|\Sigma_{i,a}|$. Its rows are indexed by $(j,a)$ with $j\in\Sigma_{i,a}$ and its columns are denoted by $\xi_{i,l}$. The set $L_i\subseteq\{1,\ldots,K_i\}$ selects $n_i$ columns with weights $\nu_{i,l}>0$.
\end{enumerate}
The selected columns and their weights fix the sources in \eqref{c:qfinal}. For each source, the certified data in \cite{proof_code} specify the response and its auxiliary potential in \eqref{resp:representation} in the same normalization as the source, together with the corresponding scalar $\widetilde\mu_j$ in the even problem.

For each row $i$, let $\lambda_m,\mu_k$ be the eigenvalues given by Lemma~\ref{c:lem:1dspec} for the weights $\rho_1,\rho_2$ of the cell $\varrho_i=TU$. A mode $(m,k)$ is admissible for the component $a$ if the cell is not symmetric or $(-1)^k=(J_{\varsigma_i})_{aa}$. Let $G_i$, $G_{2,i}$, $\Lambda_i$ and $N_i$ be defined by \eqref{c:NGG} with $\varrho=\varrho_i$, $M=M_i-\eta_iI_{s_i}$, $\Sigma_a=\Sigma_{i,a}$, $d=d_i$, $W=W_i$ and $\Lambda_*=\Lambda_{*,i}$, and let $D_i$ and $g_{i,l}$ be defined by \eqref{c:gersh1} with $N=N_i$ and $X=X_i$. For $l\in L_i$, with $\Xi_{i,l}$ given by \eqref{c:Xi}, let $B_{i,l}^2$ be the middle term of \eqref{c:vector-norm} for the vector rows and the right side of \eqref{c:scalar-norm} for the scalar rows.

\begin{table}[htbp]
\centering
\caption{The rows of the cell construction. The cells $T_kU_1$ represent the reflected pairs $\{T_kU_1,T_kU_1^-\}$. All displayed thresholds are fixed decimal constants.}
\label{c:input-cells}
\small
\begin{tabular}{lrrr}
\toprule
cell $\varrho_i$ & $W_i$ & $\Lambda_{*,i}$ & $n_i$\\
\midrule
\multicolumn{4}{l}{$\varsigma_i=-1$, $s_i=2$, $M_i=M^o$,
 $d_i=1-\gamma_0$, $\eta_i=\eta_v$}\\[2pt]
$T_0U_0$ & 3.76 & 12.5455 & 9\\
$T_0U_1$ & 0.504 & 1.7237 & 1\\
$T_1U_0$ & 4.96 & 16.6307 & 11\\
$T_1U_1$ & 0.489 & 1.6562 & 3\\
$T_2U_0$ & 0.362 & 1.3053 & 1\\
$T_2U_1$ & 0.235 & 0.854 & 1\\
\midrule
\multicolumn{4}{l}{$\varsigma_i=1$, $s_i=2$, $M_i=M^v$,
 $d_i=1-\gamma_0$, $\eta_i=\eta_v$}\\[2pt]
$T_0U_0$ & 3.71 & 12.3774608652 & 12\\
$T_0U_1$ & 0.494 & 1.7237277896 & 0\\
$T_1U_0$ & 4.72 & 15.7151282883 & 9\\
$T_1U_1$ & 0.477 & 1.623936892 & 2\\
$T_2U_0$ & 0.285 & 1.0866964506 & 0\\
$T_2U_1$ & 0.184 & 0.6501580511 & 0\\
\midrule
\multicolumn{4}{l}{$\varsigma_i=1$, $s_i=1$, $M_i=V$,
 $d_i=\beta$, $\eta_i=\eta_s$}\\[2pt]
$S_0V_0$ & 0.706 & 11.2020998964 & 6\\
$S_1V_0$ & 0.706 & 11.2110651424 & 3\\
$S_2V_1$ & 0.00238 & 0.0378391209 & 0\\
\bottomrule
\end{tabular}
\end{table}

\begin{proposition}\label{c:prop:CAP}
The following hold.
\begin{enumerate}[(i)]
\item For $(r,z)\in\R^2_+$ with $R\geq200$,
\begin{equation}\label{c:Hyp}
\begin{split}
|e_{rr}-e^{(1)}_{rr}|&\leq0.05\,\sigma R^{-2},\qquad
|e_{zz}-e^{(1)}_{zz}|\leq0.05\,\sigma R^{-2},\\
|a_{12}-a^{(1)}_{12}|&\leq0.05\,\sigma R^{-2},\qquad
|e_0-e^{(1)}_0|\leq\sigma R^{-4}.
\end{split}
\end{equation}

\item For $1\leq k\leq K_b$, $0\leq L_{1,k}\leq L_{2,k}<200$,
$h_{1,k},h_{2,k}\geq0$, and for $r\in[r_{k-1},r_k]$, $r>0$,
\begin{equation}\label{c:boxes}
\begin{aligned}
e_0(r,z)&\geq-h_{1,k}-h_{2,k} &&\text{if } |z|\leq L_{1,k},\\
e_0(r,z)&\geq-h_{2,k} &&\text{if } L_{1,k}\leq|z|\leq L_{2,k},\\
e_0(r,z)&\geq0 &&\text{if } L_{2,k}\leq|z|\leq200.
\end{aligned}
\end{equation}

\item Let $\mathcal B:=(0,200]\times[-200,200]$. Then $M^o\geq0$ on $\mathcal B\setminus\mathcal R_0$ and $V\geq0$ on $\mathcal B\setminus\mathcal R_s$.

\item For every row $i$, $M_i-\eta_iI_{s_i}\geq-W_iI_{s_i}$ on the support of $\varrho_i$.

\item For every row $i$ and $1\leq a\leq s_i$, the modes in $\Sigma_{i,a}$ are admissible for the component $a$, and
\begin{equation}\label{c:input-modes}
\lambda_m+\mu_k\geq\Lambda_{*,i}
\qquad\text{for every admissible mode }(m,k)\notin\Sigma_{i,a}.
\end{equation}

\item For every row $i$, with $\|\cdot\|$ the spectral norm,
\[
\|I-X_iX_i^T\|\leq10^{-12},\qquad \|N_i\|\leq30,\qquad
g_{i,l}\geq-\nu_{i,l}\quad(l\in L_i),\qquad
g_{i,l}\geq0\quad(l\notin L_i).
\]

\item With $i$ ranging over the rows of the indicated blocks of
Table~\ref{c:input-cells},
\begin{equation}\label{c:input-norms}
\begin{aligned}
\sum_{i,\,l\in L_i}B_{i,l}&\leq19.2,&\qquad
\sum_{i,\,l\in L_i}B_{i,l}^2&\leq19.86&&\qquad\text{(first block)},\\
\sum_{i,\,l\in L_i}B_{i,l}&\leq34,&\qquad
\sum_{i,\,l\in L_i}B_{i,l}^2&\leq61&&\qquad\text{(last two blocks)}.
\end{aligned}
\end{equation}
\end{enumerate}
\end{proposition}

\section{Rigorous evaluation of the approximate functions}\label{sec:CAP}

Subsection~\ref{sec:main-residual} and Section~\ref{sec:coercivity} reduce the analytic estimates to weighted integrals and pointwise inequalities for explicit fields. We enclose these quantities by interval arithmetic. The code and interval data are given in \cite{proof_code}. The implementation uses Julia~\cite{Julia}, \texttt{IntervalArithmetic.jl}~\cite{IntervalArithmeticJl}, and \texttt{RadiiPolynomial.jl}~\cite{RadiiPolynomialJl}.

As in \eqref{sec3eq5.2}, we consider functions $f$ of the form
\begin{equation}\label{eq : decomp general f}
    f = \chi f_{\rm in} + (1-\chi) f_{\rm ext},
\end{equation}
where $f_{in}$ is a tensor-product spline on the computational square $S_L$,  $f_{ext}$ is an asymptotic expansion in the exterior region and  $\chi$ is a cutoff function. The interior and exterior functions are defined in Subsections~\ref{sec:computational-domain} and~\ref{fixed_asymp}, respectively, and the cutoff is given by \eqref{res:cutoff}. 

\subsection{Local polynomial representations}\label{ssec:approximate solution}

On the computational domain $S_L$, the interior fields are
tensor-product splines. Their basis functions are piecewise
polynomials with compact support, so a restriction to one cell
involves only a fixed number of spline coefficients. We recall
their definition and the properties used below; see \cite{deBoor,PT}
for further details.

For a clamped knot vector $t_1\leq\cdots\leq t_m$ and a degree $p$,
the B-splines are defined recursively by
\begin{equation}\label{eq:coxdeboor}
\begin{aligned}
B_{n,0}(u)&=\mathbf 1_{[t_n,t_{n+1})}(u),\\
B_{n,q}(u)&=\frac{u-t_n}{t_{n+q}-t_n}B_{n,q-1}(u)
+\frac{t_{n+q+1}-u}{t_{n+q+1}-t_{n+1}}B_{n+1,q-1}(u),
\qquad 1\leq q\leq p.
\end{aligned}
\end{equation}
A term with zero denominator is taken to be zero. On a nonempty knot
interval $[t_j,t_{j+1}]$, only
$B_{j-p,p},\ldots,B_{j,p}$ can be nonzero. The recursion on this interval
is polynomial arithmetic, starting with $B_{j,0}=1$ on its interior and extending to the closed cell by limits.

Let $B^u_{\alpha,p}$ and
$B^v_{\beta,p}$ be the bases associated with the knot vectors in
$u$ and $v$, respectively. A tensor-product spline is specified by
coefficients $c_{\alpha\beta}$ through
\begin{equation}\label{cap:tensor-spline}
f(u,v)=\sum_{\alpha,\beta}c_{\alpha\beta}
 B^u_{\alpha,p}(u)B^v_{\beta,p}(v).
\end{equation}
The sums range over the basis indices in each variable. On a fixed
knot cell, only $(p+1)^2$ terms contribute. To evaluate products
and derivatives from these data, we convert this restriction to
the monomial basis in the coordinates
\begin{equation}\label{cap:local-coordinates}
u=\frac{u_i+u_{i+1}}2+\frac{u_{i+1}-u_i}2x,\qquad
v=\frac{v_j+v_{j+1}}2+\frac{v_{j+1}-v_j}2y,
\qquad (x,y)\in[-1,1]^2.
\end{equation}

\begin{lemma}\label{lem:tensorloc}\label{cap:derivative-models}
Let $f$ be the spline in \eqref{cap:tensor-spline}. On a knot cell
$[u_i,u_{i+1}]\times[v_j,v_{j+1}]$, its restriction is a polynomial
$\Pi(x,y)$ of bidegree at most $(p,p)$ in the coordinates
\eqref{cap:local-coordinates}. The derivatives of this polynomial
piece satisfy
\begin{equation}\label{eq:tensorloc}
\partial_u^a\partial_v^b f
=\left(\frac2{u_{i+1}-u_i}\right)^a
 \left(\frac2{v_{j+1}-v_j}\right)^b
 \partial_x^a\partial_y^b\Pi .
\end{equation}
\end{lemma}

\begin{proof}
Let $\mathcal A_u$ and $\mathcal A_v$ denote the active basis indices
on the two knot intervals. The recursion \eqref{eq:coxdeboor}
gives their polynomial restrictions
\[
B^u_{\alpha,p}(u(x))=\sum_{k=0}^p b^u_{\alpha k}x^k,
\qquad
B^v_{\beta,p}(v(y))=\sum_{l=0}^p b^v_{\beta l}y^l,
\qquad \alpha\in\mathcal A_u,\ \beta\in\mathcal A_v.
\]
Substituting these expressions in \eqref{cap:tensor-spline} and
collecting equal powers of $x,y$ gives
\begin{equation}\label{cap:spline-coefficient-conversion}
\Pi(x,y)=\sum_{k,l=0}^p C_{kl}x^ky^l,
\qquad
C_{kl}=\sum_{\alpha\in\mathcal A_u}\sum_{\beta\in\mathcal A_v}
 b^u_{\alpha k}c_{\alpha\beta}b^v_{\beta l}.
\end{equation}
Thus the coefficients of $\Pi$ are obtained from the given spline
data by a change of basis in each variable. In matrix form, let
$T_u=(b^u_{\alpha k})$, $T_v=(b^v_{\beta l})$, and let
$c_{\rm act}=(c_{\alpha\beta})_{\alpha\in\mathcal A_u,\,
\beta\in\mathcal A_v}$ be the active coefficient block. Then
$C=T_u^Tc_{\rm act}T_v$. Finally, \eqref{eq:tensorloc} follows
from the affine change of variables.
\end{proof}

To apply the differential operators globally, we also need the
polynomial pieces to join with sufficient regularity. The next
result recalls this property and the support of the basis functions.

\begin{lemma}\label{cap:spline-regularity}
The support of $B_{n,p}$ is contained in $[t_n,t_{n+p+1}]$.
At an interior knot of multiplicity $m\leq p$, it is $C^{p-m}$.
Consequently, a tensor-product spline of degree $p$ with simple
interior knots is $C^{p-1}$ on its closed computational square,
with one-sided derivatives at its boundary.

For the degree-$13$ profile splines and the cutoff
\eqref{res:cutoff}, the functions in \eqref{sec3eq5.2} are $C^4$.
Their physical lifts have the parities stated in
Subsection~\ref{sec:computational-domain}.
\end{lemma}

\begin{proof}
The support assertion follows by induction in \eqref{eq:coxdeboor}.
On nonempty knot intervals, differentiation gives
\[
B_{n,p}'=\frac{p}{t_{n+p}-t_n}B_{n,p-1}
-\frac{p}{t_{n+p+1}-t_{n+1}}B_{n+1,p-1},
\]
with the same convention for zero denominators. Induction shows that
the one-sided derivatives through order $p-m$ agree at a knot of
multiplicity $m$. Taking tensor products proves the regularity claim.
The supplied profile knots are simple in the interior, so the
profiles are $C^{12}$ there. A clamped spline need not have a smooth
extension by zero at an outer endpoint. Here this extension is made
only after multiplication by $\chi$: its derivatives through order
four vanish at the outer boundary. The term $(1-\chi)f_{\rm ext}$ is defined to be zero near the
origin. Thus $\chi f_{\rm in}+(1-\chi)f_{\rm ext}$ is $C^4$,
without requiring matching boundary data. Finally, the coordinate maps in \eqref{u:cov} are smooth functions
of $r^2,z^2$ at the axes; the extra factor $z$ gives the odd physical
fields.
\end{proof}

Lemmas~\ref{lem:tensorloc} and~\ref{cap:spline-regularity} give
the regularity of the interior fields and a polynomial formula for
each derivative on a cell. It remains to include the exterior
functions in the band $S_{band}=S_{100}\setminus S_{84}$.
In the gluing region, we use the exact identity
\begin{equation}\label{X:glue}
f=f_{\rm ext}+\chi(f_{\rm in}-f_{\rm ext}).
\end{equation}
The residuals involve second derivatives and the coercivity
coefficients involve third derivatives. All these derivatives,
including the mixed ones, are given, for $a,b\in\mathbb N_0$
with $a+b\leq3$, by
\begin{equation}\label{cap:gluing-jet-rule}
\partial_u^a\partial_v^b f
=\partial_u^a\partial_v^b f_{\rm ext}
+\sum_{i=0}^a\sum_{j=0}^b
 \binom ai\binom bj
 (\partial_u^i\partial_v^j\chi)
 \partial_u^{a-i}\partial_v^{b-j}(f_{\rm in}-f_{\rm ext}).
\end{equation}

\subsection{Local expansions and their errors}\label{ssec : local taylor}
\label{sec:tm}

Each spline piece is a polynomial. In the gluing region the exterior fields also contain powers and logarithms. We represent both by a polynomial with interval coefficients and a uniform error bound, after identifying each cell with $[-1,1]^2$ by \eqref{cap:local-coordinates}.

\begin{definition}\label{def : polynomial model}
    A local model $(\Pi,\epsilon)$ consists of a polynomial $\Pi$ 
and a nonnegative number $\epsilon$. We say that $(\Pi,\epsilon)$ encloses a function $f$ on a cell $C$ if
    \begin{equation}\label{eq : enclosure}
        |f-\Pi|\leq\epsilon\quad\hbox{on }C.
    \end{equation}
\end{definition}

In particular, an exact polynomial gives a model with $\epsilon=0$.
In the computation, the coefficients of $\Pi$ are enclosed by
intervals, and the formulas below are evaluated with outward rounding.

To control the degree when models are multiplied, we project onto
polynomials of bidegree at most $(N,N)$. If
$\Pi=\sum c_{ij}x^iy^j$, write
\begin{equation}\label{cap:coefficient-norm}
|\Pi|_1=\sum_{i,j}|c_{ij}|,\qquad
\pi^{\leq N}\Pi=\sum_{0\leq i,j\leq N}c_{ij}x^iy^j.
\end{equation}

Sums and products of the polynomial parts are formed on their interval coefficients, and the derivatives of the spline pieces are computed on the local polynomials by Lemma~\ref{lem:tensorloc}. For the analytic and glued factors, we first form the derivative expressions by \eqref{cap:gluing-jet-rule} and~\eqref{X:uderiv}, then enclose each derivative separately: a bound for $f-\Pi$ alone gives no bound for its derivatives.

\begin{lemma}\label{l:tmops}
Suppose that $(\Pi_1,\epsilon_1)$ and $(\Pi_2,\epsilon_2)$ enclose
$f_1$ and $f_2$. Then $|f_i|\leq|\Pi_i|_1+\epsilon_i$, and
$(\Pi_1\pm\Pi_2,\epsilon_1+\epsilon_2)$ encloses $f_1\pm f_2$.
A model of $f_1f_2$ is
\begin{equation}\label{cap:tm-product}
\left(\pi^{\leq N}(\Pi_1\Pi_2),\,
\epsilon_1|\Pi_2|_1+\epsilon_2|\Pi_1|_1+\epsilon_1\epsilon_2
+|(I-\pi^{\leq N})(\Pi_1\Pi_2)|_1\right).
\end{equation}
Moreover, $(\pi^{\leq N}\Pi,\epsilon+|(I-\pi^{\leq N})\Pi|_1)$ encloses $f$
whenever $(\Pi,\epsilon)$ does.
\end{lemma}

\begin{proof}
Every monomial has modulus at most one on $[-1,1]^2$. Write
$f_i=\Pi_i+(f_i-\Pi_i)$ and expand the product. Its difference from
$\pi^{\leq N}(\Pi_1\Pi_2)$ is the discarded polynomial plus
$(f_1-\Pi_1)\Pi_2+(f_2-\Pi_2)\Pi_1+(f_1-\Pi_1)(f_2-\Pi_2)$.
The triangle inequality gives the assertion.
\end{proof}

Besides the integrals considered below, the proof requires pointwise
coefficient inequalities. For these, a model gives the range enclosure
\begin{equation}\label{cap:model-range}
f([-1,1]^2)\subseteq
c_{00}+[-1,1]\left(\sum_{i+j>0}|c_{ij}|+\epsilon\right).
\end{equation}
If this range is too large, we subdivide the cell. Conversely, if $f(C)\subseteq[a,b]$, then $(c,\max\{|a-c|,|b-c|\})$ encloses $f$ for every $c\in\mathbb R$; we use such models when an analytic expansion is inconvenient, for example for $r=\sqrt{U_0(u)}$ on a cell touching the axis, and on cells that meet several pieces, with the interval hull of their ranges.

It remains to construct models for the elementary functions in the
exterior expressions. We expand powers and logarithms about a positive
constant. Given a model $(\Pi,\epsilon)$ of $g$, choose $c>0$ as the
midpoint of its constant-coefficient enclosure and put
\begin{equation}\label{cap:analytic-model-input}
\eta=\frac{g-c}{c},\qquad
q=\frac{|\Pi-c|_1+\epsilon}{c}.
\end{equation}
Lemma~\ref{l:tmops} gives $\|\eta\|_\infty\leq q$. The following
two results give explicit bounds for the omitted terms, provided
that this computable $q$ satisfies the stated convergence conditions.

\begin{lemma}\label{l:binom}
Let $c>0$ and let $g\colon C\to\R$ on a cell $C$. Suppose that
$\|(g-c)/c\|_{L^\infty(C)}\leq q<1$. For $\alpha\in\R$ and $N\geq1$,
if
\[
q\max\left\{1,\frac{|\alpha-N|}{N+1}\right\}<1,
\]
then
\begin{equation}\label{cap:tm-power}
\left\|g^\alpha-c^\alpha\sum_{k=0}^N
 \binom{\alpha}{k}\left(\frac{g-c}{c}\right)^k\right\|_{L^\infty(C)}
\leq
\frac{c^\alpha|\binom{\alpha}{N+1}|q^{N+1}}
 {1-q\max\{1,|\alpha-N|/(N+1)\}}.
\end{equation}
\end{lemma}

\begin{proof}
For $k\geq N$ the ratio $|\alpha-k|/(k+1)$ is at most
$\max\{1,|\alpha-N|/(N+1)\}$. The absolute values of the omitted
binomial terms are therefore bounded by a geometric series.
\end{proof}

\begin{lemma}\label{l:log}
With the notation of Lemma~\ref{l:binom}, suppose $\|(g-c)/c\|_{L^\infty(C)}\leq q<1$. Then
\begin{equation}\label{C:logeps}
\left\|\log g-\log c-\sum_{k=1}^N
 \frac{(-1)^{k+1}}{k}\left(\frac{g-c}{c}\right)^k\right\|_{L^\infty(C)}
\leq\frac{q^{N+1}}{(N+1)(1-q)}.
\end{equation}
\end{lemma}

\begin{proof}
Use the power series of $\log(1+(g-c)/c)$ and
$\sum_{k>N}q^k/k\leq q^{N+1}/((N+1)(1-q))$.
\end{proof}

The finite sums are evaluated using Lemma~\ref{l:tmops}, and the
displayed tails are added to their model errors. This encloses, in
particular, $(1+u)^{-m}$, $(1+v)^{-m}$, $(1+R^2)^{-m}$,
$R^{-1}$ and $\log R$ on every cell where they occur. Positivity
and convergence are checked before applying an expansion.

In the interior, the residuals \eqref{Emod} and \eqref{Emod2} are computed from their numerators $(1+u)^7(1+v)^7E_{\alpha,in}$, $\alpha\in\{\Psi,\Omega,\wt\Psi,\wt\Omega\}$, which are polynomials on each knot cell. They are formed without truncation, so that cancellations between their terms are retained, before the factor $(1+u)^{-7}(1+v)^{-7}$ is restored by Lemma~\ref{l:tmops}. In the band the same operations use the gluing formula \eqref{cap:gluing-jet-rule}, and the errors of the analytic factors are propagated by Lemma~\ref{l:tmops}.

\subsection{Weighted integrals on a cell}\label{sec:2norm}

We enclose integrals of products of two local models. The polynomial parts are integrated term by term, with their signs, and only the remainders are bounded in absolute value.

We first recall the integral of a polynomial on the normalized cell.
For $\Pi=\sum c_{ij}x^iy^j$,
\begin{equation}\label{cap:tm-integral}
\int_{[-1,1]^2}\Pi\,dx\,dy
=4\sum_{\substack{i,j\ {\rm even}}}\frac{c_{ij}}{(i+1)(j+1)}.
\end{equation}
Thus a model $(\Pi,\epsilon)$ gives the integral enclosure obtained
by adding $[-4\epsilon,4\epsilon]$ to this finite sum.

The weights arising from the physical measure separate into a
factor in each variable. Let $\rho(x),\mu(y)$ be nonnegative densities
with models $(\Pi_\rho,\epsilon_\rho)$ and
$(\Pi_\mu,\epsilon_\mu)$. Lemma~\ref{l:tmops} gives the density error
\begin{equation}\label{cap:density-error}
\|\rho\mu-\Pi_\rho\Pi_\mu\|_{L^\infty}
\leq\epsilon_{\rho\mu}:=
\epsilon_\rho|\Pi_\mu|_1+\epsilon_\mu|\Pi_\rho|_1
+\epsilon_\rho\epsilon_\mu.
\end{equation}

\begin{lemma}\label{cap:integral-lemma}
Suppose that $(\Pi_f,\epsilon_f)$ and $(\Pi_g,\epsilon_g)$ enclose
$f$ and $g$, and that the densities are as above. If
$\Pi_f=\sum c_{ij}x^iy^j$ and $\Pi_g=\sum d_{kl}x^ky^l$, then
\begin{equation}\label{cap:weighted-pairing}
\begin{aligned}
&\left|\int fg\rho\mu
-\sum_{i,j,k,l}c_{ij}d_{kl}
 \left(\int_{-1}^1x^{i+k}\Pi_\rho\,dx\right)
 \left(\int_{-1}^1y^{j+l}\Pi_\mu\,dy\right)\right|\\
&\quad\leq
4\epsilon_{\rho\mu}|\Pi_f|_1|\Pi_g|_1
+\big(\epsilon_f|\Pi_g|_1+\epsilon_g|\Pi_f|_1
+\epsilon_f\epsilon_g\big)
\left(\int_{-1}^1\Pi_\rho\,dx\int_{-1}^1\Pi_\mu\,dy
+4\epsilon_{\rho\mu}\right).
\end{aligned}
\end{equation}
All double integrals in this formula are over $[-1,1]^2$.
When $g=f$, the first term on the right can be replaced by
$\epsilon_{\rho\mu}\int_{[-1,1]^2}\Pi_f^2$.
\end{lemma}

\begin{proof}
The difference $fg-\Pi_f\Pi_g$ is bounded by
$\epsilon_f|\Pi_g|_1+\epsilon_g|\Pi_f|_1+\epsilon_f\epsilon_g$.
Integrate this error against the nonnegative density $\rho\mu$, using
\[
\int\rho\mu\leq
\int_{-1}^1\Pi_\rho\,dx\int_{-1}^1\Pi_\mu\,dy
+4\epsilon_{\rho\mu}.
\]
The remaining difference is
$\int\Pi_f\Pi_g(\rho\mu-\Pi_\rho\Pi_\mu)$.
Its modulus is at most
$4\epsilon_{\rho\mu}|\Pi_f|_1|\Pi_g|_1$; for a square it is at most
$\epsilon_{\rho\mu}\int\Pi_f^2$. Finally, expand the polynomial
pairing and separate the variables.
\end{proof}

The polynomial sum in \eqref{cap:weighted-pairing} is a finite combination of one-dimensional moments of $\Pi_\rho$ and $\Pi_\mu$, which are explicit sums of their coefficients. For $g=f$ we intersect the resulting enclosure with $[0,\infty)$.

On cells meeting the axis $z=0$, the density \eqref{res:measure} for even fields contains the unbounded factor $U_0(v)^{-1/2}$. There we integrate in the variable $s=\sqrt v$: the density becomes analytic in $(u,s)$, and the spline pieces remain polynomials in the corresponding normalized coordinate, with unchanged model errors. The derivatives are formed before this change of variables.

\subsection{The finite exterior expressions}\label{ssec : asymptotics arithmetic}

Outside the computational square, the approximate functions are
the finite expressions of Subsection~\ref{fixed_asymp}. After dividing
an odd field by $z$, every term is a product of powers of $R^{-1}$
and $\log R$ with a polynomial in $r^2/R^2$. We use this structure
to develop an arithmetic on the coefficients of the entire exterior
expression, before evaluating it on a cell or estimating its size.

For a multi-index $n=(n_1,n_2,n_3)\in\mathbb N_0^3$, write
\begin{equation}\label{cap:asymptotic-polynomial}
x=\frac{r^2}{R^2},\qquad t=R^{-1},\qquad \ell=\log R,\qquad
F(x,t,\ell)=\sum_{n\in\mathbb N_0^3}F_nx^{n_1}t^{n_2}\ell^{n_3}.
\end{equation}
Only finitely many coefficients $F_n$ are nonzero. The index $n_1$
is the angular degree, $n_2$ records the decay in $R$, and $n_3$ is
the logarithmic degree. Here $F$ represents the quotient by $z$ for
an odd field and the physical function itself for an even field.
We identify the expression with its finite array $(F_n)$, extended
by zero to all other indices, including those with a negative
component.

We first describe the algebra of these arrays. Let $F,G$ be two
such expressions, $c\in\mathbb R$, and write $m\leq n$ for
componentwise comparison. Addition and scalar multiplication act
coefficientwise; multiplication is the Cauchy product
\begin{equation}\label{cap:exterior-coefficient-arithmetic}
(F+G)_n=F_n+G_n,\qquad (cF)_n=cF_n,\qquad
(FG)_n=(F*G)_n:=\sum_{0\leq m\leq n}F_mG_{n-m}.
\end{equation}
In particular, the coefficients of a square are $F*F$. Multiplication by $x$, $t$, $\ell$ acts on the arrays by shifts of the indices, and differentiation in these variables by shifts combined with multiplication by the corresponding exponents. 
All products retain their full finite support.

We next relate this algebra to physical derivatives. Differentiating
$x,t,\ell$ with respect to $r,z$ gives
\[
\frac1r\partial_rF=t^2\mathcal A_rF,\qquad
\frac1z\partial_zF=t^2\mathcal A_zF,
\]
where
\begin{equation}\label{cap:asymptotic-derivatives}
\mathcal A_r=2(1-x)\partial_x-t\partial_t+\partial_\ell,\qquad
\mathcal A_z=-2x\partial_x-t\partial_t+\partial_\ell.
\end{equation}
The Laplacians and higher derivatives needed in the residuals and
the gluing formulas can now be formed by composing these finite
operations. We record the resulting identities.

\begin{proposition}\label{p:ops}
 Recall $L_0$ from \eqref{sf1.4}. For $\epsilon\in\{0,1\}$ define
\begin{equation}\label{cap:asymptotic-laplacian}
\mathcal L_\epsilon F
=t^2\big((4-2x)\mathcal A_rF+(2x-1+2\epsilon)\mathcal A_zF
+x\mathcal A_r^2F+(1-x)\mathcal A_z^2F\big).
\end{equation}
Then $L_0F=\mathcal L_0F$ and
$L_0(zF)/z=\mathcal L_1F$. For any $m\in\mathbb N_0$ and any finite
expression $G$,
\begin{equation}\label{cap:scaled-derivative-recursion}
\frac1r\partial_r(t^{2m}G)
=t^{2m+2}(\mathcal A_r-2mI)G,\qquad
\frac1z\partial_z(t^{2m}G)
=t^{2m+2}(\mathcal A_z-2mI)G.
\end{equation}
In particular,
\[
\left(\frac1r\partial_r\right)^2F
=t^4(\mathcal A_r^2F-2\mathcal A_rF),\qquad
\left(\frac1z\partial_z\right)^2F
=t^4(\mathcal A_z^2F-2\mathcal A_zF).
\]
\end{proposition}

\begin{proof}
Differentiate $x,t,\ell$ with respect to $r$ and $z$ to obtain
\eqref{cap:asymptotic-derivatives}. Differentiating $t^{2m}$ gives the correction $-2m$ in
\eqref{cap:scaled-derivative-recursion}. Taking $m=1$ gives the
second scaled derivative formulas.
Finally use
\[
L_0F=r^2\left(\frac1r\partial_r\right)^2F
+4\frac1r\partial_rF
+z^2\left(\frac1z\partial_z\right)^2F+\frac1z\partial_zF,
\qquad L_0(zF)/z=L_0F+2\frac{\partial_zF}{z}.
\]
\end{proof}

To combine these exterior expressions with the spline pieces in
the band, we need their $u,v$ derivatives. From \eqref{u:U0},
$\partial_u=((1+u)^3/6)(1/r)\partial_r$; repeated differentiation
gives
\begin{equation}\label{X:uderiv}
\begin{aligned}
\partial_uF&=\frac{(1+u)^3}{6}\frac1r\partial_rF,\qquad
\partial_u^2F=\frac{(1+u)^2}{2}\frac1r\partial_rF
+\frac{(1+u)^6}{36}\left(\frac1r\partial_r\right)^2F,\\
\partial_u^3F&=(1+u)\frac1r\partial_rF
+\frac{(1+u)^5}{4}\left(\frac1r\partial_r\right)^2F
+\frac{(1+u)^9}{216}\left(\frac1r\partial_r\right)^3F .
\end{aligned}
\end{equation}
The $v$ derivatives are given by the same formulas with $u,r$
replaced by $v,z$. Formula~\eqref{cap:scaled-derivative-recursion}
gives the third scaled derivatives and all mixed derivatives. We first form the finite coefficient expressions for these
derivatives, then evaluate them with the local models of $x,t,\ell$
from Subsection~\ref{sec:tm}.

On a spherical exterior the angular and radial variables separate. This gives the exterior integrals, with the weights $r^az^b$ of the residual estimates and the full $\mathbb R^5$ measure, including reflection in $z$.

\begin{proposition}\label{p:gram}
Let $a,b$ be nonnegative integers and $L_e>0$.
Suppose that $n_2>5+a+b$ whenever $(F*F)_n\neq0$, where $F*F$
is the Cauchy product in \eqref{cap:exterior-coefficient-arithmetic}. Then
\begin{equation}\label{X:quadform}
4\pi^2\iint_{\substack{r,z>0\\R\geq L_e}}
r^{3+a}z^bF^2\,dr\,dz
=4\pi^2\sum_n(F*F)_n\,
\mathcal A_{2n_1+3+a,b}\mathcal B_{4+a+b-n_2,n_3}(L_e),
\end{equation}
where
\begin{equation}\label{cap:exterior-moments}
\begin{aligned}
\mathcal A_{p,q}&=\int_0^{\pi/2}\sin^p\theta\cos^q\theta\,d\theta,\\
\mathcal B_{q,j}(L_e)&=\int_{L_e}^\infty R^q(\log R)^j\,dR
=\frac{L_e^{q+1}}{-q-1}
\sum_{i=0}^j\frac{j!}{(j-i)!}
 \frac{(\log L_e)^{j-i}}{(-q-1)^i},\qquad q<-1.
\end{aligned}
\end{equation}
\end{proposition}

\begin{proof}
Set $r=R\sin\theta$, $z=R\cos\theta$, so that
$dr\,dz=R\,dR\,d\theta$. The term indexed by $n$ in $r^{3+a}z^bF^2$ gives the angular
integral $\mathcal A_{2n_1+3+a,b}$ and the radial integral
$\mathcal B_{4+a+b-n_2,n_3}$. Repeated integration by parts gives the
displayed expression for $\mathcal B$.
\end{proof}

The angular moments $\mathcal A_{p,q}$ are evaluated by the standard reduction formulas in $p$ and $q$. The exterior of the computational square lies in a spherical exterior, so Proposition~\ref{p:gram} bounds its integrals; the same representation also bounds exterior suprema.

\begin{lemma}\label{cap:exterior-bound}
For $L>0$ such that $L_e:=\sqrt{U_0(L)}>1$, the coordinate map sends $S_L=[0,L]^2$
onto the physical square $[0,L_e]^2$. Its exterior in the
nonnegative physical quadrant is
\[
\mathcal E_L:=\{(r,z)\in[0,\infty)^2:\max\{r,z\}>L_e\}.
\]
Let $a,b\in\mathbb N_0$, and let $F(x,t,\ell)$ be a polynomial
with real coefficients as in \eqref{cap:asymptotic-polynomial},
with finitely many nonzero coefficients $F_n$. In the integrals
and suprema below, $F$ denotes its evaluation at
$(x,t,\ell)=(r^2/R^2,R^{-1},\log R)$. Define
\[
\mathfrak I_{a,b}(F;L_e)
:=4\pi^2\iint_{\substack{r,z>0\\R\geq L_e}}
 r^{3+a}z^bF^2\,dr\,dz.
\]
If $n_2>5+a+b$ whenever $(F*F)_n\neq0$, then
\begin{equation}\label{cap:exterior-domain-bound}
0\leq4\pi^2\iint_{\mathcal E_L}r^{3+a}z^bF^2\,dr\,dz
\leq\mathfrak I_{a,b}(F;L_e),
\end{equation}
where the last quantity is given by \eqref{X:quadform}.
Independently, suppose that $n_2\geq a+b$ and
$(n_2-a-b)\log L_e\geq n_3$ whenever $F_n\neq0$. Then
\begin{equation}\label{cap:exterior-supremum}
\sup_{\substack{r,z\geq0\\R\geq L_e}}|r^az^bF|
\leq\mathfrak S_{a,b}(F;L_e)
:=\sum_n|F_n|\,
 M_{2n_1+a,b}L_e^{a+b-n_2}(\log L_e)^{n_3},
\end{equation}
where
\[
M_{p,q}:=\max_{0\leq\theta\leq\pi/2}
 \sin^p\theta\cos^q\theta
=\left(\frac p{p+q}\right)^{p/2}
 \left(\frac q{p+q}\right)^{q/2}
\]
for $p,q>0$, and $M_{p,q}=1$ if either exponent is zero.
\end{lemma}

\begin{proof}
If $(r,z)\in\mathcal E_L$, then $r>L_e$ or $z>L_e$. Hence
$R>L_e$, and the integrand in
\eqref{cap:exterior-domain-bound} is nonnegative. For the pointwise
bound, write $r=R\sin\theta$, $z=R\cos\theta$ term by term.
The hypotheses imply
\[
\frac d{dR}\big(R^{a+b-n_2}(\log R)^{n_3}\big)\leq0
\qquad(R\geq L_e).
\]
Maximizing the angular factor gives $M_{2n_1+a,b}$, and summing the
absolute values proves \eqref{cap:exterior-supremum}.
\end{proof}

\section{Verification of the quantitative estimates}\label{sec:residual}\label{sec:numerical-residual}

We apply Section~\ref{sec:CAP} to prove the profile, coercivity, and response estimates used in Proposition~\ref{close:proposition}. 

Given $h,g$ with the same parity, and a  union $D$ of computational cells, we need to compute the following weighted integrals:
\begin{equation}\label{cap:shared-integrals}
\mathcal I_D(h,g):=\sum_{C\subseteq D}\iint_C hgJ\,du\,dv,
\qquad \mathcal I_D(h):=\mathcal I_D(h,h),
\end{equation}
where the density is given in
\eqref{res:measure} and repeated here for convenience:
\begin{equation}\label{X:measure}
J(u,v)=\frac{\pi^2}{9}
\frac{U_0(u)(1+u)^3(1+v)^3}{\sqrt{U_0(v)}}.
\end{equation}
It includes reflection in $z$, and the fields in the integrand
retain their physical factors $r,z$. Thus $\mathcal I_D(h)$ is the squared
$L^2(\R^5)$ norm on the corresponding reflected physical region.
We write $\mathfrak I_{a,b}(F)$ for
$\mathfrak I_{a,b}(F;L_e)$, where $L_e=\sqrt{U_0(100)}$.

\subsection{The profile bounds}\label{cap:profile-verification}

We first supply the residual and coefficient bounds of
Table~\ref{certi1}. By \eqref{res:final}, they imply the two combined
residual bounds in Table~\ref{close:profiles}. We then estimate
the corrected potentials and the test fields; these involve the
Poisson defects in addition to integrals of the prescribed profiles.

\subsubsection{Weighted residual integrals}\label{cap:residual-integrals}

We use the cleared numerators of Subsection~\ref{sec:tm} on
interior cells and \eqref{cap:gluing-jet-rule} in the band. Given a cell $C$, we apply
Lemma~\ref{res:bounds-lemma} by letting $w_C=\sup_Cw$, with $w$ as in
\eqref{res:w}. Monotonicity of $U_0$ and $k$ in \eqref{res:k}
places this supremum at the upper corner of $C$; $w(0,0)=0$.
Integration of $w\leq w_C$ on each cell gives
\begin{equation}\label{cap:stream-integrals}
\begin{aligned}
\|\sqrt w\,\mathcal R_\Omega\mathbf1_{S_{100}}\|_2^2
&\leq\sum_{C\subseteq S_{100}}w_C\mathcal I_C(\mathcal R_\Omega),\qquad \|\sqrt w\,\mathcal R_\Psi^{\rm in}\|_2^2
\leq\sum_{C\subseteq S_{84}}w_C\mathcal I_C(\mathcal R_\Psi),\\
\|\mathcal R_\Psi^{\rm band}\|_2^2
&=\mathcal I_{S_{band}}(\mathcal R_\Psi),\qquad \|\sqrt w\,\mathcal R_\Psi^{\rm band}\|_2^2
\leq\sum_{C\subseteq S_{band}}w_C\mathcal I_C(\mathcal R_\Psi).
\end{aligned}
\end{equation}

The band contributes both a weighted and an unweighted norm in
\eqref{res:far-bound}. For the even residuals the multiplier is $r$:
\begin{equation}\label{cap:even-integrals}
\|r\mathcal R_{\wOmega}\mathbf1_{S_{100}}\|_2^2
=\mathcal I_{S_{100}}(r\mathcal R_{\wOmega}),\qquad
\|r\mathcal R_{\wPsi}\|_2^2
=\mathcal I_{S_{100}}(r\mathcal R_{\wPsi}).
\end{equation}

The exterior stream defects vanish by
\eqref{res:exterior-identities}. For the other two exterior terms,
$w\leq r^2$ and the odd factor $z$ give
\begin{equation}\label{cap:exterior-norms}
\begin{aligned}
\|\sqrt w\,\mathcal R_\Omega^{\rm ext}\|_2^2
\leq\mathfrak I_{2,2}(\mathcal R_\Omega^{\rm ext}/z), \quad 
\|r\mathcal R_{\wOmega}^{\rm ext}\|_2^2
\leq\mathfrak I_{2,0}(\mathcal R_{\wOmega}^{\rm ext}).
\end{aligned}
\end{equation}

It remains to bound the coefficient factors in
Lemma~\ref{res:bounds-lemma}. Recall $\Phi$ and $I_o$ from
\eqref{res:Phi}, with $R_1=r_\ast$. For
$H=\Omega^\ast$ or $\wt\Omega^\ast$, monotonicity of $\Phi$ gives
\begin{equation}\label{cap:coefficient-integrals}
\int_{C\cap B_{900}}(r+2|z|)^2\Phi(R)^2H^2\,dy
\leq\Phi\big(\min\{900,\sup_CR\}\big)^2
\sup_C(r+2|z|)^2\mathcal I_C(H).
\end{equation}

Since $900<L_e$, the ball $B_{900}$ lies in the physical image of
$S_{100}$. Therefore
\begin{equation}\label{cap:coefficient-global-integral}
I_o(H,900)^2\leq\sum_{\substack{C\subseteq S_{100}\\C\cap B_{900}\ne\varnothing}}
\Phi\big(\min\{900,\sup_CR\}\big)^2\sup_C(r+2|z|)^2\mathcal I_C(H).
\end{equation}
A cell crossing the sphere is included in full.
For the suprema, a local model $(\Pi_H,\epsilon_H)$ gives
\begin{equation}\label{cap:cell-supremum}
\sup_C|rH|\leq\sup_Cr\,(|\Pi_H|_1+\epsilon_H).
\end{equation}

Write $H^{\rm ext}=\Omega_{\geq-11}$ for $H=\Omega^\ast$, and
$H^{\rm ext}=\wt\Omega_{\geq-9}$ for $H=\wt\Omega^\ast$.
For $A=\R^5$ or $\R^5\setminus B_{900}$, this gives
\begin{equation}\label{cap:coefficient-global-supremum}
\|rH\|_{L^\infty(A)}\leq\max\left\{
\max_{\substack{C\subseteq S_{100}\\C\cap A\ne\varnothing}}\sup_Cr\,(|\Pi_H|_1+\epsilon_H),
\|rH^{\rm ext}\|_{L^\infty(A\cap S_{\rm ext})}\right\}.
\end{equation}
The exterior contributions are bounded by
\[
\sup_{S_{\rm ext}}|r\Omega_{\geq-11}|
\leq\mathfrak S_{1,1}(\Omega_{\geq-11}/z;L_e),\qquad
\sup_{S_{\rm ext}}|r\wt\Omega_{\geq-9}|
\leq\mathfrak S_{1,0}(\wt\Omega_{\geq-9};L_e),
\]
using \eqref{cap:exterior-supremum}. These integrals and maxima give the residual and coefficient bounds in Table~\ref{certi1}. Equation~\eqref{res:final} then yields $\|\mathcal R\|_{Z'}<5.85\times10^{-12}$ and $\|\wt{\mathcal R}\|_{Z'}<1.5\times10^{-12}$. 

\subsubsection{Poisson corrections and profile norms}\label{cap:corrected-profile-norms}

The Poisson corrections in \eqref{cl3.3} satisfy
\begin{equation}\label{cap:profile-potential-corrections}
\Psi_5-\Psi^\ast=\Delta^{-1}\mathcal R_\Psi,\qquad
\wt\Psi_5-\wPsi^\ast=\Delta^{-1}\mathcal R_{\wPsi}.
\end{equation}
In particular, \eqref{robu1} and
the odd Hardy inequality give
\begin{equation}\label{cap:profile-correction}
\|\mathcal R_\Psi\|_Z^2
\leq\mathcal I_{S_{100}}(\mathcal R_\Psi)
+\|\sqrt w\,\mathcal R_\Psi^{\rm in}\|_2^2
+\|\sqrt w\,\mathcal R_\Psi^{\rm band}\|_2^2.
\end{equation}

The right-hand side gives $\|\mathcal R_\Psi\|_Z\leq1.251\times10^{-9}$.
For the even profile we also need a positive lower bound on the
corrected gradient, since it is the denominator of the projection
\eqref{robu21}.

\begin{lemma}\label{cap:profile-norm-lemma}
Define $G_{100}\geq0$ by
\begin{equation}\label{cap:normalization-integral}
G_{100}^2
:=\mathcal I_{S_{100}}(\partial_r\wt\Psi_5)
 +\mathcal I_{S_{100}}(\partial_z\wt\Psi_5)
\end{equation}
and, using the exterior derivative operators
\eqref{cap:asymptotic-derivatives}, put
\begin{equation}\label{cap:profile-exterior-energies}
G_{\rm ext}^2
:=\mathfrak I_{2,0}(t^2\mathcal A_r\wt\Psi_{\geq-7})
 +\mathfrak I_{0,2}(t^2\mathcal A_z\wt\Psi_{\geq-7}).
\end{equation}
Then
\begin{equation}\label{cap:corrected-even-gradient}
\begin{aligned}
\mathcal I_{S_{100}}(\wt\Omega^\ast)
&\leq\|\wt\Omega^\ast\|_2^2
\leq\mathcal I_{S_{100}}(\wt\Omega^\ast)
 +\mathfrak I_{0,0}(\wt\Omega_{\geq-9}),\\
\max\{0,G_{100}-\|r\mathcal R_{\wPsi}\|_2\}
&\leq\|\nabla\wPsi^\ast\|_2
\leq\sqrt{G_{100}^2+G_{\rm ext}^2}
 +\|r\mathcal R_{\wPsi}\|_2.
\end{aligned}
\end{equation}

Moreover, with $S=r\partial_r+z\partial_z$ as in \eqref{robu1},
\begin{equation}\label{cap:profile-scaling}
\begin{aligned}
\|\wt\Omega^\ast\|_Z^2
&=\|\wt\Omega^\ast\|_2^2+\|\nabla\wPsi^\ast\|_2^2,\qquad \|S\wt\Omega^\ast\|_{Z'}
\leq\|R\wt\Omega^\ast\|_2+5\|\nabla\wPsi^\ast\|_2,\\
\|R\wt\Omega^\ast\|_2^2
&\leq\mathcal I_{S_{100}}(r\wt\Omega^\ast)
 +\mathcal I_{S_{100}}(z\wt\Omega^\ast)
+\mathfrak I_{2,0}(\wt\Omega_{\geq-9})
 +\mathfrak I_{0,2}(\wt\Omega_{\geq-9}).
\end{aligned}
\end{equation}

\end{lemma}

\begin{proof}
The first bound in \eqref{cap:corrected-even-gradient} follows by
splitting the squared norm into $S_{100}$ and its exterior.
For the gradient, \eqref{res:radial-source} gives
\[
\|\nabla(\wPsi^\ast-\wt\Psi_5)\|_2
=\|\mathcal R_{\wPsi}\|_{\dot H^{-1}}
\leq\|r\mathcal R_{\wPsi}\|_2.
\]
Moreover,
\begin{equation}\label{cap:profile-energy-sandwich}
G_{100}^2\leq\|\nabla\wt\Psi_5\|_2^2\leq G_{100}^2+G_{\rm ext}^2.
\end{equation}
The triangle inequality gives \eqref{cap:corrected-even-gradient}.
Finally,
$Sh=\operatorname{div}(yh)-5h$, \eqref{res:radial-source} and
$R^2=r^2+z^2$ prove \eqref{cap:profile-scaling}.
\end{proof}

Thus the even-profile comparisons are
$\|\wt\Omega^\ast\|_2\in[1,1.01]$,
$\|\nabla\wPsi^\ast\|_2\in[0.88,0.89]$, and
$\|R\wt\Omega^\ast\|_2\leq22.55$.

\subsubsection{The drift}\label{cap:drift-test}

We estimate the $L^\infty(\R^5)$ norm of the drift
$(-r\partial_z\Psi_5,r\partial_r\Psi_5+2\Psi_5)$ required in
Table~\ref{close:profiles}. This norm multiplies the Poisson-defect
bound $\|r\mathcal R_{\wPsi}\|_2$ in
\eqref{cap:test-field-corrections}. On each cell
$C\subseteq S_{100}$, including the gluing band, let
$(\Pi_C,\epsilon_C)$ be a local model of the complete square
$r^2(\partial_z\Psi_5)^2+(r\partial_r\Psi_5+2\Psi_5)^2$.
Its computational expression is \eqref{cap:test-field-drift-model};
we form this square before taking the local range. The uniform
model error and the coefficient norm give
\begin{equation}\label{cap:global-drift-bound}
\begin{aligned}
\|(-r\partial_z\Psi_5,r\partial_r\Psi_5+2\Psi_5)\|_\infty^2
\leq\max\Big\{&\max_{C\subseteq S_{100}}(|\Pi_C|_1+\epsilon_C),\\
&\|(-r\partial_z\Psi_{\geq-9},
 r\partial_r\Psi_{\geq-9}+2\Psi_{\geq-9})\|_{L^\infty(S_{\rm ext})}^2\Big\}.
\end{aligned}
\end{equation}
The cell enclosures in \eqref{cap:global-drift-bound} and the exterior estimate \eqref{cap:exterior-drift-supremum} give
\begin{equation}\label{cap:drift-comparison}
\sup_{\R^5}\big[r^2(\partial_z\Psi_5)^2
 +(r\partial_r\Psi_5+2\Psi_5)^2\big]\leq22.45^2.
\end{equation}

\subsubsection{The corrected test fields}\label{cap:test-verification}

To estimate the joint norm of $a_\ast,b_\ast$ in
\eqref{close:test-fields}, we first replace the corrected
potentials by the prescribed (explicit) ones and define
\begin{equation}\label{cap:test-fields}
\begin{aligned}
a_5={}&\wt\Omega^\ast
+\partial_z\Psi_5(r\partial_r\wt\Psi_5+2\wt\Psi_5)
-(r\partial_r\Psi_5+2\Psi_5)\partial_z\wt\Psi_5,\\
b_5={}&\left(
\frac r2\wt\Omega^\ast-r\Omega^\ast\partial_z\wt\Psi_5,\,
\frac z2\wt\Omega^\ast+
\Omega^\ast(r\partial_r\wt\Psi_5+2\wt\Psi_5)
\right).
\end{aligned}
\end{equation}

\begin{lemma}\label{cap:test-correction-lemma}
With $C_b$ as in Lemma~\ref{nl:bilinear},
\begin{equation}\label{cap:test-field-corrections}
\begin{aligned}
\|a_\ast-a_5\|_2\leq{}&
C_b\|\mathcal R_\Psi\|_Z\|\wt\Omega^\ast\|_Z
+\|(-r\partial_z\Psi_5,r\partial_r\Psi_5+2\Psi_5)\|_\infty
                  \|r\mathcal R_{\wPsi}\|_2,\\
\|b_\ast-b_5\|_2\leq{}&\|r\Omega^\ast\|_\infty
                  \|r\mathcal R_{\wPsi}\|_2.
\end{aligned}
\end{equation}
Consequently,
\begin{equation}\label{cap:test-field-assembly}
\sqrt{\|a_\ast\|_2^2+\|b_\ast\|_2^2}
\leq\sqrt{\int_{\R^5}(a_5^2+|b_5|^2)}
+\sqrt{\|a_\ast-a_5\|_2^2+\|b_\ast-b_5\|_2^2}.
\end{equation}

\end{lemma}

\begin{proof}
The bilinear velocity estimate used in \eqref{nl4} gives
\begin{equation}\label{cap:velocity-product}
\left\|p_z(rq_r+2q)-(rp_r+2p)q_z\right\|_2
\leq C_b\|\Delta p\|_Z\|\Delta q\|_Z.
\end{equation}

This follows from \eqref{nl8}, H\"older's inequality and the
velocity $H^1$ to $L^4$ estimate, with constant
$3^{3/4}C_4^2/(4\sqrt\pi)<C_b$.

Set $p=\Delta^{-1}\mathcal R_\Psi$ and
$q=\Delta^{-1}\mathcal R_{\wPsi}$. Substituting
$\Psi^\ast=\Psi_5-p$, $\wPsi^\ast=\wt\Psi_5-q$ gives
\begin{equation}\label{cap:test-field-differences}
\begin{aligned}
a_\ast-a_5&={}-p_z(r\wPsi_r^\ast+2\wPsi^\ast)
 +(rp_r+2p)\wPsi_z^\ast -\Psi_{5,z}(rq_r+2q)+(r\Psi_{5,r}+2\Psi_5)q_z,\\
b_\ast-b_5&={}r\Omega^\ast(q_z,-q_r-2q/r).
\end{aligned}
\end{equation}
Apply \eqref{cap:velocity-product} to the terms containing $p$. The remaining
terms are bounded using the drift supremum,
$\|(q_r+2q/r,q_z)\|_2=\|\nabla q\|_2$ from
\eqref{close:shifted-gradient}, and
$\|\nabla q\|_2\leq\|r\mathcal R_{\wPsi}\|_2$.
This proves \eqref{cap:test-field-corrections}; the triangle
inequality in $L^2\oplus L^2$ gives
\eqref{cap:test-field-assembly}.
\end{proof}

The integral of $a_5^2+|b_5|^2$ can be enclosed directly
from the prescribed fields. Their cleared numerators give
\begin{equation}\label{cap:test-field-recovery}
a_5=\frac{N_{a_5}}{(1+u)^7(1+v)^7},\qquad
b_{5,r}=\frac{rN_{b_{5,r}}}{(1+u)^7(1+v)^7},\qquad
b_{5,z}=\frac{zN_{b_{5,z}}}{(1+u)^7(1+v)^7}.
\end{equation}

In particular, the bounded-domain integral is the integral of
\begin{equation}\label{cap:test-field-joint-integrand}
(a_5^2+|b_5|^2)J
=\big(N_{a_5}^2+r^2N_{b_{5,r}}^2+z^2N_{b_{5,z}}^2\big)
 (1+u)^{-14}(1+v)^{-14}J.
\end{equation}

The full-space bound is
\begin{equation}\label{cap:test-field-global-integral}
\begin{aligned}
\int_{\R^5}(a_5^2+|b_5|^2)
\leq{}&\mathcal I_{S_{100}}(a_5)
 +\mathcal I_{S_{100}}(b_{5,r})+\mathcal I_{S_{100}}(b_{5,z})\\
&+\mathfrak I_{0,0}(a_5^{\rm ext})
 +\mathfrak I_{2,0}(b_{5,r}^{\rm ext}/r)
 +\mathfrak I_{0,2}(b_{5,z}^{\rm ext}/z).
\end{aligned}
\end{equation}
The finite exterior expressions are
\eqref{cap:exterior-test-field-polynomials}.
Together with \eqref{cap:test-field-corrections}, the preceding
profile and drift bounds give
\begin{equation}\label{cap:test-field-correction-size}
\sqrt{\|a_\ast-a_5\|_2^2+\|b_\ast-b_5\|_2^2}\leq1.36\times10^{-10}.
\end{equation}

Interval evaluation of \eqref{cap:test-field-global-integral} yields
\[
\sqrt{\int_{\R^5}(a_5^2+|b_5|^2)}+1.36\times10^{-10}\leq14.4.
\]
Together with \eqref{cap:test-field-assembly}, this proves the corrected test-field bound and completes the proof of Proposition~\ref{PropCentral2}.

\subsection{Weak coercivity and the selected functionals}\label{cap:coercivity-verification}

By Theorem~\ref{c:thm:main}, weak coercivity follows from the
estimates of Proposition~\ref{c:prop:CAP}, namely pointwise bounds, a
complementary spectral gap, finite matrix inequalities, and norms
of the selected functionals. These estimates concern the prescribed
potential $\Psi_5$, and the correction to $\Psi^\ast$ is controlled by
\eqref{c:profile-loss}.

\subsubsection{Coefficient bounds and confinement}\label{cap:coefficient-positivity}

The coefficients \eqref{mos2} are determined by $\Psi_5$;
in particular $e_0=\partial_z\Delta_{\R^5}\Psi_5$ involves its
third derivatives. The formulas
of Appendix~\ref{app:coercivity-derivatives}
therefore apply to the prescribed potential, independently of
$\Omega^\ast$.

In the odd form, the box bounds \eqref{c:boxes} control the
negative part of $e_0$ through \eqref{c:zeta}. For $r\geq200$ we use the exterior box of Proposition~\ref{c:prop:boxes}(a).
Let $C=C_r\times C_z$ be the intersection of a cell with $\{0<r<200,\ z\geq0\}$ and suppose it has nonempty interior. Then
\begin{equation}\label{cap:zeta-cell-bound}
\zeta(r,z)\leq\bar\zeta_C:=\frac4{\pi^2}
\max_{\substack{k:\,\operatorname{int}(C_r)\cap(r_{k-1},r_k)\ne\varnothing}}
\sum_{p=1}^2h_{p,k}L_{p,k}^2\mathbf1_{\{\inf C_z<L_{p,k}\}}
\qquad((r,z)\in\operatorname{int}C).
\end{equation}

Thus, if $M^o_{zz}=a_{22}-\zeta$ as in \eqref{c:Mdef},
\begin{equation}\label{cap:coercivity-matrix-lower-bound}
M^o\geq
\begin{pmatrix}a_{11}&a_{12}\\a_{12}&a_{22}-\bar\zeta_C\end{pmatrix}
\qquad\hbox{on the interior of }C.
\end{equation}
The box construction and its boundary convention are given in
Appendix~\ref{app:coercivity-boxes}. To verify that this lower matrix and the matrices $M_i-\eta_iI+W_iI$ are nonnegative, we use
\begin{equation}\label{cap:positive-matrix-test}
\begin{pmatrix}a&b\\b&c\end{pmatrix}\geq0
\quad\Longleftrightarrow\quad a\geq0,\ c\geq0,\ ac-b^2\geq0.
\end{equation}
For the scalar rows the corresponding condition is $V-\eta_s+W_i\geq0$.
We verify these inequalities by range enclosures on each cell, and for
the determinant we enclose the range of a local model of $ac-b^2$ itself.

The far-field condition \eqref{c:Hyp} starts at $R=200$.
For each bounded spline or gluing cell $C$ meeting this region,
the required comparisons are
\begin{equation}\label{cap:bounded-far-field-tests}
\sup_{C\cap\{R\geq200\}}R^2|e-e^{(1)}|\leq0.05\sigma
\qquad(e\in\{e_{rr},e_{zz},a_{12}\}),\qquad \sup_{C\cap\{R\geq200\}}R^4|e_0-e_0^{(1)}|\leq\sigma.
\end{equation}
On the infinite exterior the corresponding
bounds are the finite sums
\eqref{cap:coercivity-exterior-comparisons}.
Together with the confinement estimates for $M^o$ and $V$, these bounds prove parts~(i)--(iv) of Proposition~\ref{c:prop:CAP}.

\subsubsection{The complementary spectral gap}\label{cap:spectral-gaps}

For each row $i$, the lower bound
$M_i-\eta_iI\geq-W_iI$ permits the use of
Theorem~\ref{c:thm:refined} once every omitted admissible mode has
$\lambda_m+\mu_k\geq\Lambda_{*,i}$ and
$d_i\Lambda_{*,i}>W_i$. The latter inequality follows from
Table~\ref{c:input-cells}. The former reduces to finitely many
comparisons as follows.

\begin{lemma}\label{cap:gap-lemma}
Let $\lambda_m,\mu_k$ be the nonnegative increasing eigenvalues of
a row of Table~\ref{c:input-cells}, with validated indices
$0\leq m\leq M$, $0\leq k\leq K$. Suppose that
\[
\lambda_M\geq\Lambda_*,\qquad \mu_K\geq\Lambda_*,
\]
and that every omitted admissible pair in this rectangle satisfies
$\lambda_m+\mu_k\geq\Lambda_*$. Then the same inequality holds for
every omitted admissible pair.
\end{lemma}

\begin{proof}
Outside the rectangle at least one index is larger than its last
validated index. Monotonicity and nonnegativity give the assertion.
\end{proof}

We consider only modes admissible for component $a$, as defined in
Subsection~\ref{c:ss:CAP}. Recall that $n_i=|L_i|$ is the number of
selected sources in row $i$ of Table~\ref{c:input-cells}.
When $n_i>0$, the sets $\Sigma_{i,a}$ and the selected coefficients
in \eqref{c:Xi} are fixed, since they define the sources
\eqref{c:qfinal}.

For a fixed row $i$, write the validated enclosures as
\[
\lambda_m\in[\lambda_m^-,\lambda_m^+],\qquad
\mu_k\in[\mu_k^-,\mu_k^+],\qquad
\Lambda_{*,i}\in[\Lambda_{*,i}^-,\Lambda_{*,i}^+].
\]
The eigenvalue enclosures are given in Appendix~\ref{app:eigenfunctions},
and $[\Lambda_{*,i}^-,\Lambda_{*,i}^+]$ is the outward rounding of the
decimal value of $\Lambda_{*,i}$ in Table~\ref{c:input-cells}. To apply Lemma~\ref{cap:gap-lemma}, we verify
\begin{equation}\label{cap:gap-interval-tests}
\lambda_M^-\geq\Lambda_{*,i}^+,\qquad
\mu_K^-\geq\Lambda_{*,i}^+,\qquad
\lambda_m^-+\mu_k^-\geq\Lambda_{*,i}^+.
\end{equation}
The last inequality is required for every admissible
$(m,k)\notin\Sigma_{i,a}$ with $0\leq m\leq M$, $0\leq k\leq K$.
These comparisons imply \eqref{c:input-modes} for every omitted
admissible mode, including those outside the validated rectangle.

If $n_i=0$, no source is prescribed on this row, and we choose
\[
\Sigma_{i,a}=\{(m,k)\text{ admissible for }a:\lambda_m+\mu_k<\Lambda_{*,i}\}.
\]
For a validated pair, membership in $\Sigma_{i,a}$ is decided by comparing the interval enclosures of $\lambda_m+\mu_k$ and $\Lambda_{*,i}$.
The inequalities \eqref{cap:gap-interval-tests} hold on every row of Table~\ref{c:input-cells}. This proves part~(v) of Proposition~\ref{c:prop:CAP}.

\subsubsection{Signed Gram integrals and finite matrices}\label{cap:finite-matrices}

We next enclose the finite matrix in Theorem~\ref{c:thm:refined}.
For a mode $j=(m,k)$, recall $e_j=\varphi_m\psi_k$ and
$\Lambda_j=\lambda_m+\mu_k$ from \eqref{c:modes}.
The matrices $G_i$, $G_{2,i}$, $\Lambda_i$ and $N_i$ are those of Subsection~\ref{c:ss:CAP}, defined by \eqref{c:NGG} with $M=M_i-\eta_iI$; in $G_{2,i}$ the square $(M_i-\eta_iI)^2$ is the pointwise matrix square.

For $p\in\{1,2\}$, let $(\Pi_{ab}^{(p)},\epsilon_{ab}^{(p)})$
be a local model on $C=C_r\times C_z$, with coordinates $x(r),y(z)$,
satisfying
\begin{equation}\label{cap:gram-coefficient-models}
\sup_{(r,z)\in C}\left|
\big[(M_i-\eta_iI)^p\big]_{ab}
-\Pi_{ab}^{(p)}(x(r),y(z))\right|\leq\epsilon_{ab}^{(p)}.
\end{equation}
These models are formed by the arithmetic of Lemma~\ref{l:tmops}.
For retained modes $(m,k)$ in component $a$ and $(n,l)$ in
component $b$, recall $e_{mk}=\varphi_m\psi_k$ and
$\varrho_i=T_iU_i$. Separation of variables gives, for both $p=1,2$,
\begin{equation}\label{cap:gram-moments}
\begin{aligned}
&\int_C\varrho_i\Pi_{ab}^{(p)}(x(r),y(z))e_{mk}e_{nl}r^3\,dr\,dz\\
&\quad=\sum_{\alpha,\beta}(\Pi_{ab}^{(p)})_{\alpha\beta}
\left(\int_{C_r}T_i\varphi_m\varphi_n x(r)^\alpha r^3\,dr\right)
\cdot\left(\int_{C_z}U_i\psi_k\psi_l y(z)^\beta\,dz\right).
\end{aligned}
\end{equation}
Cauchy--Schwarz bounds the error by
\begin{equation}\label{cap:gram-error}
\begin{aligned}
&\left|\int_C\varrho_i
\left[\big[(M_i-\eta_iI)^p\big]_{ab}-\Pi_{ab}^{(p)}(x(r),y(z))\right]
e_{mk}e_{nl}r^3\,dr\,dz\right|\\
&\quad\leq\epsilon_{ab}^{(p)}
\left(\int_{C_r}T_i\varphi_m^2r^3\,dr
      \int_{C_r}T_i\varphi_n^2r^3\,dr\right)^{1/2}
\left(\int_{C_z}U_i\psi_k^2\,dz
      \int_{C_z}U_i\psi_l^2\,dz\right)^{1/2}.
\end{aligned}
\end{equation}
The modes are the exact normalized modes; the moment enclosures
use \eqref{cap:tm-integral}.

The selected columns $\xi_{i,l}$, $l\in L_i$, define $\Xi_{i,l}$ and the sources in \eqref{c:qfinal}. We complete these columns with approximate eigenvectors of the compression of $N_i$ to their orthogonal complement. Since the completed matrix $X_i$ is only approximately orthogonal, we use the following lemma.

\begin{lemma}\label{cap:matrix-lemma}
Denote by $\|\cdot\|_F$ the Frobenius norm and by $\|\cdot\|$ the spectral norm.
For a real symmetric matrix $N$ and a square real matrix $X$, put
\[
\tau_F=\|I-XX^T\|_F,\qquad D=X^TNX.
\]

If $\tau_F<1$, then
\begin{equation}\label{cap:matrix-norm}
\|I-XX^T\|\leq\tau_F,\qquad
\|N\|\leq\min\left\{
\max_j\sum_k|N_{jk}|,
\frac{\max_j\sum_k|D_{jk}|}{1-\tau_F}\right\}.
\end{equation}

\end{lemma}

\begin{proof}
The spectral norm is bounded by the Frobenius norm. Also
$XX^T\geq(1-\tau_F)I$, so $X$ is invertible and
$\|X^{-1}\|\leq(1-\tau_F)^{-1/2}$. Use
$N=X^{-T}DX^{-1}$ and the row-sum bounds for the two symmetric
matrices.
\end{proof}

For $N=N_i$, $X=X_i$ and $D=X_i^TN_iX_i$, the comparisons are
\begin{equation}\label{cap:matrix-comparisons}
\begin{gathered}
\tau_F\leq10^{-12},\qquad\min\left\{\max_j\sum_k|N_{jk}|,\frac{\max_j\sum_k|D_{jk}|}{1-\tau_F}\right\}\leq30,\\
D_{ll}-\sum_{k\ne l}|D_{lk}|\geq
\begin{cases}-\nu_{i,l},&l\in L_i,\\0,&l\notin L_i,\end{cases}
\qquad (N=N_i,\ X=X_i).
\end{gathered}
\end{equation}

Put $w_{i,l}=\nu_{i,l}$ for $l\in L_i$ and $w_{i,l}=0$ otherwise. Over the fifteen rows of Table~\ref{c:input-cells}, interval evaluation of \eqref{cap:matrix-comparisons} gives
\[
\min_{i,l}(g_{i,l}+w_{i,l})\geq 2.2640\times10^{-6},\qquad
\max_i\|N_i\|\leq15.72930,\qquad \max_i\|I-X_iX_i^T\|_F\leq9.73\times10^{-14}.
\]

On the even row $T_0U_1$ we have $\zeta=0$ and $M^v=M^o+cI$, where $c=\widetilde\lambda_{app}-0.224$. After identifying the retained modes, the even Gram matrices are therefore $G+cI$ and $G_2+2cG+c^2I$. In \eqref{c:NGG} we then use the complementary gap $\Lambda_*$ of the even row.

These bounds prove part~(vi) of Proposition~\ref{c:prop:CAP}.

\subsubsection{Norms of the selected functionals}\label{cap:scalar-functionals}

To complete Proposition~\ref{c:prop:CAP}, we need the prescribed
quantities $B_{i,l}$ in \eqref{c:input-norms}.
For vector rows their definition and orthonormality give
\begin{equation}\label{cap:vector-functional-bound}
\|q_{(i,l)}\|_{Z'}^2\leq B_{i,l}^2
=\nu_{i,l}\int_{Q_i}\varrho_i^2|\Xi_{i,l}|^2r^3drdz
\leq\nu_{i,l}\sum_{a=1}^2\sum_{j\in\Sigma_{i,a}}
                  |(\xi_{i,l})_{j,a}|^2.
\end{equation}
The parity and reflection factors are included in \eqref{c:vector-norm}.

For scalar rows, \eqref{c:scalar-norm} instead requires a radial
primitive. Its explicit form explains both the bounded integral
and its noncompact tail.

\begin{lemma}\label{cap:source-primitive-lemma}
Let $T$ be a radial plateau cutoff of \eqref{c:plateau},
supported on $[a,b]$ with $0<a<b$,
and let $(\lambda_m,\varphi_m)_{m\geq0}$ be its normalized Neumann
eigenpairs for the weight $r^3T(r)$, as in
Lemma~\ref{c:lem:1dspec}. Set
$H_T(r):=\int_0^r t^3T(t)\,dt$ and choose
$\varphi_0=H_T(b)^{-1/2}>0$. Extend $T\varphi_m$ by zero outside
$[a,b]$, and define
\begin{equation}\label{cap:scalar-primitive}
p_m(r):=\frac1{r^3}\int_0^r t^3(T\varphi_m)(t)\,dt
\quad(r>0),\qquad p_m(0):=0.
\end{equation}
Then $p_m=0$ for $r\leq a$ and
\[
\frac1{r^3}\partial_r(r^3p_m)=T\varphi_m.
\]
The primitives have the explicit forms
\begin{equation}\label{cap:radial-primitives}
p_0(r)=\frac{H_T(r)}{r^3\sqrt{H_T(b)}},\qquad
p_m(r)=-\frac{T(r)\varphi_m'(r)}{\lambda_m}
\quad(a<r<b,\ m\geq1).
\end{equation}
For $r\geq b$, they satisfy
\[
p_0(r)=\frac{\sqrt{H_T(b)}}{r^3},\qquad p_m(r)=0\quad(m\geq1).
\]
\end{lemma}

\begin{proof}
The divergence identity follows from the defining integral.
For $m\geq1$, Lemma~\ref{c:lem:1dspec}(iii) gives
$-(r^3T\varphi_m')'=\lambda_m r^3T\varphi_m$ with zero endpoint
flux. Integration from $a$ to $r$ gives
\eqref{cap:radial-primitives}; zero flux at $b$ makes the full
integral vanish, so $p_m=0$ beyond $b$.
For $m=0$, insert the constant $\varphi_0=H_T(b)^{-1/2}$ in
\eqref{cap:scalar-primitive}. Since $H_T=0$ below $a$ and
$H_T(r)=H_T(b)$ above $b$, this gives the remaining formulas.
\end{proof}

Apply Lemma~\ref{cap:source-primitive-lemma} with $T=T_i$,
supported on $[r_a,r_b]$, and let $[z_a,z_b]$ be the support of $U_i$.
Recall the primitive in \eqref{c:scalar-norm}:
\begin{equation}\label{cap:scalar-functional-primitive}
H_{i,l}(r,z):=\int_0^r t^3\varrho_i(t,z)\Xi_{i,l}(t,z)\,dt
=r^3U_i(z)\sum_{(m,k)\in\Sigma_{i,1}}(\xi_{i,l})_{(m,k),1}p_m(r)\psi_k(z).
\end{equation}
Define the bounded and exterior parts of the prescribed bound
$B_{i,l}^2$ by
\begin{equation}\label{cap:scalar-functional-parts}
B_{i,l,\mathrm{in}}^2:=\nu_{i,l}\int_{z_a}^{z_b}\int_0^{r_b}\frac{|H_{i,l}(r,z)|^2}{r^3}\,dr\,dz,\qquad B_{i,l,\mathrm{ext}}^2:=\nu_{i,l}\int_{z_a}^{z_b}\int_{r_b}^{\infty}\frac{|H_{i,l}(r_b,z)|^2}{r^3}\,dr\,dz.
\end{equation}
Then \eqref{c:scalar-norm} is precisely
\begin{equation}\label{cap:scalar-functional-decomposition}
\|q_{(i,l)}\|_{Z'}^2\leq B_{i,l}^2=B_{i,l,\mathrm{in}}^2+B_{i,l,\mathrm{ext}}^2.
\end{equation}
To evaluate these parts, set
\begin{equation}\label{cap:primitive-gram-matrices}
J_{mn}=\int_0^{r_b}p_m(r)p_n(r)r^3\,dr,\qquad
K_{kk'}=\int_{z_a}^{z_b}U_i(z)^2\psi_k(z)\psi_{k'}(z)\,dz.
\end{equation}
Expanding the complete square in \eqref{cap:scalar-functional-parts} gives
\begin{equation}\label{cap:scalar-primitive-quadratic}
B_{i,l,\mathrm{in}}^2=\nu_{i,l}\sum_{\substack{(m,k)\in\Sigma_{i,1}\\(n,k')\in\Sigma_{i,1}}}
(\xi_{i,l})_{(m,k),1}(\xi_{i,l})_{(n,k'),1}J_{mn}K_{kk'}.
\end{equation}
For $r\geq r_b$, only $p_0$ survives, so
\begin{equation}\label{cap:scalar-terminal-primitive}
H_{i,l}(r,z)=H_{i,l}(r_b,z)=\sqrt{H_{T_i}(r_b)}\,U_i(z)
\sum_{k:\,(0,k)\in\Sigma_{i,1}}(\xi_{i,l})_{(0,k),1}\psi_k(z).
\end{equation}
Since $\int_{r_b}^{\infty}r^{-3}dr=(2r_b^2)^{-1}$, $0\leq U_i\leq1$, and the axial modes are orthonormal in $L^2(U_i\,dz)$, the exterior part satisfies
\begin{equation}\label{cap:coercivity-scalar-tail}
\begin{aligned}
B_{i,l,\mathrm{ext}}^2
&=\frac{\nu_{i,l}H_{T_i}(r_b)}{2r_b^2}
\int_{z_a}^{z_b}U_i^2
\left|\sum_{k:\,(0,k)\in\Sigma_{i,1}}
(\xi_{i,l})_{(0,k),1}\psi_k\right|^2\,dz\\
&\leq\frac{\nu_{i,l}H_{T_i}(r_b)}{2r_b^2}
\int_{z_a}^{z_b}U_i
\left|\sum_{k:\,(0,k)\in\Sigma_{i,1}}
(\xi_{i,l})_{(0,k),1}\psi_k\right|^2\,dz\\
&=\frac{\nu_{i,l}H_{T_i}(r_b)}{2r_b^2}
\sum_{k:\,(0,k)\in\Sigma_{i,1}}|(\xi_{i,l})_{(0,k),1}|^2.
\end{aligned}
\end{equation}
Thus the full source bound is
\begin{equation}\label{cap:scalar-functional-bound}
\|q_{(i,l)}\|_{Z'}\leq B_{i,l}\leq
\left(B_{i,l,\mathrm{in}}^2+\frac{\nu_{i,l}H_{T_i}(r_b)}{2r_b^2}
\sum_{k:\,(0,k)\in\Sigma_{i,1}}|(\xi_{i,l})_{(0,k),1}|^2\right)^{1/2}.
\end{equation}
The polynomial formula for $H_{T_i}$ in
Appendix~\ref{cap:appendix-primitive-arithmetic} and the exact-mode
models evaluate \eqref{cap:scalar-primitive-quadratic}. The matrix
$K$ always uses the whole axial support.

Evaluating the integrals in \eqref{cap:vector-functional-bound} and \eqref{c:scalar-norm}, including the exterior term of the scalar radial primitive, gives
\[
\begin{array}{c|cc}
 &\sum_{i,l} B_{i,l}&\sum_{i,l}B_{i,l}^2\\ \hline
\text{odd}&<18.877&<19.586\\
\text{even}&<33.792&<60.208
\end{array}
\]
These inequalities prove part~(vii) of Proposition~\ref{c:prop:CAP} and imply \eqref{c:qbounds}. The potential correction \eqref{c:profile-loss}, together with the
bound on $\|\mathcal R_\Psi\|_Z$ in Table~\ref{close:profiles}, gives
\begin{equation}\label{cap:coercivity-correction-margin}
\gamma_0-C_b\|\mathcal R_\Psi\|_Z\geq\gamma_0-1.001\times10^{-10}>\frac1{10}.
\end{equation}
This is the constant required by Proposition~\ref{close:proposition}.

\subsection{The response estimates}\label{sec:responses}\label{cap:functional-enclosures}

To prove the last two rows of Table~\ref{close:responses}, we estimate the errors in \eqref{close:response-identities} and the norms of the responses using their spline representations.

\subsubsection{Response functions and Poisson defects}\label{resp:construction}

Let $\Phi_{j,in},\Psi_{j,in},\wt\Phi_{j,in},\wt\Psi_{j,in}$
be the degree-$9$ tensor-product splines specified in \cite{proof_code}.
With the cutoff \eqref{res:cutoff}, we define
\begin{equation}\label{resp:representation}
\begin{aligned}
\Phi_j&=E\left(\frac{z\chi\Phi_{j,in}}{(1+R^2)^4}\right),&
\Psi_j&=E\left(\frac{z\chi\Psi_{j,in}}{(1+R^2)^2}\right),\\
\wt\Phi_{j,app}&=E\left(\frac{\chi\wt\Phi_{j,in}}{(1+R^2)^3}\right),&
\wt\Psi_j&=E\left(\frac{\chi\wt\Psi_{j,in}}{1+R^2}
 +(1-\chi)\wt\Psi_{j,ext}\right).
\end{aligned}
\end{equation}

Here $E$ is the lift \eqref{sf8}, and $\wt\Psi_{j,ext}$ is given
by \eqref{resp:multipoles-xt}, with prescribed real coefficients
$a_{j,0},a_{j,1},a_{j,2},a_{j,3}$.
The exterior expression satisfies
\[
|\wt\Psi_{j,ext}|=O(R^{-3}),\qquad
|\nabla\wt\Psi_{j,ext}|=O(R^{-4})\qquad(R\to\infty).
\]
Since $1-\chi=0$ on $S_{84}$, the factor $1-\chi$ removes the
singularity at the origin. Thus
$\Psi_j,\wt\Psi_j\in\dot H^1(\R^5)$.
The raw responses $\Phi_j,\wt\Phi_{j,app}$ are compactly supported.
For either raw response $f$, \eqref{robu1} and
$Sf=\operatorname{div}(yf)-5f$ give
\begin{equation}\label{resp:norm-identities}
\|f\|_Z^2=\|f\|_2^2+\|\nabla\Delta^{-1}f\|_2^2,
\qquad
\|Sf\|_{Z'}\leq\|Rf\|_2+5\|\nabla\Delta^{-1}f\|_2.
\end{equation}
The Hardy bound \eqref{res:radial-source} makes the right sides
finite, so both raw responses belong to $Z_0$ with their stated
parities. The normalized even response is
$\wt\Phi_j=P_0\wt\Phi_{j,app}$.

For fixed $\wt\mu_j\in\R$, define the defects using the
prescribed potentials by
\begin{equation}\label{resp:residuals}
\begin{aligned}
E_j={}&\mathcal O_{\R^5}\Phi_j
+\mathcal N_5(\Psi_5,\Phi_j)+\mathcal N_5(\Psi_j,\Omega^\ast)-q_j,\\
\wt E_j={}&(\mathcal O_{\R^5}-\wt\lambda_{app})\wt\Phi_{j,app}
+\mathcal N_5(\Psi_5,\wt\Phi_{j,app})
+\mathcal N_5(\wt\Psi_j,\Omega^\ast)
-\wt\mu_j\wt\Omega^\ast-\wt q_j,\\
E_j^\psi={}&\Delta\Psi_j-\Phi_j,
\qquad \wt E_j^\psi=\Delta\wt\Psi_j-\wt\Phi_{j,app}.
\end{aligned}
\end{equation}

The sources are the exact functions \eqref{c:qfinal}, using
Proposition~\ref{prop : spectral data} for the exact eigenfunctions
and Lemma~\ref{cap:local-mode-lemma} for their numerical plateau models.
Since the explicit potentials lie in $\dot H^1$, the Poisson
defects give the exact identities
\begin{equation}\label{resp:potential}
\Delta^{-1}\Phi_j=\Psi_j-\Delta^{-1}E_j^\psi,
\qquad
\Delta^{-1}\wt\Phi_{j,app}=\wt\Psi_j-\Delta^{-1}\wt E_j^\psi.
\end{equation}

Indeed, the difference between the two sides is harmonic in
$\dot H^1$ and therefore has zero gradient. These identities are
the starting point of both error estimates below.

\subsubsection{Odd response errors}\label{resp:odd-estimate}

Using $\Psi^\ast=\Psi_5-\Delta^{-1}\mathcal R_\Psi$ in
\eqref{resp:potential}, the exact odd error is
\begin{equation}\label{resp:odd-error-identity}
e_j=E_j-\mathcal L_3\Delta^{-1}E_j^\psi
       -\mathcal N_5(\Delta^{-1}\mathcal R_\Psi,\Phi_j).
\end{equation}

\begin{lemma}\label{resp:lemma}
The odd response and error satisfy
\begin{equation}\label{resp:bound}
\begin{aligned}
&\|e_j\|_{Z'}\leq\|\sqrt w\,E_j\|_2+\|r\Omega^\ast\|_\infty\|\sqrt w\,E_j^\psi\mathbf1_{\rm in}\|_2\\
&\quad+\left(I_o(\Omega^\ast,900)^2\|E_j^\psi\mathbf1_{\rm band}\|_2^2+\|r\Omega^\ast\|_{L^\infty(\R^5\setminus B_{900})}^2\|\sqrt w\,E_j^\psi\mathbf1_{\rm band}\|_2^2\right)^{1/2}+C_b\|\mathcal R_\Psi\|_Z\|\Phi_j\|_2,\\
&\|\Phi_j\|_Z^2\leq\|\Phi_j\|_2^2+\left(\|\nabla\Psi_j\|_2+\|\sqrt w\,E_j^\psi\|_2\right)^2.
\end{aligned}
\end{equation}

\end{lemma}

\begin{proof}
The source $E_j^\psi$ is odd and supported in $S_{100}$.
Split it into $S_{84}$ and $S_{band}$. Apply \eqref{res:near}
to the interior piece and \eqref{res:far-bound} to the band,
whose support is contained in $\{R\geq r_\ast\}$.
Equation~\eqref{nl4} bounds the last term of
\eqref{resp:odd-error-identity}. This proves the first inequality.
For the norm, \eqref{resp:potential} and the odd Hardy inequality give
\[
\|\nabla\Delta^{-1}\Phi_j\|_2
\leq\|\nabla\Psi_j\|_2+\|\sqrt w\,E_j^\psi\|_2.
\]
Insert this bound in \eqref{resp:norm-identities}.
\end{proof}

\subsubsection{The even projection and error}\label{resp:even-estimate}

Since $\wt e_j$ in \eqref{close:response-identities} is the error of the projected response $\wt\Phi_j=P_0\wt\Phi_{j,app}$, we also need a bound for the coefficient $m_0(\wt\Phi_{j,app})$ in \eqref{robu22}. This coefficient involves the corrected potential $\wPsi^\ast$ in addition to the prescribed potential $\wt\Psi_5$.

\begin{lemma}\label{resp:even-lemma}
Assume $\|\nabla\wPsi^\ast\|_2>0$. Then
\begin{equation}\label{resp:m0}
|m_0(\wt\Phi_{j,app})|
\leq\frac{
 |\langle\wt\Phi_{j,app},\wt\Psi_5\rangle| 
 {}+\big(\|\nabla\wt\Psi_j\|_2+\|r\wt E_j^\psi\|_2\big)
       \|r\mathcal R_{\wPsi}\|_2
 }
 {\|\nabla\wPsi^\ast\|_2^2}.
\end{equation}

For $\wt\Phi_j=P_0\wt\Phi_{j,app}$,
\begin{equation}\label{resp:even-assembly}
\begin{aligned}
\|\wt e_j\|_{Z'}&\leq{}\|r\wt E_j\|_2
 +\|r\Omega^\ast\|_\infty\|r\wt E_j^\psi\|_2
 +C_b\|\mathcal R_\Psi\|_Z\|\wt\Phi_{j,app}\|_2 +|m_0(\wt\Phi_{j,app})|\,\|\wt{\mathcal R}\|_{Z'},\\
\|\wt\Phi_j\|_{Z_0}&\leq{}
 \|\wt\Phi_{j,app}\|_Z+\delta\|S\wt\Phi_{j,app}\|_{Z'} +|m_0(\wt\Phi_{j,app})|
 \left(\|\wt\Omega^\ast\|_Z+\delta\|S\wt\Omega^\ast\|_{Z'}\right).
\end{aligned}
\end{equation}

Here $\delta=10^{-6}$. The unprojected norms follow from
\eqref{resp:norm-identities}, with $f=\wt\Phi_{j,app}$, and
\[
\|\nabla\Delta^{-1}\wt\Phi_{j,app}\|_2
\leq\|\nabla\wt\Psi_j\|_2+\|r\wt E_j^\psi\|_2.
\]
\end{lemma}

\begin{proof}
Since $\Delta\wPsi^\ast=\wt\Omega^\ast$, the denominator of
\eqref{robu21} is
$\langle\wt\Omega^\ast,\wPsi^\ast\rangle
=-\|\nabla\wPsi^\ast\|_2^2$.
Using $\wPsi^\ast=\wt\Psi_5-\Delta^{-1}\mathcal R_{\wPsi}$
and integrating by parts gives the signed identity
\begin{equation}\label{resp:normalizer-identity}
m_0(\wt\Phi_{j,app})
=-\frac{\langle\wt\Phi_{j,app},\wt\Psi_5\rangle
 +\langle\nabla\Delta^{-1}\wt\Phi_{j,app},
          \nabla\Delta^{-1}\mathcal R_{\wPsi}\rangle}
 {\|\nabla\wPsi^\ast\|_2^2}.
\end{equation}
Cauchy--Schwarz, \eqref{res:radial-source}, and
\eqref{resp:potential} imply \eqref{resp:m0}.
The quantities appearing in \eqref{resp:m0} satisfy
\begin{equation}\label{resp:normalizer-inputs}
\langle\wt\Phi_{j,app},\wt\Psi_5\rangle
=\mathcal I_{S_{100}}(\wt\Phi_{j,app},\wt\Psi_5),
\qquad\|\nabla\wPsi^\ast\|_2^2\geq0.88^2>0,
\end{equation}
where the second inequality uses Table~\ref{close:profiles}. The pairing
is enclosed as the sum of the signed cell integrals.

For the equation error, use
$\mathcal S_{\wt\lambda_{app}}\wt\Omega^\ast=-\wt{\mathcal R}$
and \eqref{resp:potential}. The projection produces
\begin{equation}\label{resp:even-error-identity}
\wt e_j={}\wt E_j-\mathcal L_3\Delta^{-1}\wt E_j^\psi
 -\mathcal N_5(\Delta^{-1}\mathcal R_\Psi,\wt\Phi_{j,app}) +m_0(\wt\Phi_{j,app})\,\wt{\mathcal R}.
\end{equation}
Apply the even case of \eqref{res:near} and \eqref{nl4} to obtain
the error bound in \eqref{resp:even-assembly}.
The norm bound follows from \eqref{resp:norm-identities} and
$\wt\Phi_j=\wt\Phi_{j,app}-m_0(\wt\Phi_{j,app})\wt\Omega^\ast$.
This also gives $\wt\Phi_j\in Z_{0,\ast}^{even}$.
\end{proof}

In particular the normalized response generally has an exterior:
\begin{equation}\label{resp:projected-exterior}
\wt\Phi_j^{\rm ext}
=-m_0(\wt\Phi_{j,app})\wt\Omega_{\geq-9}.
\end{equation}
Its contribution is included in the global profile norms in
\eqref{resp:even-assembly}. The raw response
$\wt\Phi_{j,app}$ is the compactly supported field.

\subsubsection{Complete response integrals}\label{resp:global-integrals}

The odd fields $\Phi_j,\Psi_j,E_j,E_j^\psi$ and the even fields $\wt\Phi_{j,app},\wt E_j^\psi$ vanish outside $S_{100}$. Hence their norms in Lemmas~\ref{resp:lemma} and~\ref{resp:even-lemma} are cell integrals; for instance, $\|\nabla\Psi_j\|_2^2=\mathcal I_{S_{100}}(\partial_r\Psi_j)+\mathcal I_{S_{100}}(\partial_z\Psi_j)$ and, since $R^2=r^2+z^2$, $\|R\wt\Phi_{j,app}\|_2^2=\mathcal I_{S_{100}}(r\wt\Phi_{j,app})+\mathcal I_{S_{100}}(z\wt\Phi_{j,app})$. The weighted norms are bounded with the cell weights $w_C=\sup_Cw$ as in \eqref{cap:stream-integrals}, separately on $S_{84}$ and $S_{band}$.

The even potential $\wt\Psi_j$ is not compactly supported. Outside
$S_{100}$ its equation residual and Poisson defect are
\begin{equation}\label{resp:exterior}
\wt E_j^{\rm ext}
=\mathcal N_5(\wt\Psi_{j,ext},\Omega_{\geq-11})
-\wt\mu_j\wt\Omega_{\geq-9},\qquad
(\wt E_j^\psi)^{\rm ext}=0.
\end{equation}
Consequently,
\begin{equation}\label{resp:global-even-integrals}
\begin{aligned}
\|r\wt E_j\|_2^2
&\leq\mathcal I_{S_{100}}(r\wt E_j)
 +\mathfrak I_{2,0}(\wt E_j^{\rm ext}),\\
\|\nabla\wt\Psi_j\|_2^2
&\leq\mathcal I_{S_{100}}(\partial_r\wt\Psi_j)
 +\mathcal I_{S_{100}}(\partial_z\wt\Psi_j)
 +\mathfrak I_{2,0}(t^2\mathcal A_r\wt\Psi_{j,ext})
 +\mathfrak I_{0,2}(t^2\mathcal A_z\wt\Psi_{j,ext}).
\end{aligned}
\end{equation}
The finite polynomial for $\wt E_j^{\rm ext}$ is
\eqref{resp:exterior-polynomial}. Both exterior contributions in
\eqref{resp:global-even-integrals} enter before square roots are
taken. Appendix~\ref{cap:appendix-response-arithmetic} gives the
local formulas used to evaluate these integrals.

\subsubsection{The source weights and final sums}\label{resp:weighted-sums}

The $B_{i,l}$ bounds used for coercivity also bound the norms of
the exact sources \eqref{c:qfinal}. For the response sums we may
improve an individual weight using the scalar Hardy estimate.
If $j$ enumerates $(i,l)$, define
\begin{equation}\label{resp:source-weights}
\begin{aligned}
\beta_j&:=\min\left\{B_{i,l},
 \left(\sum_{C\subseteq S_{100}}w_C\mathcal I_C(q_j)\right)^{1/2}
 \right\},&&\|q_j\|_{Z'}\leq\beta_j,\\
\wt\beta_j&:=\min\left\{B_{i,l},
 \mathcal I_{S_{100}}(r\wt q_j)^{1/2}\right\},
 &&\|\wt q_j\|_{Z'}\leq\wt\beta_j.
\end{aligned}
\end{equation}
All sources are supported in $S_{100}$. With these weights, the interval enclosures for Lemmas~\ref{resp:lemma} and~\ref{resp:even-lemma} give
\begin{equation}\label{resp:final-comparisons}
\begin{aligned}
\sum_{j=1}^{26}\beta_j\|e_j\|_{Z'}&<4.790\times10^{-4},&
\sum_{j=1}^{32}\wt\beta_j\|\wt e_j\|_{Z'}&<7.131\times10^{-5},\\
\sum_{j=1}^{26}\beta_j\|\Phi_j\|_Z&<15.244,&
\sum_{j=1}^{32}\wt\beta_j\|\wt\Phi_j\|_{Z_0}&<141.786.
\end{aligned}
\end{equation}
Since $\|q_j\|_{Z'}\leq\beta_j$ and $\|\wt q_j\|_{Z'}\leq\wt\beta_j$, these inequalities complete the proof of Table~\ref{close:responses}.

\appendix

\section{Computation of the eigenfunctions}\label{app:eigenfunctions}

The cell construction of Section~\ref{sec:coercivity} requires the eigenvalues and the normalized eigenfunctions of a weighted Neumann problem. In particular, the numerical approximations must be assigned the correct mode indices: an enclosure of a zero of the differential equation does not by itself identify its position in the spectrum. We construct the regular solutions on overlapping series pieces and match them through a finite system of equations, whose zero we enclose by a contraction argument. Counting the zeros of the resulting eigenfunction identifies the mode, and a final normalization yields the spectral data used in Section~\ref{sec:residual}.

\subsection{Regular solutions and the matching equations}

Let $P[a,x_2,x_3,b]$ be the plateau function of \eqref{c:plateau}. We consider either
\[
\rho(x)=x^3P[a,x_2,x_3,b](x),\qquad a>0,
\]
or $\rho(x)=P[a,x_2,x_3,b](x)$, and seek a nonzero solution of
\begin{equation}\label{eq : eigenvalue problem with rho}
-(\rho\varphi')'=\lambda\rho\varphi\quad\text{on }(a,b),
\qquad \varphi'(a)=\varphi'(b)=0.
\end{equation}
The regularity, simplicity and ordering of these eigenpairs were established in Lemma~\ref{c:lem:1dspec}. Here we give enclosures which retain that ordering. The mode $m=0$ is known exactly:
\[
\lambda_0=0,\qquad
\varphi_0(x)=\left(\int_a^b\rho(x)\,dx\right)^{-1/2}.
\]
The integral is an explicit polynomial integral. We therefore restrict the validation below to $\lambda>0$.

On the plateau $[x_2,x_3]$, equation~\eqref{eq : eigenvalue problem with rho} is explicit. For the radial weight $x^3P$ its solutions are linear combinations of
\[
u_\lambda(x)=\frac{J_1(\sqrt\lambda x)}x,
\qquad v_\lambda(x)=\frac{Y_1(\sqrt\lambda x)}x,
\]
whereas for the weight $P$ we take
\[
u_\lambda(x)=\cos(\sqrt\lambda x),\qquad
v_\lambda(x)=\sin(\sqrt\lambda x).
\]

At $a$ and $b$, the weight vanishes to order five. We use a series centered at each endpoint to select the regular solution, and a second series on each remaining part of the ramp. More precisely, choose
\[
a<x_1<x_2\le x_3<x_4<b,
\]
and put $\delta_a=x_1-a$, $\delta_b=b-x_4$, $\bar x_a=(x_1+x_2)/2$, $w_a=(x_2-x_1)/2$, $\bar x_b=(x_3+x_4)/2$, and $w_b=(x_4-x_3)/2$. The decomposition is
\begin{equation}\label{eq:decomp}
\varphi(x)=
\begin{cases}
\displaystyle\varphi_a(x)=\sum_{n\ge0}a_n\left(\frac{x-a}{\delta_a}\right)^n,&a\le x\le x_1,\\[6pt]
\displaystyle\psi_a(x)=\sum_{n\ge0}c^a_n\left(\frac{x-\bar x_a}{w_a}\right)^n,&x_1\le x\le x_2,\\[6pt]
A u_\lambda(x)+B v_\lambda(x),&x_2\le x\le x_3,\\[4pt]
\displaystyle\psi_b(x)=\sum_{n\ge0}c^b_n\left(\frac{x-\bar x_b}{w_b}\right)^n,&x_3\le x\le x_4,\\[6pt]
\displaystyle\varphi_b(x)=\sum_{n\ge0}b_n\left(\frac{b-x}{\delta_b}\right)^n,&x_4\le x\le b.
\end{cases}
\end{equation}
The endpoint series use the local variable $y\in[0,1]$, while the other two series use $y\in[-1,1]$. On the right endpoint piece the physical derivative is $-\delta_b^{-1}\partial_y$; this sign is retained in the matching conditions.

At an endpoint, we write
\[
\rho(a+\delta_a y)=\delta_a^5y^5Q(y)
\quad\text{or}\quad
\rho(b-\delta_b y)=\delta_b^5y^5Q(y),
\]
where $Q(0)\ne0$, and expand $q=Q'/Q=\sum_{j\ge0}q_jy^j$. In the local variable, the equation for the endpoint series is
\[
y\varphi_u''+(5+yq)\varphi_u'+\delta_u^2\lambda y\varphi_u=0,
\qquad u\in\{a,b\},
\]
where primes in this display denote differentiation with respect to $y$. Thus, with $u_n=a_n$ or $b_n$,
\begin{equation}\label{eq : equations for un}
u_1=0,\qquad
n(n+4)u_n=-\sum_{j=0}^{n-2}(n-1-j)q_j u_{n-1-j}
                 -\delta_u^2\lambda u_{n-2},\qquad n\ge2.
\end{equation}
The value $u_0$ is free. We denote the solution sequence by $\mathcal S^u(u_0,\lambda)$; linearity gives
\begin{equation}\label{def : solution operator Su}
\mathcal S^u(u_0,\lambda)=u_0\mathcal S^u(1,\lambda).
\end{equation}

On a centered ramp piece, $\rho(\bar x_u+w_uy)$ is nonzero near $[-1,1]$. With
\[
q(y)=\frac{\partial_y\rho(\bar x_u+w_uy)}{\rho(\bar x_u+w_uy)}
     =\sum_{j\ge0}q_jy^j,
\]
the corresponding recurrence is
\begin{equation}\label{eq : equations for cn}
n(n-1)c^u_n=-\sum_{j=0}^{n-2}(n-1-j)q_jc^u_{n-1-j}
                    -w_u^2\lambda c^u_{n-2},\qquad n\ge2.
\end{equation}
Here $c^u_0,c^u_1$ are free. Its solution operator satisfies
\begin{equation}\label{def : solution operator Cu}
\mathcal C^u(c^u_0,c^u_1,\lambda)
=c^u_0\mathcal C^u(1,0,\lambda)+c^u_1\mathcal C^u(0,1,\lambda).
\end{equation}
Consequently, for fixed $\lambda$ each piece in \eqref{eq:decomp} is determined by its displayed amplitudes.

It remains to match the values and the physical first derivatives at the four junctions. Let
\[
\theta=(\lambda,a_0,c^a_0,c^a_1,A,B,c^b_0,c^b_1,b_0)\in\R^9.
\]
Since \eqref{eq : eigenvalue problem with rho} is linear, we may normalize the solution by $a_0=1$. The matching map is
\begin{equation}\label{app:matching-map}
g(\theta)=
\begin{pmatrix}
\varphi_a(x_1)-\psi_a(x_1)\\
\psi_a(x_2)-A u_\lambda(x_2)-B v_\lambda(x_2)\\
A u_\lambda(x_3)+B v_\lambda(x_3)-\psi_b(x_3)\\
\psi_b(x_4)-\varphi_b(x_4)\\
\varphi_a'(x_1)-\psi_a'(x_1)\\
\psi_a'(x_2)-A u_\lambda'(x_2)-B v_\lambda'(x_2)\\
A u_\lambda'(x_3)+B v_\lambda'(x_3)-\psi_b'(x_3)\\
\psi_b'(x_4)-\varphi_b'(x_4)\\
a_0-1
\end{pmatrix}.
\end{equation}
It is affine in the last eight coordinates:
\begin{equation}\label{eq : g affine}
g(\theta)=L(\lambda)(\theta_2,\ldots,\theta_9)^T-e_9,
\qquad L(\lambda)\in\R^{9\times8}.
\end{equation}
The series estimates in Subsection~\ref{sec : solution operators} justify two termwise differentiations. Thus every zero of $g$ gives a nonzero $C^2$ solution of \eqref{eq : eigenvalue problem with rho}. Indeed, the matching makes $\varphi$ and $\varphi'$ continuous, and the differential equation then matches $\varphi''$ at the interior junctions. The endpoint series give $\varphi'(a)=\varphi'(b)=0$ and $\varphi(a)=1$.

\subsection{Uniform bounds for the series solutions}\label{sec : solution operators}

To validate the matching map, we need both the solution operators and their derivatives with respect to $\lambda$ on an interval. The recurrences give the finite coefficients directly. We now bound the remaining coefficients, uniformly in that interval.

First we enclose $q$. Since the ramp cutoff is the polynomial in \eqref{c:S}, the nonzero roots of the endpoint polynomial $Q$ consist of the affine images of the four roots of
\[
126-420t+540t^2-315t^3+70t^4,
\]
and, for a radial weight, the triple root corresponding to $x=0$. For a centered ramp piece one must also include the order-five root at the adjacent endpoint. The four roots of the displayed quartic are enclosed in disjoint complex discs by the same contraction argument used below in Lemma~\ref{lem : fixed point g}. Since the degree is four, these are all of its roots, so that we obtain the complete list of roots for each polynomial weight.

Choose $\nu>1$ smaller than the modulus of every root $\zeta$ in that list. Counting multiplicities, the logarithmic derivative gives
\[
q_j=-\sum_\zeta\zeta^{-j-1},
\]
and therefore
\begin{equation}\label{eq : G}
\sum_{j\ge s}|q_j|\nu^j\le
G_s:=\sum_\zeta\frac{(\nu/|\zeta|)^s}{|\zeta|-\nu},\qquad s\ge0.
\end{equation}
The root enclosures give interval upper bounds for the right-hand side. In particular, the regular coefficient $q$ is analytic on a disc larger than the local interval of evaluation.

We use the weighted sequence space
\[
\ell^1_\nu=\left\{U=(u_n)_{n\ge0}:\|U\|_\nu:=\sum_{n\ge0}|u_n|\nu^n<\infty\right\},
\]
and write $\pi^{\le N}$ and $\pi^{>N}$ for its coefficient projections. For the endpoint recurrence with $u_0=1$, define
\begin{equation}\label{eq : f lambda}
(f_\lambda(U))_n=
\begin{cases}
u_0-1,&n=0,\\
u_1,&n=1,\\
\displaystyle u_n+\frac{\sum_{k=1}^{n-1}kq_{n-1-k}u_k+\delta_u^2\lambda u_{n-2}}{n(n+4)},&n\ge2.
\end{cases}
\end{equation}
Its linear part is $I+K$, where $K$ is strictly lower triangular. For each $\lambda$, the recurrence \eqref{eq : equations for un} determines the finite sequence
\[
\bar U=(1,0,u_2,\ldots,u_N,0,\ldots).
\]
Its first $N+1$ equations are satisfied exactly. The omitted residual is bounded by
\begin{equation}\label{def : Y infinity}
\begin{aligned}
\|\pi^{>N}f_\lambda(\bar U)\|_\nu
\le Y_\infty(\bar U):={}&
\frac{\nu}{(N+1)(N+5)}\sum_{k=1}^N k|u_k|\nu^kG_{N-k}\\
&+\delta_u^2|\lambda|\left(
\frac{\nu^{N+1}|u_{N-1}|}{(N+1)(N+5)}+
\frac{\nu^{N+2}|u_N|}{(N+2)(N+6)}\right).
\end{aligned}
\end{equation}

\begin{lemma}\label{lem : radii recurrence}
Let $N\ge2$, let $q$ satisfy \eqref{eq : G}, and suppose
\[
Z_1:=\frac{\nu G_0}{N+5}
      +\frac{\delta_u^2|\lambda|\nu^2}{(N+1)(N+5)}<1.
\]
Then the endpoint recurrence with $u_0=1$ has a unique solution $U^*\in\ell^1_\nu$, and
\[
\pi^{\le N}U^*=\bar U,\qquad
\|U^*-\bar U\|_\nu\le r_U:=\frac{Y_\infty(\bar U)}{1-Z_1}.
\]
For the centered recurrence \eqref{eq : equations for cn}, the same conclusion holds for either initial pair $(1,0)$ or $(0,1)$, with $\delta_u$ replaced by $w_u$, and with $N+5,N+6$ replaced by $N,N+1$ in these bounds.
\end{lemma}

\begin{proof}
The finite recurrence fixes the first $N+1$ coefficients uniquely. On $\pi^{>N}\ell^1_\nu$, the remaining equations are equivalent to
\[
\pi^{>N}U=-\pi^{>N}K\bar U-\pi^{>N}K\pi^{>N}U.
\]
For a column indexed by $k$, the coefficient involving $q_j$ is $kq_j/[n(n+4)]$, with $n=k+1+j>N$. Since $k\le n-1$,
\[
\frac{k}{n(n+4)}\le\frac1{N+5},\qquad
\frac1{n(n+4)}\le\frac1{(N+1)(N+5)}.
\]
The weighted column sums, together with \eqref{eq : G}, give
\[
\|\pi^{>N}K\pi^{>N}\|_\nu\le Z_1.
\]
For the finite forcing $\pi^{>N}K\bar U$, the summation for the $k$-th coefficient starts at $j=N-k$; the terms containing $\lambda$ occur only at $n=N+1,N+2$. This gives \eqref{def : Y infinity}. The affine map on the tail is a contraction, and its fixed point has norm at most $Y_\infty/(1-Z_1)$. It gives the unique full solution, since every solution must have the same finite coefficients. For the centered recurrence the denominators are $n(n-1)$, so the two bounds become $1/N$ and $1/[N(N+1)]$, respectively.
\end{proof}

When $\lambda$ ranges over an interval, forward substitution with interval arithmetic encloses the first $N+1$ coefficients for every such $\lambda$, and \eqref{def : Y infinity} and $Z_1$, evaluated with intervals, give a uniform tail bound.

The weighted norm also controls physical derivatives. For $i=0,1,2$, let
\[
c_0=1,\qquad c_i=\left(\frac{i}{e\log\nu}\right)^i\quad(i=1,2).
\]
Since $n!/(n-i)!\le n^i$, every $V\in\ell^1_\nu$ satisfies
\begin{equation}\label{eq : evaluation}
\sup_{|y|\le1}\left|\sum_{n\ge i}\frac{n!}{(n-i)!}V_ny^{n-i}\right|
\le c_i\|V\|_\nu.
\end{equation}
On a physical series piece of scale $h\in\{\delta_a,w_a,w_b,-\delta_b\}$, the corresponding bound has the additional factor $|h|^{-i}$. In particular, the series constructed above converge with two derivatives on their closed pieces.

We also need the $\lambda$-derivatives of these solutions. The coefficients of the endpoint solution satisfy
\begin{equation}\label{eq : recurrence dot u}
\begin{aligned}
n(n+4)\partial_\lambda u_n={}&-\sum_{k=1}^{n-1}kq_{n-1-k}\partial_\lambda u_k
-\delta_u^2\lambda\partial_\lambda u_{n-2}-\delta_u^2u_{n-2},\qquad n\ge2,\\
\partial_\lambda u_0={}&\partial_\lambda u_1=0.
\end{aligned}
\end{equation}
The finite coefficients are again obtained by forward substitution. Analytic dependence on $\lambda$ follows from the uniformly convergent Neumann series for the tail operator. The omitted part of $\partial_\lambda U^*$ is bounded by
\begin{equation}\label{app:parameter-tail}
\begin{aligned}
\|\partial_\lambda U^*-\partial_\lambda\bar U\|_\nu
\le\frac1{1-Z_1}\Bigg[&Y_\infty(\partial_\lambda\bar U)
+\delta_u^2\left(
\frac{\nu^{N+1}|u_{N-1}|}{(N+1)(N+5)}+
\frac{\nu^{N+2}|u_N|}{(N+2)(N+6)}\right)\\
&+\frac{\delta_u^2\nu^2 r_U}{(N+1)(N+5)}\Bigg].
\end{aligned}
\end{equation}
Indeed, the first term bounds the homogeneous recurrence residual, the next two finite terms are the contributions of $\partial_\lambda f_\lambda(\bar U)$ at $N+1,N+2$, and the last term bounds the forcing from $U^*-\bar U$. The centered pieces have the same replacements as in Lemma~\ref{lem : radii recurrence}. Combining \eqref{app:parameter-tail} with \eqref{eq : evaluation} gives the values and the physical first derivatives of the solution operators and their $\lambda$-derivatives at every junction.

\subsection{Validation of the matching equations}\label{sec : Dg}

We can now enclose the matching map and its derivative. From \eqref{eq : g affine},
\[
Dg(\theta)=
\big(\,\partial_\lambda L(\lambda)(\theta_2,\ldots,\theta_9)^T\ \big|\ L(\lambda)\,\big).
\]
The series entries of $L$ and $\partial_\lambda L$ are supplied by the previous subsection. For the plateau entries we differentiate the explicit formulas. In the radial case they follow from the recurrences of the Bessel functions~\cite[\S10.6]{DLMF}.  The trigonometric case is elementary, and the second physical derivatives follow from the differential equation.

Point values of the Bessel functions are enclosed by arbitrary-precision ball arithmetic \cite{Arb}. Interval arguments are handled by the mean value theorem and the identities $2J_n'=J_{n-1}-J_{n+1}$ and $2Y_n'=Y_{n-1}-Y_{n+1}$. The required derivative bounds follow from the fact that $J_n(x)^2+Y_n(x)^2$ decreases on $(0,\infty)$, a consequence of Nicholson's integral representation \cite[\S10.9]{DLMF}. This also gives uniform bounds on the entire plateau; in the trigonometric case we use $|\sin|,|\cos|\le1$.

A zero of $g$ is then obtained from the following finite-dimensional criterion.
\begin{lemma}\label{lem : fixed point g}
Let $\bar\theta\in\R^9$, let $M\in\R^{9\times9}$, and suppose $g$ is continuously differentiable on a neighborhood of the closed ball $B_{r^*}(\bar\theta)$. Assume
\[
\|Mg(\bar\theta)\|_\infty\le\epsilon,\qquad
\sup_{\theta\in B_{r^*}(\bar\theta)}\|I-MDg(\theta)\|_\infty\le Z<1.
\]
If $r=\epsilon/(1-Z)\le r^*$, then $g$ has a unique zero in $B_{r^*}(\bar\theta)$, and this zero belongs to $B_r(\bar\theta)$.
\end{lemma}
\begin{proof}
The derivative of $\theta\mapsto\theta-Mg(\theta)$ has norm at most $Z$. Moreover, for $r\le s\le r^*$ and $\theta\in B_s(\bar\theta)$,
\[
\|\theta-Mg(\theta)-\bar\theta\|_\infty\le Zs+\epsilon\le s.
\]
The contraction mapping theorem applies on both closed balls of radii $r$ and $r^*$. The resulting fixed point is unique in the larger ball. Finally, $MDg(\bar\theta)$ is invertible since its distance from the identity is less than one, so $M$ is invertible and this fixed point is a zero of $g$.
\end{proof}

We choose $\bar\theta$ as a numerical solution of the matching equations and $M$ as a numerical inverse of $Dg(\bar\theta)$. These numerical choices only specify the point and the preconditioner. The quantities $\epsilon$ and $Z$ are then enclosed with intervals. To compute $Z$, we apply Lemma~\ref{lem : radii recurrence} and \eqref{app:parameter-tail} uniformly for
\[
\lambda\in[\bar\lambda-r^*,\bar\lambda+r^*],
\qquad
\theta_j\in[\bar\theta_j-r^*,\bar\theta_j+r^*],\quad2\le j\le9.
\]
This gives an enclosure of every entry of $Dg$ on the required ball. Lemma~\ref{lem : fixed point g} thus supplies an eigenvalue with
\[
|\lambda-\bar\lambda|\le r.
\]

\subsection{Normalization and identification of the eigenvalue}\label{sec : normalization}

Let $\theta^*$ be the zero just obtained, and let $\varphi$ be its eigenfunction with $\varphi(a)=1$. Let $\tilde\varphi$ be the piecewise approximation given by fixed numerical polynomials on the four series pieces and by $\bar A u_{\bar\lambda}+\bar B v_{\bar\lambda}$ on the plateau. The exact eigenfunction $\varphi$ is enclosed by evaluating the solution operators for $\lambda\in[\bar\lambda-r,\bar\lambda+r]$ and the amplitudes in the corresponding parameter intervals.

For instance, on the left centered piece the weighted coefficient error is bounded by
\[
\big\|\theta^*_3\pi^{\le N}\mathcal C^a(1,0,\lambda)
      +\theta^*_4\pi^{\le N}\mathcal C^a(0,1,\lambda)-\tilde\varphi\big\|_\nu
+|\theta^*_3|r_{\mathcal C^a(1,0)}+|\theta^*_4|r_{\mathcal C^a(0,1)}.
\]
Here each radius is furnished by Lemma~\ref{lem : radii recurrence}, uniformly over the eigenvalue enclosure. Formula~\eqref{eq : evaluation} converts this into errors for the value and the first two physical derivatives. The other series pieces are treated in the same way. On the plateau, for $i=0,1,2$, the mean value theorem gives
\[
\sup_{[x_2,x_3]}|\varphi^{(i)}-\tilde\varphi^{(i)}|
\le r\sum_{w\in\{u,v\}}\sup_{\substack{x\in[x_2,x_3]\\|\lambda-\bar\lambda|\le r}}|w_\lambda^{(i)}(x)|
+r\left(|\bar A|\sup|\partial_\lambda u_\lambda^{(i)}|
       +|\bar B|\sup|\partial_\lambda v_\lambda^{(i)}|\right),
\]
where the last two suprema are over the same set.

We first identify the mode using the unnormalized eigenfunction. Let $\epsilon$ be an upper bound for the errors of its value and first derivative on all pieces. Subdivide the five physical pieces and bisect until each closed interval $I$ satisfies at least one of
\[
\tilde\varphi>\epsilon,\qquad
\tilde\varphi<-\epsilon,\qquad
\tilde\varphi'>\epsilon,\qquad
\tilde\varphi'<-\epsilon
\quad\text{throughout }I.
\]
Since these inequalities are verified in interval arithmetic and the errors are at most $\epsilon$, the true eigenfunction has the sign of $\tilde\varphi$ on an interval of the first two kinds, and its derivative has the sign of $\tilde\varphi'$ on an interval of the last two kinds. Consecutive derivative intervals are joined only when this derivative sign agrees. The first and the last intervals must be sign intervals, and adjacent sign intervals with opposite signs must be separated by a derivative interval. Every maximal derivative interval is therefore a region where the true eigenfunction is strictly monotone. It contains exactly one zero when the neighboring signs are opposite, and none when they agree. These are all the zeros.

The approximation need not agree exactly across its numerical junctions. The argument applies to the true $C^1$ eigenfunction, which has the same value and derivative error bounds on both sides. In particular, a derivative interval crossing a junction is still a region of strict monotonicity when the recorded signs agree. If the resulting count is $m$, Lemma~\ref{c:lem:1dspec} identifies the eigenvalue as $\lambda_m$.

It remains to impose weighted normalization. Define
\[
\mathcal N=\int_a^b\rho\varphi^2\,dx.
\]
On a series piece, the integral of the polynomial approximation is evaluated exactly by integrating its coefficients, with its physical Jacobian. If $|\varphi-\tilde\varphi|\le e$ on the piece, its error is at most
\[
\left(2e\sup|\tilde\varphi|+e^2\right)\int\rho\,dx.
\]
On the plateau, the true solution has the explicit amplitudes enclosed by $\theta^*$, and $\rho\varphi^2$ has an explicit primitive, given by Lommel's integral~\cite[\S10.22]{DLMF} in the radial case and by an elementary formula in the trigonometric case. Evaluating these primitives at $x_2,x_3$ and adding the four series integrals encloses $\mathcal N$. The computation checks that its lower endpoint is positive.

The normalized eigenfunction is $\varphi_m=\varphi/\sqrt{\mathcal N}$. Choose a numerical approximation $\tilde s$ to $\mathcal N^{-1/2}$ and round the coefficients of $\tilde s\tilde\varphi$ to specify $\tilde\psi_m$. For each physical derivative,
\[
\sup|\varphi_m^{(i)}-\tilde\psi_m^{(i)}|
\le\mathcal N^{-1/2}\sup|\varphi^{(i)}-\tilde\varphi^{(i)}|
   +\sup|\mathcal N^{-1/2}\tilde\varphi^{(i)}-\tilde\psi_m^{(i)}|,
\qquad i=0,1,2,
\]
where the second term includes coefficient rounding and the enclosure of the normalization factor. Since $\varphi(a)=1$, this normalization gives $\varphi_m(a)>0$.

When the weight is even, Lemma~\ref{c:lem:1dspec} gives $\varphi_m(-x)=(-1)^m\varphi_m(x)$. We replace the numerical approximation by
\[
\frac{\tilde\psi_m(x)+(-1)^m\tilde\psi_m(-x)}2,
\]
and round its coefficients again. Parity projection does not enlarge any of the three uniform derivative error bounds; the final rounding errors are then added. In this way the delivered approximation has the exact parity of the true mode.

\subsection{The spectral data}\label{sec : results}

We apply this procedure to the weights in Table~\ref{tab:current-weights}. 

\begin{table}[htbp]
\centering
\caption{Weights of the cell construction and the highest validated mode index. The constant mode is known exactly.}
\label{tab:current-weights}
\begin{tabular}{lcr}
\toprule
weight & plateau parameters $(a,x_2,x_3,b)$ & highest index\\
\midrule
$T_0r^3$ &(0.5,1,15,17)&19\\
$T_1r^3=S_1r^3$ &(15,17,29,31)&21\\
$T_2r^3$ &(29,31,42,45)&6\\
$S_0r^3$ &(0.1,0.2,15,17)&18\\
$S_2r^3$ &(29,31,152,160)&9\\
$U_0$ &(-7,-5,5,7)&18\\
$U_1$ &(5,7,20,23)&8\\
$V_0$ &(-40,-35,35,40)&84\\
$V_1$ &(-200,-180,180,200)&24\\
\bottomrule
\end{tabular}
\end{table}

\begin{proposition}\label{prop : spectral data}
Let $\rho$ be a weight in Table~\ref{tab:current-weights} and $1\le m\le M$, where $M$ is its highest index. The $m$-th eigenvalue of \eqref{eq : eigenvalue problem with rho} and the eigenfunction normalized by $\int_a^b\rho\varphi_m^2\,dx=1$ and $\varphi_m(a)>0$ satisfy
\[
|\lambda_m-\bar\lambda_m|\le3.3\times10^{-31},\qquad
\sup_{x\in I}|\varphi_m^{(i)}(x)-\tilde\psi_m^{(i)}(x)|\le1.1\times10^{-24},\quad i=0,1,2,
\]
for every closed piece $I$ of \eqref{eq:decomp}. At a junction, the derivative bounds hold for the approximations from both adjoining pieces.
\end{proposition}
\begin{proof}
There are $207$ nonconstant eigenpairs in the table. Their interval bounds satisfy
\[
Z\le4.1\times10^{-12}<1,\qquad
\frac{\epsilon}{1-Z}\le3.3\times10^{-31}<10^{-20}=:r^*.
\]
Lemma~\ref{lem : fixed point g} gives an eigenfunction and its eigenvalue enclosure. The zero count identifies the eigenvalue as $\lambda_m$. The positive normalization integral and the derivative error estimates in Subsection~\ref{sec : normalization}, including parity projection and coefficient rounding, give the three stated derivative bounds. The sign follows from the initial condition $a_0=1$.
\end{proof}

For the nine constant modes the positive normalized eigenfunction is $(\int_a^b\rho(x)\,dx)^{-1/2}$, which we enclose directly.

\section{Auxiliary formulas for the verification}\label{app:verification-formulas}

\subsection{Formulas for the drift and test fields}\label{app:profile-arithmetic}

\subsubsection{The drift}

Using \eqref{u:cov} and \eqref{u:norm}, we obtain
\begin{equation}\label{cap:test-field-drift-model}
\begin{aligned}
r^2(\partial_z\Psi_5)^2+(r\partial_r\Psi_5+2\Psi_5)^2
={}&r^2\left(\Psi_{app}^{(c)}
        +\frac{6z^2}{(1+v)^3}\partial_v\Psi_{app}^{(c)}\right)^2\\
&+z^2\left(2\Psi_{app}^{(c)}
        +\frac{6r^2}{(1+u)^3}\partial_u\Psi_{app}^{(c)}\right)^2.
\end{aligned}
\end{equation}

For $L_e=\sqrt{U_0(100)}$, \eqref{cap:asymptotic-derivatives} and
\eqref{cap:exterior-supremum} give the following bound, with $I$
the identity and the odd quotient understood as in
\eqref{cap:asymptotic-polynomial}:
\begin{equation}\label{cap:exterior-drift-supremum}
\begin{aligned}
&\|(-r\partial_z\Psi_{\geq-9},
 r\partial_r\Psi_{\geq-9}+2\Psi_{\geq-9})\|_{L^\infty(S_{\rm ext})}^2\\
&\quad\leq
\mathfrak S_{1,0}\left(
 (I+(1-x)\mathcal A_z)\frac{\Psi_{\geq-9}}z;L_e\right)^2
+\mathfrak S_{0,1}\left(
 (2I+x\mathcal A_r)\frac{\Psi_{\geq-9}}z;L_e\right)^2.
\end{aligned}
\end{equation}

\subsubsection{The test-field numerators}

The numerators in \eqref{cap:test-field-recovery} are obtained by inserting \eqref{u:cov} and \eqref{u:norm} in \eqref{cap:test-fields}. They are polynomials, on each cell, in the spline pieces of $\Psi_{app}^{(c)},\wt\Psi_{app}^{(c)},\Omega_{app}^{(c)},\wt\Omega_{app}^{(c)}$ and their first derivatives.

On the exterior, \eqref{cap:asymptotic-derivatives} gives
\begin{equation}\label{cap:exterior-test-field-polynomials}
\begin{aligned}
a_5^{\rm ext}={}&\wt\Omega_{\geq-9}
 +\left((I+(1-x)\mathcal A_z)\frac{\Psi_{\geq-9}}z\right)
       (2I+x\mathcal A_r)\wt\Psi_{\geq-7}\\
&\quad-(1-x)
 \left((2I+x\mathcal A_r)\frac{\Psi_{\geq-9}}z\right)
       \mathcal A_z\wt\Psi_{\geq-7},\\
b_{5,r}^{\rm ext}/r={}&\tfrac12\wt\Omega_{\geq-9}
 -(1-x)\frac{\Omega_{\geq-11}}z\,
                 \mathcal A_z\wt\Psi_{\geq-7},\\
b_{5,z}^{\rm ext}/z={}&\tfrac12\wt\Omega_{\geq-9}
 +\frac{\Omega_{\geq-11}}z\,
                 (2I+x\mathcal A_r)\wt\Psi_{\geq-7}.
\end{aligned}
\end{equation}

\subsection{Coercivity coefficients and selected functionals}\label{app:coercivity-arithmetic}

\subsubsection{Derivatives of the prescribed potential}\label{app:coercivity-derivatives}

For $\Psi_5=z\Psi_{app}^{(c)}(u,v)$, the derivatives entering \eqref{mos2} follow from \eqref{u:cov} by the chain rule; for example,
\[
r\partial_r\partial_z\Psi_5=\frac{6r^2}{(1+u)^3}\Big(\partial_u\Psi_{app}^{(c)}+\frac{6z^2}{(1+v)^3}\partial_u\partial_v\Psi_{app}^{(c)}\Big).
\]
Since $\Delta\Psi_5=z\Luv\Psi_{app}^{(c)}$, the third derivatives enter only through
\[
e_0=\Luv\Psi_{app}^{(c)}+\frac{6z^2}{(1+v)^3}\partial_v\big(\Luv\Psi_{app}^{(c)}\big),
\]
where $\partial_v\Luv$ is obtained by differentiating the coefficients of \eqref{u:Luv} in $v$.

\subsubsection{Boxes and pointwise matrix bounds}\label{app:coercivity-boxes}

We construct the boxes \eqref{c:boxes} from lower bounds for
$e_0$ obtained by \eqref{cap:model-range}. Let
$0=r_0<\cdots<r_{K_b}=200$, $K_b=2482$, be the radial grid, and
$0=z_0<\cdots<z_J=200$ the vertical grid. Their transformed coordinates
contain the spline knots; the radial partition is bisected until it
has $2482$ columns. Outward terminal evaluation cells cover the exact
physical endpoints $200$. For $1\leq k\leq K_b$ and $1\leq j\leq J$,
let $\ell_{k,j}$ be a lower bound satisfying
\[
\ell_{k,j}\leq\inf_{[r_{k-1},r_k]\times[z_{j-1},z_j]}e_0.
\]

The levels $L_{p,k}$ and heights $h_{p,k}$ are chosen from these lower bounds so that \eqref{c:boxes} holds; the inequalities are checked on every cell after rounding.

The partition includes spline knots, box jumps and support boundaries.
At a grid boundary, each adjacent one-sided box piece is enclosed
on its closure.

Near $r=0$, let $r_+$ bound the radial coordinate on the cell.
With $\gamma_0$ as in \eqref{c:parameters}, we retain the positive
contribution of the singular terms using
\[
\frac{3(1-\gamma_0)}{r^2}\geq\frac{3(1-\gamma_0)}{r_+^2},
\qquad \frac{0.035}{r^2}\geq\frac{0.035}{r_+^2}
\quad(0<r\leq r_+).
\]

\subsubsection{Exterior coefficient bounds}\label{app:coercivity-exterior}

For an odd exterior potential $zF(x,t,\ell)$ on $r,z\geq0$, in the variables of \eqref{cap:asymptotic-polynomial},
Proposition~\ref{p:ops} and \eqref{mos2} give its coefficient fields as
\begin{equation}\label{cap:coercivity-exterior}
\begin{aligned}
e_{rr}={}&-x\big(\mathcal A_rF
 +(1-x)(\mathcal A_r\mathcal A_zF-2\mathcal A_zF)\big)
 -2\big(F+(1-x)\mathcal A_zF\big),\\
e_{zz}={}&x\big(\mathcal A_rF
 +(1-x)(\mathcal A_r\mathcal A_zF-2\mathcal A_zF)\big)
 +F+(1-x)\mathcal A_zF,\\
a_{12}={}&\frac{\sqrt{x(1-x)}}2
 \big(4\mathcal A_rF+x(\mathcal A_r^2F-2\mathcal A_rF)
 -3\mathcal A_zF-(1-x)(\mathcal A_z^2F-2\mathcal A_zF)\big),\\
e_0={}&\mathcal L_1F+(1-x)\mathcal A_z(\mathcal L_1F).
\end{aligned}
\end{equation}
The four expressions in \eqref{cap:coercivity-exterior} are linear
in $F$. Evaluating them at $F=\Psi_{\geq-9}/z-\sigma t^2$ therefore
gives, respectively, $e_{rr}-e_{rr}^{(1)}$, $e_{zz}-e_{zz}^{(1)}$,
$a_{12}-a_{12}^{(1)}$ and $e_0-e_0^{(1)}$, with the leading
coefficients of \eqref{c:lead1}. We subtract these leading terms
before taking absolute values.

The coefficient differences $e_{rr}-e_{rr}^{(1)}$,
$e_{zz}-e_{zz}^{(1)}$, $e_0-e_0^{(1)}$, and
$(a_{12}-a_{12}^{(1)})/\sqrt{x(1-x)}$ are finite polynomials.
For the last expression the quotient denotes the polynomial obtained
by removing the displayed angular factor in
\eqref{cap:coercivity-exterior}. Since $\sqrt{x(1-x)}\leq1/2$,
\eqref{cap:exterior-supremum}, applied after multiplication by
$t^{-2}$ or $t^{-4}$, gives the following tests for \eqref{c:Hyp}:
\begin{equation}\label{cap:coercivity-exterior-comparisons}
\begin{aligned}
\sup_{R\geq L_e}R^2|e_{rr}-e_{rr}^{(1)}|
&\leq\mathfrak S_{0,0}\big(t^{-2}(e_{rr}-e_{rr}^{(1)});L_e\big)\leq0.05\sigma,\\
\sup_{R\geq L_e}R^2|e_{zz}-e_{zz}^{(1)}|
&\leq\mathfrak S_{0,0}\big(t^{-2}(e_{zz}-e_{zz}^{(1)});L_e\big)\leq0.05\sigma,\\
\sup_{R\geq L_e}R^2|a_{12}-a_{12}^{(1)}|
&\leq\frac12\mathfrak S_{0,0}\left(
 t^{-2}\frac{a_{12}-a_{12}^{(1)}}{\sqrt{x(1-x)}};L_e\right)\leq0.05\sigma,\\
\sup_{R\geq L_e}R^4|e_0-e_0^{(1)}|
&\leq\mathfrak S_{0,0}\big(t^{-4}(e_0-e_0^{(1)});L_e\big)\leq\sigma.
\end{aligned}
\end{equation}

\subsubsection{Taylor models on the plateaux}\label{cap:appendix-mode-arithmetic}

Let $g$ be the plateau piece of the numerical approximation
$\tilde\psi_m$ in Proposition~\ref{prop : spectral data}, with
$\bar\lambda=\bar\lambda_m>0$. With prescribed normalized and
rounded amplitudes $\alpha,\beta$, its values at a radial plateau
point $c>0$ are
\begin{equation}\label{cap:plateau-initial-values}
g(c)=\frac{\alpha J_1(\sqrt{\bar\lambda} c)
              +\beta Y_1(\sqrt{\bar\lambda} c)}c,
\qquad
g'(c)=\frac{\sqrt{\bar\lambda}}c
 \big(\alpha J_0(\sqrt{\bar\lambda} c)+\beta Y_0(\sqrt{\bar\lambda} c)\big)
 -\frac{2g(c)}c.
\end{equation}

For an axial plateau the corresponding values are
$g(c)=\alpha\cos(\sqrt{\bar\lambda} c)+\beta\sin(\sqrt{\bar\lambda} c)$ and
$g'(c)=\sqrt{\bar\lambda}[-\alpha\sin(\sqrt{\bar\lambda} c)+\beta\cos(\sqrt{\bar\lambda} c)]$.
The radial values require only Bessel functions of orders zero and
one, which we enclose by the following finite sums. Write $\gamma_E$
for Euler's constant and
$\mathsf H_k=\sum_{j=1}^k j^{-1}$, with $\mathsf H_0=0$.

\begin{lemma}\label{cap:bessel-seed-lemma}
For $n\in\{0,1\}$ and $\xi>0$, put
\[
\tau_{n,k}=\frac{(-1)^k(\xi/2)^{2k+n}}{k!(k+n)!},\qquad
\eta_{n,k}=\mathsf H_k+\mathsf H_{k+n}-2\gamma_E.
\]
The Bessel series \cite[\S\S10.2,10.8]{DLMF} give
\begin{equation}\label{cap:bessel-seed-series}
J_n(\xi)=\sum_{k\geq0}\tau_{n,k},\qquad
Y_n(\xi)=-\frac{2\mathbf1_{n=1}}{\pi\xi}
+\frac2\pi\log(\xi/2)J_n(\xi)
-\frac1\pi\sum_{k\geq0}\eta_{n,k}\tau_{n,k}.
\end{equation}
For $K\in\mathbb N_0$, define the finite sums
\[
J_n^{[K]}(\xi):=\sum_{k=0}^K\tau_{n,k},\qquad
S_n^{[K]}(\xi):=\sum_{k=0}^K\eta_{n,k}\tau_{n,k}.
\]
If
\[
q=\frac{\xi^2}{4(K+2)(K+n+2)}<1,
\]
their tails satisfy
\begin{equation}\label{cap:bessel-seed-tails}
\begin{aligned}
|J_n-J_n^{[K]}|&\leq T_J:=\frac{|\tau_{n,K+1}|}{1-q},\\
\left|\sum_{k>K}\eta_{n,k}\tau_{n,k}\right|
&\leq T_Y:=|\tau_{n,K+1}|
\left(\frac{|\eta_{n,K+1}|}{1-q}
+\frac{2q}{(K+2)(1-q)^2}\right).
\end{aligned}
\end{equation}
The corresponding approximation
\[
Y_n^{[K]}(\xi):=-\frac{2\mathbf1_{n=1}}{\pi\xi}
+\frac2\pi\log(\xi/2)J_n^{[K]}(\xi)-\frac1\pi S_n^{[K]}(\xi)
\]
satisfies
\begin{equation}\label{cap:bessel-Y-tail}
|Y_n-Y_n^{[K]}|\leq\frac{2|\log(\xi/2)|T_J+T_Y}{\pi}.
\end{equation}
\end{lemma}

\begin{proof}
The ratio of successive absolute terms after the first omitted term
is at most $q$. Also
$|\eta_{n,K+1+m}|\leq|\eta_{n,K+1}|+2m/(K+2)$.
Summing $q^m$ and $mq^m$ proves \eqref{cap:bessel-seed-tails};
the assertion for $Y_n$ follows from \eqref{cap:bessel-seed-series}.
\end{proof}

Evaluating these sums and bounds outward gives the initial values
in \eqref{cap:plateau-initial-values}. The other coefficients follow
from the radial equation $g''+3g'/r+\bar\lambda g=0$. Substituting
$g(r)=\sum_{n\geq0}a_n(r-c)^n$, $c>0$, and comparing coefficients
after multiplication by $r$ gives, for $n\geq0$,
\begin{equation}\label{cap:mode-recurrence}
c(n+2)(n+1)a_{n+2}
=-(n+1)(n+3)a_{n+1}-\bar\lambda c a_n-\bar\lambda a_{n-1},
\qquad a_{-1}=0.
\end{equation}

For axial plateaux the equation is $g''+\bar\lambda g=0$, and the
recurrence becomes $(n+2)(n+1)a_{n+2}=-\bar\lambda a_n$.
These recurrences provide the polynomial parts of the local models.
The next estimate controls what is discarded from each derivative.

\begin{lemma}\label{cap:local-mode-lemma}
Let $\bar\lambda\geq0$, $c>0$, and let $g$ solve
\[
g''(r)+\frac3r g'(r)+\bar\lambda g(r)=0
\]
in a neighborhood of $c$. Suppose that $N\in\mathbb N_0$ and
$0<h<\rho<c$.
Put $\kappa=\max\{1,\sqrt{\bar\lambda}\}$ and
\begin{equation}\label{cap:mode-disc}
B=\max\{|g(c)|,|g'(c)|/\kappa\}
\exp\left(\rho\max\left\{\kappa,
\frac{\bar\lambda}{\kappa}+\frac3{c-\rho}\right\}\right),
\qquad q=h/\rho.
\end{equation}
Write $g(r)=\sum_{n\geq0}a_n(r-c)^n$. The degree-$N$ Taylor
polynomials of $g,g',g''$ satisfy
\begin{equation}\label{cap:mode-tail}
\begin{aligned}
\sup_{|r-c|\leq h}\left|g(r)-\sum_{n=0}^Na_n(r-c)^n\right|
&\leq\frac{Bq^{N+1}}{1-q},\\
\sup_{|r-c|\leq h}\left|g'(r)-\sum_{n=0}^N(n+1)a_{n+1}(r-c)^n\right|
&\leq\frac{B(N+2)q^{N+1}}{\rho(1-q)^2},\\
\sup_{|r-c|\leq h}\left|g''(r)-\sum_{n=0}^N(n+2)(n+1)a_{n+2}(r-c)^n\right|
&\leq\frac{B(N+3)(N+2)q^{N+1}(1+q)}{\rho^2(1-q)^3}.
\end{aligned}
\end{equation}

If instead $g$ solves $g''+\bar\lambda g=0$ near a real center
$c$, the same conclusion holds for $0<h<\rho$, with no restriction
$\rho<c$ and with $3/(c-\rho)$ omitted from $B$.
\end{lemma}

\begin{proof}
The radial system has analytic coefficients on $|\zeta-c|\leq\rho$,
so its solution extends holomorphically to that disc. The axial
system has entire coefficients.
Apply Gronwall's inequality along complex segments to the system
for $(g,g'/\kappa)$, with the maximum norm. Its coefficient norm is
bounded by the maximum in \eqref{cap:mode-disc}. Hence
$|g(\zeta)|\leq B$ on $|\zeta-c|\leq\rho$, and Cauchy's estimate gives
$|a_n|\leq B\rho^{-n}$. For $d=0,1,2$,
\[
\sup_{|r-c|\leq h}\left|g^{(d)}(r)
-\sum_{n=0}^N\frac{(n+d)!}{n!}a_{n+d}(r-c)^n\right|
\leq\frac B{\rho^d}\sum_{n>N}\frac{(n+d)!}{n!}q^n.
\]
Summing these tails, with $n+1\leq(N+2)(n-N)$ and
$(n+2)(n+1)\leq(N+3)(N+2)(n-N)^2$ for $n>N$,
proves \eqref{cap:mode-tail}.
\end{proof}

For a cell contained in $[c-h,c+h]$, we choose the disc radius by
\begin{equation}\label{cap:mode-radius-choice}
\rho=\min\{c/2,(N+1)/\kappa\}\quad\hbox{(radial)},\qquad
\rho=\max\{2h,(N+1)/\kappa\}\quad\hbox{(axial)},
\end{equation}
with downward rounding, and use the model only when $h<\rho$.
The radial choice also keeps the disc away from the singular point
$r=0$ of the plateau equation.

\subsubsection{Polynomial primitive of the radial cutoff}\label{cap:appendix-primitive-arithmetic}

For $T=P[a,n_1,n_2,b]$, the primitive $H_T(r)=\int_0^rt^3T(t)\,dt$ of Lemma~\ref{cap:source-primitive-lemma} is a piecewise polynomial. It vanishes for $r\leq a$, equals $H_T(n_1)+(r^4-n_1^4)/4$ on the plateau $[n_1,n_2]$, and is constant for $r\geq b$. On a ramp, parametrized as $r=c+hy$, $y\in[0,1]$, from its outer endpoint $c\in\{a,b\}$, so that $T(r)=\chi(y)$ with $\chi$ as in~\eqref{c:S}, we have $H_T(r)-H_T(c)=\int_0^y(c+ht)^3\chi(t)h\,dt$; this polynomial is computed exactly from the coefficients of $\chi$.

\subsection{The response defects}\label{cap:appendix-source-response}\label{cap:appendix-response-arithmetic}

\subsubsection{Derivatives on the bounded domain}

The responses \eqref{resp:representation} contain the factors $(1+R^2)^{-s}$, where $1+R^2=1+U_0(u)+U_0(v)$ and $\partial_aU_0(a)=(1+a)^3/3$ for $a\in\{u,v\}$. Their derivatives are computed after multiplication by $(1+R^2)^d$ with $d\geq s+2$, so that all terms are polynomials on each cell; for instance,
\[
(1+R^2)^d\partial_a\frac{\Pi}{(1+R^2)^s}=(1+R^2)^{d-s-1}\Big[(1+R^2)\Pi_a-\frac s3(1+a)^3\Pi\Big].
\]
The common factors are $(1+u)^7(1+v)^7$ and $(1+R^2)^d$,
with $d=6,5$ for $E_j,\wt E_j$, and $d=4,3$ for
$E_j^\psi,\wt E_j^\psi$, respectively. The cutoff terms use
\eqref{cap:gluing-jet-rule}, and the physical factor $z$ is inserted before
subtracting $q_j$ in the odd equation \eqref{resp:residuals}.

\subsubsection{The exterior of the even response}

The function $\wt\Psi_{j,ext}$ in \eqref{resp:representation} is
\begin{equation}\label{resp:multipoles-xt}
\begin{aligned}
\wt\Psi_{j,ext}={}&a_{j,0}t^3+a_{j,1}t^5(4-5x)
+a_{j,2}t^7(8-28x+21x^2)\\
&+a_{j,3}t^9(64-432x+792x^2-429x^3).
\end{aligned}
\end{equation}

By Proposition~\ref{p:ops}, $\Delta\wt\Psi_{j,ext}=0$ for $R>0$.
Using \eqref{cap:asymptotic-derivatives} in \eqref{resp:exterior} gives
\begin{equation}\label{resp:exterior-polynomial}
\begin{aligned}
\wt E_j^{\rm ext}={}&
\left(x\mathcal A_r\wt\Psi_{j,ext}+2\wt\Psi_{j,ext}\right)
\cdot\left[\frac{\Omega_{\geq-11}}z
 +(1-x)\mathcal A_z\!\left(\frac{\Omega_{\geq-11}}z\right)\right]\\
&-x(1-x)(\mathcal A_z\wt\Psi_{j,ext})
 \mathcal A_r\!\left(\frac{\Omega_{\geq-11}}z\right)
 -\wt\mu_j\wt\Omega_{\geq-9}.
\end{aligned}
\end{equation}

\end{document}